\pdfoutput=0

\documentclass[12pt,twoside]{amsart}
\usepackage{amssymb,amsmath}
\usepackage{enumerate}

\newtheorem{thm}{Theorem}[section]
\newtheorem{cor}[thm]{Corollary}
\newtheorem{lem}[thm]{Lemma}

\theoremstyle{definition}

\newtheorem{rem}[thm]{Remark}

\numberwithin{equation}{section}  

\def\vx#1{{\scriptstyle\scriptstyle\scriptstyle
\scriptstyle\circ}\!\!_{_{_{#1}}}\!\!}
\def\vn{{\scriptstyle\scriptstyle\scriptstyle
\scriptstyle\circ}\!\!\!_{_{_{n}}}\!\!}
\def\vl#1{{\scriptstyle\scriptstyle\scriptstyle
\scriptstyle\circ}\!\!\!\!\!_{_{_{#1}}}\!\!\!\!\!}
\def\vo#1{{\scriptstyle\scriptstyle\scriptstyle
\scriptstyle\circ}{#1}\!\!}

\def\edg{-\!\!\!-}
\def\dedg{\,{\scriptstyle{=\!\!\!=\!\!\!=\!\!=}}}
\def\tedg{\,{\scriptstyle{\equiv\!\!\!\equiv\!\!\!\equiv\!\!\equiv}}}
\def\qedg{\,\,^{_=}_{^=}\!\!^{_=}_{^=}\!\!^{_=}_{^=}}
\def\vedgv#1#2{\!^{^{\big|\!\!\!^{^{^{^{^{\vo{#1}\!\!\!^{#2}}}}}}}}\!\!}

\def\redg{}

\def\um{{\phantom{\vx{}-}}}
\def\dm{{\phantom{\vx{a}\edg}}}
\def\tm{{\phantom{\vx{a}\edg\vx{}-}}}
\def\qm{{\phantom{\vx{a}\edg\vx{a}\edg}}}
\def\cm{{\phantom{\vx{a}\edg\vx{a}\edg\vx{}-}}}
\def\sm{{\phantom{\vx{a}\edg\vx{a}\edg\vx{a}\edg}}}
\def\stm{{\phantom{\vx{a}\edg\vx{a}\edg\vx{a}\edg\vx{}-}}}
\def\om{{\phantom{\vx{a}\edg\vx{a}\edg\vx{a}\edg\vx{a}\edg}}}
\def\vm{{\phantom{\vx{a}\edg\vx{a}\edg\vx{a}\edg\vx{a}\edg\vx{}-}}}
\def\num{{\phantom{..}}}
\def\ndm{{\phantom{2}}}
\def\ntm{{\phantom{2..}}}
\def\nqm{{\phantom{22}}}
\def\ncm{{\phantom{22..}}}
\def\nsm{{\phantom{222}}}
\def\nstm{{\phantom{222..}}}

\def\nnm{{\phantom{2222..}}}
\def\ndcm{{\phantom{22222}}}

\def\a{\alpha}
\def\d{\delta}

\def\tn{\tilde n}

\makeatletter
\DeclareRobustCommand*\cal{\@fontswitch\relax\mathcal}
\makeatother

\def\Cal{\cal}
\def\roman{\rm}
\input xy
\xyoption{all}

\def\compo{\,{\scriptstyle\circ}\,}

\def\N{{\Bbb{N}}}
\def\Z{{\Bbb{Z}}}

\def\R{{\Bbb{R}}}
\def\C{{\Bbb{C}}}
\def\sy{{\Cal{S}}}

\def\U{{\Cal U}}
\def\Udr{{\Cal U}^{Dr}}

\def\I{{\Cal I}}

\def\ti{t_i}
\def\uno{{1\!\!\!1\!\!\!1\!\!1\!\!\!1\!\!\!1}}

\def\alg{U}
\def\m{m}
\def\udj{{\Cal U}_q^{DJ}}

\def\ddefi#1{

{{\bf Definition #1.}

}~}
\def\nota#1{

{{\bf Notation #1.}

}~}
\def\lem#1{

{{\bf Lemma #1.}

}~}
\def\prop#1{

{{\bf Proposition #1.}

}~}
\def\teo#1{

{{\bf Theorem #1.}

}~}
\def\cor#1{

{{\bf Corollary #1.}

}~}
\def\dim{

{{\bf Proof:} }}
\def\sqr#1#2{{\vcenter{\vbox{\hrule height .#2pt
                             \hbox{\vrule width .#2pt height#1pt\kern#1pt
                                   \vrule width .#2pt}
                             \hrule height .#2pt}}}}

\def\rem#1{

{{\bf Remark #1.}

}}

\def\1{{\uno}}

\def\intr{\S 0}

\def\prel{\S 1}

\def\prle{\S 2}

\def\drjq{\S 3}

\def\sbal{1}
\def\QNTALG{2}
\def\TTTUQ{3}
\def\grgr{4}
\def\dmpbfn{5}
\def\embd{6}
\def\brgrg{7}

\def\drrl{\S 4}

\def\didi{1}
\def\drcpjdef{2}
\def\ecp{3}
\def\resp{4}
\def\iz{5}
\def\zig{6}
\def\calz{7}
\def\dij{8}
\def\ipd{9}
\def\coefeq{10}
\def\calv{11}
\def\stbcalv{12}
\def\cpbijr{13}
\def\bri{14}
\def\cntb{15}
\def\vcm{16}
\def\rulc{17}
\def\rsfgcal{18}
\def\tilfg{19}
\def\bofg{20}
\def\fgad{21}
\def\txe{22}
\def\scpt{23}
\def\pspol{24}
\def\drdef{25}
\def\greq{26}
\def\prnpr{27}

\def\drmr{\S 5}

\def\aaa{1}
\def\scmp{2}
\def\scpm{3}
\def\nnmm{4}
\def\sercol{5}
\def\sercod{6}
\def\skq{7}
\def\xuit{8}
\def\relser{9}
\def\sercom{10}
\def\sercot{11}
\def\mmnn{12}
\def\crnt{13}
\def\tcalb{14}

\def\drcmp{\S 6}

\def\tilu{1}
\def\bvqz{2}
\def\greh{3}
\def\tiht{4}
\def\genqz{5}
\def\bthi{6}
\def\sepid{7}
\def\idstsep{8}
\def\qgr{9}
\def\lpdgr{10}
\def\tilpm{11}
\def\rtilpm{12}
\def\eqvdd{13}
\def\fuf{14}
\def\gendef{15}
\def\genbdef{16}
\def\disceta{17}

\def\darl{\S 7} 

\def\rkhgr{1}
\def\fibdef{2}
\def\rhxom{3}
\def\xxxtr{4}
\def\rxxom{5}
\def\boto{6}
\def\xdost{7}
\def\ttih{8}
\def\tind{9}
\def\ttibd{10}
\def\thxde{11}
\def\txh{12}
\def\thibd{13}
\def\fth{14}
\def\sop{15}
\def\boh{16}
\def\dnd{17}
\def\solocom{18}
\def\grem{19}

\def\dautemb{\S 8}

\def\zl{1}
\def\omdef{2}
\def\ombdef{3}
\def\xidef{4}
\def\xibdef{5}
\def\tidef{6}
\def\tibdef{7}
\def\aucom{8}

\def\dlem{\S 9}

\def\relid{1}
\def\rgen{2}
\def\bbb{3}
\def\ccc{4}
\def\ddd{5}
\def\eee{6}
\def\fff{7}
\def\ggg{8}
\def\cmh{9}
\def\ctxd{10}
\def\sztmdz{11}
\def\dvssdz{12}
\def\sztm{13}
\def\dvss{14}
\def\ancid{15}
\def\iii{16}
\def\hhh{17}
\def\ctepg{18}
\def\nrelm{19}
\def\stabid{20}

\def\cnat{\S 10}

\def\rku{1}
\def\tsmd{2}
\def\prku{3}
\def\rdxd{4}
\def\cyi{5}
\def\misp{6}
\def\xitd{7}
\def\tldg{8}
\def\hmbdefnt{9}

\def\serrel{\S 11}

\def\mij{1}
\def\cid{2}
\def\tds{3}
\def\tat{4}
\def\sds{5}
\def\tts{6}
\def\sernul{7}
\def\fimbd{8}
\def\finbdef{9}
\def\srsp{10}
\def\ctes{11}
\def\sdkt{12}
\def\sdt{13}
\def\std{14}
\def\cstd{15}
\def\sug{16}
\def\inttr{17}
\def\ssrr{18}
\def\sspz{19}
\def\thsdg{20}

\def\hmfr{\S 12}

\def\defo{1}
\def\defx{2}
\def\defpsi{3}
\def\semi{4}
\def\bdf{5}
\def\scr{6}
\def\bendef{7}
\def\ttst{8}
\def\injfin{9}
\def\ezo{10}
\def\eso{11}

\def\ntt{\S 13}

\def\aka{\cite{1}}
\def\beck{\cite{2}}
\def\bbk{\cite{3}}
\def\cpsld{\cite{4}}
\def\cpu{\cite{5}}
\def\cp{\cite{6}}
\def\cpy{\cite{7}}
\def\damcina{\cite{8}}
\def\drld{\cite{9}} 
 \def\drr{\cite{10}}
\def\fh{\cite{11}}
 
\def\hern{\cite{13}}
\def\iwam{\cite{14}}
\def\jm{\cite{15}}
\def\jing{\cite{16}}
\def\jzt{\cite{17}}
\def\jzh{\cite{18}}
\def\kac{\cite{19}}
\def\levsb{\cite{20}}
\def\lusz{\cite{21}}
\def\matsu{\cite{22}}
\def\nak{\cite{23}}

\begin{document}



\title{DRINFELD REALIZATION OF 
AFFINE QUANTUM
ALGEBRAS: THE RELATIONS.}


\author{Ilaria Damiani}


\maketitle

{\hskip 5 truecm
{\it{To women.}}}

{\hskip 5 truecm
{\it{Especially to those who do not have}}} 

{\hskip 5 truecm
{\it{even the opportunity to imagine}}} 

{\hskip 5 truecm
{\it{how much they would like mathematics}}}

{\hskip 5 truecm
{\it{and to those who are forced to forget it.}}}

\vskip .5 truecm

\begin{abstract}
In this paper the structure of the Drinfeld realization $\Udr_q$ of affine quantum algebras (both untwisted and twisted) is described in details, and its defining relations are studied and simplified. As an application, a homomorphism $\psi$ from this realization to the Drinfeld and Jimbo presentation $\U_q^{DJ}$ is provided, and proved to be surjective.
\end{abstract}

\vskip .5 truecm
{\bf{\intr.\ {\bf INTRODUCTION.}}}
\vskip .5truecm

Let $X_{\tilde n}^{(k)}$ be a Dynkin diagram of affine type, $\U_q^{DJ}=\U_q^{DJ}(X_{\tilde n}^{(k)})$ the quantum algebra introduced by Drinfeld and Jimbo (see \drr\ and \jm), $\Udr_q=\Udr_q(X_{\tilde n}^{(k)})$ its Drinfeld realization (see \drld).

This paper has two main goals: describing in details the structure of the Drinfeld realization $\Udr_q$ arriving at sharply simplifying its defining relations; and constructing a (surjective) homomorphism $\psi$ from this realization to the Drinfeld and Jimbo presentation $\U_q^{DJ}$, as a step to provide a complete proof that $\U_q^{DJ}$ and $\Udr_q$  are isomorphic, so that they are indeed different presentations of the same $\C(q)$-algebra $\U_q=\U_q(X_{\tilde n}^{(k)})$ (see \drld).
\vskip .3truecm
Understanding the isomorphism between $\U_q^{DJ}$ and $\Udr_q$ stated by Drinfeld in \drld\ has important applications in the study of the representation theory of affine quantum algebras: using this result the finite dimensional irreducible representations of affine quantum algebras are classified in \cpsld, \cpu\ and \cp; a geometrical realization (through the quiver varieties) of finite dimensional representations is constructed in \nak\ for the untwisted simply laced cases. 

The interest of the twisted case resides not only in that it is a generalization of the untwisted frame. Actually twisted algebras appear quite naturally while studying the untwisted setting, due to the fact that the transposition of matrices establishes a duality among the affine Cartan matrices through which untwisted Cartan matrices can correspond to twisted ones; more precisely simply laced untwisted 
matrices and the matrices of type $A_{2n}^{(2)}$ are self-dual, while transposition operates on the remaining affine Cartan matrices by interchanging untwisted and twisted ones. This observation is important and concrete because of results like those in \cpy, where the quantum symmetry group of the affine Toda field theory associated to an untwisted affine Kac-Moody algebra is proved to be the quantum algebra associated to the dual Kac-Moody algebra; and in \fh, where the authors conjecture in general, and prove for the Kirillov-Reshetikhin modules, that there exists 
a duality between the representations of an untwisted affine quantum algebra and those of the dual quantum algebra.

\vskip .3truecm
Much work has already been done in 
the direction of understanding Drinfeld's theorem: in \beck\ all the relations are proved in the untwisted case. 
Notice that this does not yet imply that $\psi$ is an isomorphism: indeed the argument for the injectivity should be completed with the proof of the existence of a basis of the integer form, necessary to conclude that the injectivity at 1 implies the injectivity at $q$-level; this point is not discussed and as far as I understand non-trivial.

For the twisted case there are several partial results: in \aka\ the author studies the case $A_2^{(2)}$, constructing  $\psi$ following \beck, but the proof that it is well defined is incomplete; a contribution to this proof is given in \hern. 

In \jing, \jzh\ and \jzt, the authors construct a homomorphism from $\U_q^{DJ}(X_{\tilde n}^{(k)})$ to $\Udr_q(X_{\tilde n}^{(k)})$ (the inverse of $\psi$) following the theorem stated by Drinfeld in \drld, that is by means of the $q$-commutators. 
In \jing\ the author gives some details in the untwisted case, sketching the proof of relations $[E_0,F_i]=0$ ($i\in I_0$) in case $A_3^{(1)}$, of the Serre relation $E_0E_1^2-(q+q^{-1})E_1E_0E_1+E_1^2E_0=0$ in case $A_n^{(1)}$ (noticing that the Serre relations involving just indices in $I_0$ are trivial, but the other Serre relations involving $E_0$ are not studied, for instance $E_1E_0^2-(q+q^{-1})E_0E_1E_0+E_0^2E_1=0$ is missing) and of relations $[E_0,F_0]={K_0-K_0^{-1}\over q_0-q_0^{-1}}$ in cases $A_n^{(1)}$ and $C_2^{(1)}$; but a strategy for generalizing these arguments is not presented, and the twisted case is just stated to be similar.
In \jzh\ the authors concentrate on the twisted case, but their work is again incomplete since the Serre relations involving indices $i\neq j\in I_0$ are treated, erroneously, as in the untwisted case, and for the other relations the authors present some examples: the commutation between $E_0$ and $F_i$ ($i\in I$) is studied in cases $A_2^{(2)}$ and $D_4^{(3)}$; some Serre relations (but not all of them) involving $E_0$ are studied in the cases $A_{2n-1}^{(2)}$ and $D_4^{(3)}$; and again a strategy for generalizing these computations is not shown. It should be noticed a mistake in the connection between the data of the finite Dynkin diagram and its non trivial automorphism on one hand and the twisted affine Dynkin diagram on the other hand, which has consequences in the following paper \jzt.
Finally in \jzt\ the authors want to fill the gap about the Serre relations involving the indices $i,j\in I_0$ such that $a_{ij}<-1$ (in the twisted case), and they use a case by case approach: but the Drinfeld  relations are misunderstood, and stated to imply relations not holding in this algebra. 

\vskip .3 truecm

These difficulties suggest the need to pay some care in understanding the Drinfeld realization, which is the aim of the present paper: the definition of the homomorphism $\psi$ from the Drinfeld realization to the Drinfeld and Jimbo presentation of affine quantum algebras becomes then a simple consequence  of this analysis, and that $\psi$ is surjective is also proved.

\vskip .3 truecm
In sections \prel\ and \prle\ we recall the notions of Dynkin diagram, Weyl group and root system, and their properties needed in developing the arguments of the following sections: in particular it is recalled how untwisted and twisted affine Dynkin diagrams, Weyl groups and root systems are connected to finite ones, and their classification and basic properties. 

\vskip .15 truecm

Section \drjq\ is again a section where some preliminary material (about the presentation of Drinfeld and Jimbo of the (affine) quantum algebras $\U_q^{DJ}$) is summarized. 

In definitions \QNTALG\ and \TTTUQ\ and in remark \grgr\ we recall the definition of $\U_q^{DJ}$, its main structures (the $Q$-gradation, the triangular decomposition, the antiautomorphisms $\Omega$ and $\Xi$, the braid group action, the embedding of the finite quantum algebra in the affine one, the root vectors $E_{\alpha}$) and properties (commutation among the (anti)automorphisms, connection between braid group action and root vectors, Poincar\'e-Birkhoff-Witt basis, Levendorskii-Soibelman formula). 

We recall also the embeddings $\varphi_i$ of the rank 1 quantum algebras $\U_q^{DJ}(A_1^{(1)})$ or $\U_q^{DJ}(A_2^{(2)})$ in the general quantum algebra $\U_q^{DJ}(X_{\tilde n}^{(k)})$ and their properties of commutation and injectivity (definition \embd\ and remark \brgrg). They will play a role in the comparison between the Drinfeld realization and the Drinfeld-Jimbo presentation of section \hmfr\ (theorem \hmfr.\bendef).

\vskip .15 truecm
In section \drrl\ we give the definition of the Drinfeld realization of affine quantum algebras (both untwisted and twisted, see \drld), discussing and translating the relations in a more explicit form, easier for the purpose of this paper. Even if it is just a reformulation, it seems useful to give the details, since they are not always clear in the literature.
\vskip .15 truecm
In section \drmr\ some notations are fixed in order to simplify the analysis of the relations. Also some relations are reformulated in terms of $q$-commutators, and some new relations as the Serre relations $(S^{\pm})$ and other similar ($(T2^{\pm})$ and ($(T3^{\pm})$) are introduced, which will play an important role in the following sections (\cnat\ and \serrel).
\vskip .15 truecm

In section \drcmp\ the main structures on ${\Cal U}_q^{Dr}$ are introduced: the $Q$-gradation; the homomorphisms $\phi_i$, underlining the role of the two affine Drinfeld realizations of rank one, $A_1^{(1)}$ and  $A_2^{(2)}$, which embed in any other  Drinfeld realization, each embedding depending on the choice of a vertex of the (``finite part'' of the) Dynkin diagram;
the  antiautomorphism $\Omega$, describing the correspondence between ``positive'' and  ``negative'' vectors $X_{i,r}^{\pm}$; the automorphisms $\Theta$ and $t_i$ (for each $i\in I_0$), which summarize  several symmetries (reflection around zero and translations) among the ``positive'' vectors. Actually these structures are defined on the algebra $\bar{\Cal U}_q^{Dr}$ (which is also defined in this section), of which the Drinfeld realization is a quotient, and the proof that they induce analogous structures on ${\Cal U}_q^{Dr}$ is quickly concluded in section \dautemb, through the discussion of section \darl.

\vskip .15 truecm

In section \darl\ the algebra $\tilde{\Cal U}_q^{Dr}$, which is an algebra (already introduced in the previous section) intermediate between $\bar{\Cal U}_q^{Dr}$ and ${\Cal U}_q^{Dr}$, is studied in details. In particular a first set of relations is simplified: the most important remarks are that the relations $(HX^{\pm})$ can be replaced by the much easier $(HXL^{\pm})$, see proposition \sop\ (they are much easier not only because they are a smaller set of relations, but mainly because they can be expressed just in terms of $q$-commutation of the generators $X_{i,r}^{\pm}$'s of $\bar{\Cal U}_q^{Dr}$, without using the $H_{i,r}$'s, see remark \solocom); and that the relations  $(HH)$ are also redundant, see proposition \boh. But also the other relations are studied and interpreted while discussing how the structures on $\bar{\Cal U}_q^{Dr}$ (see section \drcmp) induce analogous structures on $\tilde{\Cal U}_q^{Dr}$, see remarks \xdost\ and \tind.

\vskip .15 truecm
Section \dautemb\ is a short and simple section where the structures defined on $\bar{\Cal U}_q^{Dr}$ and induced on $\tilde{\Cal U}_q^{Dr}$ are proved to pass also to ${\Cal U}_q^{Dr}$; this simple analysis is carried out explicitly,  
fixing some notations, in order to use it in further considerations, especially in section \dlem.

\vskip .15 truecm

In section \dlem\ it is now possible to start concentrating on the simplification of the relations defining $\Udr_q$ over $\tilde{\Cal U}^{Dr}_q$: these are the relations involving just the $X_{i,r}^{+}$'s or just the $X_{i,r}^{-}$'s, and there is a correspondence between the two cases thanks to the action of $\tilde \Omega$. The main result of this section is that the dependence of these relations on parameters $(r_1,...,r_l)\in\Z^l$ ($l\in\Z_+$) is redundant: we can indeed just restrict to the same relations indexed by $(r,...,r)\in\Z^l$ where $r\in\Z$ (the ``constant parameter'' relations), so that the dependence on $\Z^l$ is reduced to a dependence on an integer $r$ (see lemmas \dvssdz\ and \dvss, proposition \ancid\ and corollary \nrelm); on the other hand, thanks to the action of the $\tilde t_i$'s, this situation can be again simplified just analyzing the relations relative to $(0,...,0)$ (see remark \ggg).

\vskip .15 truecm

Thanks to the results of section \dlem\ the study of the relations defining $\Udr_q$ can be pushed forward: in section \cnat\ further dependences among the relations are proved (propositions \rku\ and
\rdxd, corollary \misp\ and remark \xitd). These results are summarized in theorem \tldg\ and in corollary \hmbdefnt, where a ``minimal'' set of relations is provided.

\vskip .15 truecm

The last step of this analysis is the study of the Serre relations, performed in section \serrel: here the relations $(XD^{\pm})$-$(S3^{\pm})$ are proved to depend, in the case of rank bigger that 1, on the (``constant parameter'') Serre relations, and these are viceversa proved to depend on the relations $(XD^{\pm})$-$(S3^{\pm})$ also in the cases in which it is not tautologically evident ($k>1$, $a_{ij}<-1$). Theorem \ssrr\ and corollary \sspz\ state the final result of this study, and are the main tool for constructing the homomorphism $\psi$ and for proving that it is well defined, see section \hmfr.

\vskip .15 truecm

Section \hmfr\ is devoted to construct a homomorphism $\psi$ from $\Udr_q$ to $\U_q^{DJ}$ and to prove that it is well defined and surjective. 

In definition \defpsi\ $\tilde\psi:\tilde {\Cal U}^{Dr}_q\to\udj$ is defined, following \beck. It just requires some care in the determination of the sign $o$ (notation \defo\ and remark \defx).

The results of section \serrel\ and the correspondence, described in proposition \semi, between the (anti)automorphisms constructed on $\Udr_q$ and those already known on $\U_q^{DJ}$ make the goal of proving that $\tilde\psi$ induces $\psi$ on  $\Udr_q$ trivial in the cases of rank bigger than 1, that is in all cases different from $A_1^{(1)}$ and $A_2^{(2)}$ (theorem \bdf).

We give two different arguments to solve the cases of rank one (theorem \bendef). The first one is based on the direct computation of the simple commutation relation between $E_1$ and $E_{\delta+\alpha_1}$ in $\U_q^{DJ}(A_1^{(1)})$ and $\U_q^{DJ}(A_2^{(2)})$ (lemma \scr). The second one is a straightforward corollary of the result in the case of rank bigger than 1, once one recalls the embeddings (see remark \drjq.\brgrg) of the rank~1 quantum algebras in the general quantum algebras.

A proof that $\psi$ is surjective is provided in theorem \eso: it makes use of the correspondence between the automorphisms $t_i$ on $\Udr_q$ and the automorphisms $T_{\lambda_i}$ on $\U_q^{DJ}$ and among the $\Omega$'s (remark \ttst), and of the braid group action on $\U_q^{DJ}$. 

Theorem 
\eso\ would suggest also how to define the inverse of $\psi$.

\vskip .15 truecm

An index of the notations used in the paper is listed in the appendix, section~\ntt.

\vskip .4truecm
I'm deeply thankful to David Hernandez for proposing me to work again on the twisted affine quantum algebras: I abandoned them too many years ago, and would have neither planned nor dared to approach them again if he had not encouraged and motivated me. 

I take this occasion to make explicit my gratitude to Corrado De Concini, my 
{\it {maestro}}: for his always caring presence (even when he did not approve my choices) in the vicissitudes of my relationship with mathematics, 
and for his belief (undeserved yet helpful) he made me always feel. Not accidentally, the idea of this work was born at a conference for his $60^{th}$ birthday.

To Andrea Maffei I owe much: because we have been sharing reflections and projects about mathematics and our work since we were students till our adult life; because he is a rare, precious intellectual; and because he is (and has been in this occasion) ready to listen to and help with big and small problems, how specific they can be.
But I owe him even more: his always personal points of view and his friendship.

Eleonora Ciriza is for me more than a colleague, than a mathematician, than an unreplaceable friend: she is all this together. Her support and advice are deep-rooted in a way of being in the world that opened my mind and my life beyond the borders of my own experience.

I do not know if I would have ever arrived at the end of this paper without Salvatore, who had the difficult role of indicating me the purpose of concluding this work as a priority. In particular during the drawing up of the paper, he had to fight hard against my resistance to cut the myriad of other ``priorities'' which absorb much of my concentration and time, and against my delaying attitude of, as Penelope, always undoing what I have done. I thank him for believing in the importance of my work in my and our life.

\vskip .5truecm

{\bf{\prel.\ {\bf PRELIMINARIES: DYNKIN DIAGRAMS.}}}

\vskip .5 truecm
For the preliminary material in this section see \bbk\ and \kac.

A  Dynkin diagram $\Gamma$ of finite or affine type is the datum $(I,A)$ of its set of indices $I$ and its Cartan matrix $A=(a_{ij})_{i,j\in I}\in{\Cal{M}}_{n\times n}(\Z)$ with the following properties: 

i) $a_{ii}=2$ $\forall i\in I$; 

ii) $a_{ij}\leq 0$ $\forall i\neq j\in I$; 

iii) $a_{ij}=0\Leftrightarrow a_{ji}=0$; 

iv) the determinants of all the proper principal minors of $A$ are positive, and $det(A)\geq 0$ ($\Gamma$ is of finite type if $det(A)> 0$ and of affine type if $det(A)= 0$);

$\Gamma$ is said to be indecomposable if furthermore: 

v) if $I=I^{\prime}\cup I^{\prime\prime}$ with $I^{\prime}\cap I^{\prime\prime}=\emptyset$ and $I^{\prime},\ I^{\prime\prime}\neq\emptyset$ then $\exists i^{\prime} \in I^{\prime}$, $i^{\prime\prime}\in I^{\prime\prime}$ such that $a_{i^{\prime} i^{\prime\prime}}\neq 0$.

Between the vertices $i\neq j\in I$ there are $max\{|a_{ij}|,|a_{ji}|\}$ edges, with 
an arrow pointing at $i$ if $|a_{ij}|>|a_{ji}|$; vertices, edges and arrows uniquely determine $\Gamma$.

A Dynkin diagram automorphism of $\Gamma$ is a map $\chi\!:\! I\to\! I$ such that $a_{\chi(i)\chi(j)}=~a_{ij}$ $\forall i,j\in I$. 

It is universally known that  these data are classified (see \bbk); the type of the indecomposable finite data is denoted by $X_{\# I}$ ($X=A,B,C,D,E,F,G$).

\vskip .3 truecm

In this preliminary section we recall the construction and classification of the indecomposable Dynkin diagrams of affine type due to Kac (see \kac) and fix the general notations used in the paper.

Let $\tilde\Gamma$ be an indecomposable Dynkin diagram of finite type, with set of vertices $\tilde I$ 
($\#\tilde I=\tilde n$) and Cartan matrix  $\tilde A=(\tilde a_{i^{\prime}j^{\prime}})_{i^{\prime},j^{\prime}\in\tilde I}$. To $X_{\tilde n}$ it is attached: 

a) the root lattice $\tilde Q=\oplus_{i^{\prime}\in\tilde I}\Z\tilde\alpha_{i^{\prime}}$;

b) the Weyl group $\tilde W\subseteq Aut(\tilde Q)$ generated by the reflections $\{\tilde s_{i^{\prime}}|i^{\prime}\in\tilde I\}$ where $\tilde s_{i^{\prime}}$ is defined by $\tilde s_{i^{\prime}}(\tilde\alpha_{j^{\prime}})=\tilde\alpha_{j^{\prime}}-\tilde a_{i^{\prime}j^{\prime}}\tilde\alpha_{i^{\prime}}$ ($i^{\prime},j^{\prime}\in\tilde I$); 

c) the (uniquely determined up to a scalar factor) $\tilde W$-invariant bilineáar form $(\cdot|\cdot)$ on $\tilde Q$, which induces a positive definite scalar product on $\R\otimes_{\Z}\tilde Q=\oplus_{i^{\prime}\in\tilde I}\R\tilde\alpha_{i^{\prime}}$;

d) the root system $\tilde\Phi\subseteq\tilde Q$, which is the $\tilde W$-orbit of the set $\{\tilde\alpha_{i^{\prime}}|i^{\prime}\in\tilde I\}$ and is also characterized by the property $\tilde\Phi=\{\tilde\alpha\in\tilde Q|\exists i^{\prime}\in\tilde I$ such that $(\tilde\alpha|\tilde\alpha)=(\tilde\alpha_{i^{\prime}}|\tilde\alpha_{i^{\prime}})\}$.

A Dynkin diagram automorphism $\chi$ induces an orthogonal transformation $\chi$ of  $(\tilde Q,(\cdot|\cdot))$ ($\chi(\tilde\alpha_{i^{\prime}})=\tilde\alpha_{\chi(i^{\prime})}$), and we have that $\chi\compo\tilde s_{i^{\prime}}=\tilde s_{\chi(i^{\prime})}\compo\chi$, $\chi(\tilde\Phi)=\tilde\Phi$.

\vskip .3 truecm

Consider the datum $(X_{\tilde n},\chi)$, with $\chi$ Dynkin diagram automorphism of $X_{\tilde n}$ and let $k$ be the order of $\chi$.
It is well known (see \kac) that to this datum it is possible to attach an indecomposable Dynkin diagram of affine type $\Gamma$ and an indecomposable subdiagram of finite type $\Gamma_0\hookrightarrow\Gamma$ with the following properties: 

I) the set of vertices $I$ of $\Gamma$ and $I_0$ of $\Gamma_0$ are $I_0=\tilde I/\chi$ (the set of $\chi$-orbits in $\tilde I$: $\forall i^{\prime}\in \tilde I$ denote by $\bar {i^{\prime}}\in I_0$ the $\chi$-orbit of $i^{\prime}$) and $I=I_0\cup\{0\}$; we shall denote by $n$ the cardinality of $I_0$ and by $\{1,...,n\}$ the set $I_0$ (so that $I=\{0,1,...,n\}$); 

II) the Cartan matrix $A_0$ of $\Gamma_0$ is connected to $\tilde A$ through the relation
$$a_{\bar {i^{\prime}}\bar {j^{\prime}}}=2{\sum_{u\in\Z/k\Z}\tilde a_{\chi^u(i^{\prime})j^{\prime}}\over\sum_{u\in\Z/k\Z}\tilde a_{\chi^u(i^{\prime})i^{\prime}}};$$
remark in particular that if $k=1$ we have $I_0=\tilde I$ and $A_0=\tilde A$, hence $\Gamma_0=\tilde\Gamma$; 

III) the root lattice $Q_0=\oplus_{i\in I_0}\Z\alpha_i$ of $\Gamma_0$ naturally embeds in the root lattice $Q=~\oplus_{i\in I}\Z\alpha_i$ of $\Gamma$; their positive subsets are $Q_{0,+}=\sum_{i\in I_0}\N\alpha_i$ and $Q_{+}=\sum_{i\in I}\N\alpha_i$;

IV) the highest root $\vartheta_0$ of $\Gamma_0$ is characterized by the properties that 
$\vartheta_0\in\Phi_0$ (the root system of $\Gamma_0$) and  $\vartheta_0-\alpha\in Q_{0,+}$ $\forall\alpha\in\Phi_0$; it has also the property that $(\vartheta_0|\vartheta_0)\geq(\alpha|\alpha)$ $\forall\alpha\in\Phi_0$;

V) the highest shortest root $\vartheta_0^{(s)}$ of $\Gamma_0$ is characterized by the properties that 
$\vartheta_0^{(s)}\in\Phi_0$,   $(\vartheta_0^{(s)}|\vartheta_0^{(s)})\leq(\alpha|\alpha)$ $\forall\alpha\in\Phi_0$ and $\vartheta_0^{(s)}-\alpha\in Q_{0,+}$ $\forall\alpha\in\Phi_0$ such that $(\alpha|\alpha)=(\vartheta_0^{(s)}|\vartheta_0^{(s)})$; 

VI) the Cartan matrix $A$ of $\Gamma$ extends $A_0$: $A=(a_{ij})_{i,j\in I}$, with
$$a_{00}=2\ {\roman{and}},\ \forall i\in I_0,\ a_{0i}=-2{(\theta|\alpha_i)\over(\theta|\theta)},\ a_{i0}=-2{(\alpha_i|\theta)\over(\alpha_i|\alpha_i)},$$ where $\theta=\begin{cases} \vartheta_0&{\roman{if}}\ k=1\cr 2\vartheta_0^{(s)}&{\roman{if}}\ X_{\tilde n}=A_{2n}\ {\roman{and}}\ \chi\neq id\cr \vartheta_0^{(s)}&{\roman{otherwise.}}\end{cases} $
\vskip .3 truecm

The type of the Dynkin diagram $\Gamma$ thus constructed is denoted by $X_{\tilde n}^{(k)}$ (indeed it does not depend on $\chi$ but just on $k$), and it is well known (see \kac) that this construction provides a classification of the indecomposable affine Dynkin diagrams, that we list in the following table. 

The
labels under the vertices fix an identification between $I$
and $\{0,1,...,n\}$ such that $I_0$ corresponds to
$\{1,...,n\}$. For each type we also recall the coefficients
$r_i$ (for $i\in I_0$)
in the expression
$\theta=\sum_{i\in I_0}r_i\a_i$ (remark that we correct here a missprint in \damcina: the coefficient $r_n$ for case $A_{2n-1}^{(2)}$).

$$\sm(\Gamma,I)\vm\ \ \ 
(r_1,...,r_n)\ \leqno{X_{\tn=\tn(n)}^{(k)}\
\ \ \  n}$$
${\underline{{\phantom{\hskip 13.5 truecm}}}}$
$$\sm\vx{0}\dedg\vx{1}\om\ \ \ \ \ \ 
\nsm(1)\nstm
\leqno{\ \ A_1^{(1)}\ \ \ \ \ \ \ \ 1}
$$
$$\tm\vx{1}\,^{^{\diagup\!\!^{^{\diagup\!\!^{^{\diagup
}}}}}}
\!\!\!\!\!\!\!\!\!\!\!\edg\vx{2}\,.\,.\!
\!^{^{
{\phantom{\Big|}}
\!\!\!^{^{^{^{^{\vx{0}\!\!\!^{}}}}}}}}\!\!
\,\,.\vl{n-1}\edg\redg
\!\!\!^{^{\diagdown\!\!\!\!\!\!\!\!^{^{\diagdown
\!\!\!\!\!\!\!\!^{^{\diagdown}}}}}}
\,\,\,\vn
\cm\ \ \ \ 
\ndcm(1,...,1)\leqno{\ \
A_n^{(1)}\ \ \ \
\ >1}
$$
$$\tm\vx{1}\!<\!\!\!\dedg\vx{2}\edg\vx{3}\,...\vl{n-2}\edg\vl{n-1}
\vedgv{0}{}\edg\vn\qm\ \ \ \ \ 
\ndm\!(2,...,2,1)\ncm
\leqno{\ \ B_n^{(1)}\ \ \ \ \ 
>2}
$$
$$\dm\vx{1}\dedg\!\!>\!\!\vx{2}\edg\vx{3}\,...\vl{n-1}\edg
\vn\!<\!\!\!\dedg\vx{0}\qm\ \ \ \ \ \ 
(1,2,...,2)\ncm
\leqno{\ \ C_n^{(1)}\ \ \ \ \ 
>1}
$$
$$\dm\vx{2}\edg\vx{3}\vedgv{1}{}
\edg\vx{4}\,...\vl{n-2}
\edg\vl{n-1}\vedgv{0}{}\edg\vn\tm\ \ \ \ \ \ \ \ 
(1,1,2,...,2,1)\num
\leqno{\ \ D_n^{(1)}\ \ \ \ \ 
>3}
$$
$$\tm\vx{2}\edg\vx{3}\edg\vx{4}
\vedgv{1}{\vedgv{0}{}}\!\edg\vx{5}\edg\vx{6}\tm\
\ \ \ \ \ \ntm(2,1,2,3,2,1)\num
\leqno{\ \ E_6^{(1)}\ \ \ \ \ \ \ \ 
6}
$$
$$\um\vx{0}\edg\vx{2}\edg\vx{3}\edg\vx{4}\vedgv{1}{}
\edg\vx{5}\edg\vx{6}\edg\vx{7}\um\
\ \ \ \ \ \num(2,2,3,4,3,2,1)
\leqno{\ \ E_7^{(1)}\ \ \ \ \ \ \ \ 
7}
$$
$$\vx{2}\edg\vx{3}\edg\vx{4}\vedgv{1}{}\edg\vx{5}\edg\vx{6}
\edg\vx{7}\edg\vx{8}\edg\vx{0}\ \
\ \ \ \ (3,2,4,6,5,4,3,2)
\leqno{\ \ E_8^{(1)}\ \ \ \ \ \ \ \ 
8}
$$
$$\um\vx{1}\edg\vx{2}\!<\!\!\!\dedg
\vx{3}\edg\vx{4}\edg\vx{0}\qm\ \ \ \ 
\ncm(2,4,3,2)\ndm
\leqno{\ \ F_4^{(1)}\ \ \ \ \ \ \ \ 
4}
$$
$$\tm\vx{1}\!<\!\!\!\tedg\vx{2}\edg\vx{0}\stm\
\ \ \ \nsm(3,2)\nsm
\leqno{\ \ G_2^{(1)}\ \ \ \ \ \ \ \ 
2}
$$
$$\qm
\vx{1}\!<\!\!\!\qedg\vx{0}\om\ \ \ 
\nstm(2)\nstm
\leqno{\ \ A_2^{(2)}\ \ \ \ \ \ \ \ 
1}
$$
$$\um\vx{1}\!<\!\!\!\dedg\vx{2}\edg\vx{3}\,...
\vl{n-1}
\edg\vn\!<\!\!\!\dedg\vx{0}\qm\ \ \ \ 
\ndm(2,...,2)\nstm
\leqno{\ \ A_{2n}^{(2)}\ \ \ \ \ 
>1}
$$
$$\dm\vx{1}\dedg\!\!>\!\!\vx{2}\edg\vx{3}\,...\vl{n-2}\edg\vl{n-1}
\vedgv{0}{}\edg\vn\tm\ \ \ \ \ 
\ntm(1,2,...,2,1)\nnm
\leqno{\ \ A_{2n-1}^{(2)}\ \ \ 
>2}
$$
$$\um\vx{1}\!<\!\!\!
\dedg\vx{2}\edg\vx{3}\,...\vl{n-1}
\edg\vn\dedg\!\!>\!\!\vx{0}\tm\ \ \ \ 
\ncm(1,...,1)\nstm
\leqno{\ \ D_{n+1}^{(2)}\ \ \ \ 
>1}
$$
$$\dm\vx{0}\edg\vx{1}\edg
\vx{2}\!<\!\!\!\dedg\vx{3}\edg\vx{4}\qm\
\ \ \ 
\nqm(2,3,2,1)\ncm
\leqno{\ \ E_6^{(2)}\ \ \ \ \ \ \ \ 
4}
$$
$$\tm\vx{0}\edg\vx{1}\!<\!\!\!\tedg\vx{2}\stm\
\ \ \ \nsm(2,1)\nsm
\leqno{\ \ D_4^{(3)}\ \ \ \ \ \ \ \ 
2}
$$
\vskip .5truecm

{\bf{\prle.\ {\bf PRELIMINARIES: WEYL GROUP and ROOT SYSTEM.}}}

\vskip .5 truecm

The following structures of the affine Weyl group and root system (see 
\bbk, \iwam, \kac, \matsu) will be used in the paper: 

i) the Weyl group $W_0=<s_i|i\in I_0>\subseteq Aut(Q_0)$ of $\Gamma_0$ acts on $Q$ by
$s_i(\alpha_j)=\alpha_j-a_{ij}\alpha_i$ $\forall i\in I_0,j\in I$ and this action extends to the Weyl group $W=<s_i|i\in I>\subseteq Aut(Q)$ of $\Gamma$ by $s_0(\alpha_i)=\alpha_i-a_{0i}\alpha_0$ $\forall i\in I$; 

ii) the $W$-invariant bilinear form $(\cdot|\cdot)$ on $Q$ 
induces a positive semidefinite symmetric bilinear form on $\R\otimes_{\Z}Q$: it is obviously positive definite on $\R\otimes_{\Z}Q_0$, and has kernel generated by $\delta=\alpha_0+\theta=\sum_{i\in I}r_i\alpha_i\in Q$ where $r_0=1$ always;

iii)  $(\cdot|\cdot)$ 
can be uniquely normalized in such a way that there is a diagonal matrix $D=diag(d_i|i\in I)$ with $1\in\{d_i|i\in I_0\}\subseteq\{d_i|i\in I\}\subseteq\Z_+$ and $(\alpha_i|\alpha_j)=d_ia_{ij}$ $\forall i,j\in I$; $\forall i\in I$, $w\in W$ set $d_{w(\alpha_i)}=d_i$;

iv) $\forall i\in I_0$ define $\tilde d_i=\begin{cases} 
1&{\roman{if}}\ k=1\ {\roman{or}}\ X_{\tilde n}^{(k)}=A_{2n}^{(2)}\cr d_i&{\roman{otherwise;}}\end{cases} $

v) the weight lattice $\hat P\subseteq\R\otimes_{Z}Q_0$ is $\hat P=\oplus_{i\in I_0}\Z\lambda_i$, where $\forall i\in I_0$ $\lambda_i\in\R\otimes_{\Z}Q_0$ is defined by $(\lambda_i|\alpha_j)=\tilde d_i\delta_{ij}$ $\forall j\in I_0$; $Q_0$ naturally embeds in $\hat P$, which provides a $W$-invariant action on $Q$ by 
$x(\alpha)=\alpha-(x|\alpha)\delta$ $\forall x\in \hat P,\ \alpha\in Q$;

vi) as subgroups of $Aut(Q)$ we have $W\leq \hat P\rtimes W_0$; $\hat W=\hat P\rtimes W_0$ is called the extended Weyl group of $\Gamma$ and we have also $\hat W=W\rtimes{\Cal{T}}$, where 
${\Cal{T}}=Aut(\Gamma)\cap\hat  W$;

vii) the extended braid group $\hat{\Cal{B}}$ is the group generated by $\{T_w|w\in\hat W\}$ with relations $T_wT_{w^{\prime}}=T_{ww^{\prime}}$ whenever $l(ww^{\prime})=l(w)l(w^{\prime})$, where $l:\hat W\to\N$ is defined by
$$l(w)=min\{r\in\N|\exists i_1,...,i_r\in I,\tau\in{\Cal{T}}\ {\roman{such\ that}}\ w=s_{i_1}\cdot...\cdot s_{i_r}\tau\};$$
set $T_i=T_{s_i}$ $\forall i\in I$; recall that $l(\sum_{i\in I_0}m_i\lambda_i)=\sum_{i\in I_0}m_il(\lambda_i)$ if $m_i\in \N$ $\forall i\in I_0$;

viii) the root system $\Phi$ of $\Gamma$ decomposes into the union of the sets $\Phi^{{\roman{re}}}$ of the real roots and $\Phi^{{\roman{im}}}$ of the imaginary roots, where $\Phi^{{\roman{re}}}$ is the $W$-orbit in $Q$ of the set $\{\alpha_i|i\in I\}$ and $\Phi^{{\roman{im}}}=\{m\delta|m\in\Z\setminus\{0\}\}$; the set of positive roots is $\Phi_+=\Phi\cap Q_+$;

ix) the multiplicity of the root $\alpha\in\Phi$ is 1 if $\alpha$ is real and $\#\{i\in I_0|\tilde d_i|m\}$ if $\alpha=m\delta$ ($m\in\Z\setminus\{0\}$); the set $\hat{\Phi}$ of roots with multiplicity is $\hat{\Phi}=\Phi^{{\roman{re}}}\cup\hat{\Phi}^{{\roman{im}}}$ where $\hat{\Phi}^{{\roman{im}}}=\{(m\delta,i)|i\in I_0,m\in\Z\setminus\{0\}\},\tilde d_i|m\}$; the set of positive roots with multiplicities is $\hat{\Phi}_+=\Phi^{{\roman{re}}}_+\cup\hat{\Phi}^{{\roman{im}}}_+=(\Phi_+\cap\Phi^{\roman{re}})\cup\{(m\delta,i)\in\hat{\Phi}|m>0\}$;

x) choose a sequence $\iota:\Z\ni r\mapsto\iota_r\in I$ such that $\forall i\in I_0$
$s_{\iota_1}\cdot... \cdot s_{\iota_{_{N_i}}}\tau_i=\sum_{j=1}^i \lambda_j$ and $\forall r\in Z$ $\iota_{r+N_n}=\tau_n(\iota_r)$,  where $N_i=\sum_{j=1}^i l(\lambda_j)$ and $\tau_i\in{\Cal T}$; then $\iota$ induces a map
$$\Z\ni r\mapsto w_r\in W\ \ {\roman{defined\ by}}\ \ w_r=\begin{cases}  s_{\iota_1}\cdot... \cdot s_{\iota_{r-1}}&{\roman {if}}\ r\geq 1\cr s_{\iota_0}\cdot... \cdot s_{\iota_{r+1}}&{\roman {if}}\ r\leq 0\end{cases} $$
and a bijection
$$\Z\ni r\mapsto\beta_r=w_r(\alpha_{\iota_r})\in\Phi^{{\roman{re}}}_+;$$

xi) the total ordering $\preceq$ of $\hat{\Phi}_+$ defined by
$$\beta_r\preceq\beta_{r-1}\preceq(\tilde m\delta,i)\preceq(m\delta,j)\preceq(m\delta,i)\preceq\beta_{s+1}\preceq\beta_s$$
$$\forall r\leq 0,\ s\geq 1,\ \tilde m>m>0,\ j\leq i\in I_0$$
induces on $\Phi_+$ a convex ordering: if $\alpha=\sum_{r=1}^M\gamma_r$ with $M>1$, $\gamma_1\preceq...\preceq\gamma_M$ and $\alpha,\gamma_r\in\Phi_+$ $\forall r=1,...,M$, then either $\gamma_1\prec\alpha$ or $\gamma_r\in\Phi^{{\roman{im}}}$ $\forall r=1,...,M$.

\vskip .5 truecm

{\bf{\drjq.\  PRELIMINARIES:\ the\ DRINFELD-\!JIMBO\ PRESENTATION\ }}$\U_q$.
\vskip .5 truecm

In this section we recall the definition of the quantum algebra $\U_q$ introduced by Drinfeld and Jimbo (see \drr\ and \jm), and the structures  and results (see \beck, \damcina, \levsb,\lusz) needed in \hmfr. First of all recall some notations.
\vskip .3 truecm

\nota{\sbal}
i) For all $i\in I_0$ we denote by $q_i$ the element $q_i=q^{d_i}\in\C(q)$.

ii) Consider the ring $\Z[x,x^{-1}]$. Then for all $m,r\in\Z$ the elements $[m]_x$, $[m]_x!$ ($m\geq 0$) and ${m\brack r}_x$ ($m\geq r\geq 0$) are defined respectively by $[m]_x={x^m-x^{-m}\over x-x^{-1}}$, $[m]_x!=\prod_{s=1}^m[s]_x$ and ${m\brack r}_x={[m]_x!\over[r]_x![m-r]_x!}$,
which all lie in $\Z[x,x^{-1}]$. 

iii) Consider the field $\C(q)$ and, given $v\in\C(q)\setminus\{0\}$, the natural homomorphism $\Z[x,x^{-1}]\to\C(q)$ determined by the condition $x\mapsto v$; then for all $m,r\in\Z$ the elements $[m]_v$, $[m]_v!$ ($m\geq 0$) and ${m\brack r}_v$ ($m\geq r\geq 0$) denote the images in $\C(q)$ respectively of the elements $[m]_x$, $[m]_x!$ and ${m\brack r}_x$ .

\vskip .3 truecm
\ddefi{\QNTALG}
Let $\Gamma=(I,A)$ be a Dynkin diagram of finite or affine type. 

i) 
The (Drinfeld-Jimbo) quantum algebra of type $\Gamma$ is the $\C(q)$-algebra $\U_q\!=~\!\!\U_q(\Gamma)$ generated by
$$\{E_i,F_i,K_i^{\pm 1}|i\in I\}$$ 
with relations: 
$$K_iK_i^{-1}=1=K_i^{-1}K_i,\ \ K_iK_j=K_jK_i\ \ 
\forall i,j\in I,$$ 
$$K_iE_j=q_i^{a_{ij}}E_jK_i,\ \   K_iF_j=q_i^{-a_{ij}}F_jK_i
\ \ \ 
\forall i,j\in I,$$ 
$$[E_i,F_j]=\delta_{ij}{K_i-K_i^{-1}\over q_i-q_i^{-1}}\ \ 
\forall i,j\in I,$$ 
$$\sum_{u=0}^{1-a_{ij}}{1-a_{ij}\brack u}_{q_i} 
E_i^uE_jE_i^{1-a_{ij}-u}=0\ \ \forall i\neq j\in I,$$ 
$$\sum_{u=0}^{1-a_{ij}}{1-a_{ij}\brack u}_{q_i} 
F_i^uF_jF_i^{1-a_{ij}-u}\ \ \forall i\neq j\in I;$$ 
the last two sets of relations are called the Serre relations. 

If $\Gamma$ is affine of type $X_{\tilde n}^{(k)}$ we also set: 

ii) $\U_q^{DJ}=\U_q^{DJ}(X_{\tilde n}^{(k)})=\U_q(\Gamma)$, to stress the distinction of this affine quantum algebra from its Drinfeld realization;

iii) $\U_q^{fin}=\U_q^{fin}(X_{\tilde n}^{(k)})=\U_q(\Gamma_0)$ (see section \prel,I).
\vskip .3 truecm

\ddefi{\TTTUQ}
Recall that  
$\U_q$ is endowed with the following structures: 

i) the $Q$-gradation $\U_q=\oplus_{\alpha\in Q}\U_{q,\alpha}$ determined by the
conditions:
$$E_i\in\U_{q,\a_i},\ \ F_i\in\U_{q,-\a_i},\ \ K_i^{\pm
1}\in\U_{q,0}\ \ \forall i\in I;\ \ \  
\U_{q,\a}\U_{q,\beta}\subseteq\U_{q,\a+\beta}\ \ \forall\a,\beta\in Q;$$

ii) the triangular decomposition: 
$\U_q\cong\U_q^{-}\otimes\U_q^{0}\otimes\U_q^{+}$, where 
$\U_q^{-}$, $\U_q^{0}$ and $\U_q^{+}$ are the subalgebras of $\U_q$ 
generated respectively by 
$\{E_i|i\in I\}$, $\{K_i^{\pm 1}|i\in I\}$ and 
$\{F_i|i\in I\}$; in particular $$\U_{q,\alpha}\cong\bigoplus_{\beta,\gamma\in Q_+:\atop \gamma-\beta=\alpha}\U_{q,-\beta}^{-}\otimes\U_q^{0}\otimes\U_{q,\gamma}^{+}\ \  \forall \alpha\in Q$$ 
where $\U_{q,\alpha}^{\pm}=\U_{q,\alpha}\cap\U_{q}^{\pm}$ $\forall \alpha\in Q$;

iii) the $\C$-anti-linear anti-involution 
$\Omega:\U_q\rightarrow\U_q$ defined by 
$$\Omega(q)= q^{-1};\ \ \Omega(E_i)=F_i,\ \  
\Omega(F_i)=E_i,\ \ 
\Omega(K_i)= K_i^{-1}\ \forall i\in I;$$ 

iv) the $\C(q)$-linear anti-involution 
$\Xi:\U_q\rightarrow\U_q$ defined
by
$$\Xi(E_i)=E_i,\ \  
\Xi(F_i)= F_i,\ \ 
\Xi(K_i)=K_i^{-1}\ \forall i\in I_;$$ 

v) the braid group action defined by 
$$T_i(K_j)=K_j
K_i^{-a_{ij}}
\ \ \forall i,j\in I,$$
$$T_i(E_i)=-F_iK_i,\ \ T_i(F_i)=-K_i^{-1}E_i \ \ \forall i\in I,$$
$$T_i(E_j)=\sum_{r=0}^{-a_{ij}}(-1)^{r-a_{ij}}q_i^{-r}
E_i^{(-a_{ij}-r)}E_jE_i^{(r)},\ \ T_i(F_j)=\Omega(T_i(E_j)) \ \forall i\neq j\in I$$
where $\forall m\in\N$ $E_i^{(m)}={E_i^m\over[m]_{q_i}!}$;

vi) a natural $Aut(\Gamma)$-action: $\tau(K_i)=K_{\tau(i)}$, $\tau(E_i)=E_{\tau(i)}$, $\tau(F_i)=F_{\tau(i)}$ for all $\tau\in Aut(\Gamma)$, $i\in I_0$; if $\Gamma$ is affine then setting $T_{\tau}=\tau$  extends the braid group action to an extended braid group action;

vii) if $\Gamma\hookrightarrow\Gamma^{\prime}$ is a Dynkin diagram embedding then the $\C$-homomorphism $$\varphi_{\Gamma,\Gamma^{\prime}}:\U_q(\Gamma)\to\U_q(\Gamma^{\prime})$$ is naturally defined by $$q\mapsto q^{min\{d_i^{\prime}|i\in I\}},\ \ \ K_i^{\pm 1}\mapsto K_i^{\pm 1},\ \ \ E_i\mapsto E_i,\ \ \ F_i\mapsto F_i\ \ \ (i\in I);$$
in particular if $\Gamma$ is of affine type  $\varphi=\varphi_{\Gamma_0,\Gamma}:\U_q^{fin}\to\udj$ is a 
$\C(q)$-homomorphism;

viii) positive and negative root vectors $E_{\alpha}\in\U_{q,\alpha}^{DJ,+}$ and $F_{\alpha}=\Omega(E_{\alpha})\in\U_{q,-\alpha}^{DJ,-}$ ($\alpha\in\hat{\Phi}_+$) such that if $\Gamma$ is of affine type 
$E_{\beta_r}=T_{w_r}(E_{\iota_r})$ if $r\geq 1$, $E_{\beta_r}=T_{w_r^{-1}}^{-1}(E_{\iota_r})$ if $r\leq 0$, and 
$E_{(\tilde d_ir\delta,i)}=-E_{\tilde d_ir\delta-\alpha_i}E_i+q_i^{-2}E_iE_{\tilde d_ir\delta-\alpha_i}$ if $r>0$, $i\in I_0$.

\vskip .3 truecm
\rem{\grgr}~
We have that: 

i) $\Omega\Xi=\Xi\Omega$, $\Omega T_i=T_i\Omega$ $\forall i\in I$ and $\Omega\tau=\tau\Omega$ $\forall\tau\in {\Cal T}$; 

ii) $\Xi T_i=T_i^{-1}\Xi$ $\forall i\in I$ and $\Xi\tau=\tau\Xi$ $\forall\tau\in {\Cal T}$; 

moreover if $\Gamma$ is of affine type:

iii) $\varphi$ commutes with $\Omega$, $\Xi$ and $T_i$ ($i\in I_0$).

iv) in cases $A_1^{(1)}$ and $A_2^{(2)}$ $\Xi T_1T_{\lambda_1}=T_{\lambda_1}^{-1}\Xi T_1$ (recall that $T_{\lambda_1}=T_0 T_{\tau}=T_{\tau}T_1$, where $<\tau>=Aut(\Gamma)$, in case $A_1^{(1)}$ and $T_{\lambda_1}=T_0T_1$ in case $A_2^{(2)}$);

v) $T_w(\U_{q,\alpha}^{DJ})=\U_{q,w(\alpha)}^{DJ}$ $\forall w\in\hat{W}$, $\alpha\in Q$; 

vi) $T_w(E_i)\in\U_{q,w(\a_i)}^{DJ,+}$ if $w\in \hat W$
and $i\in I$ are such that $w(\a_i)\in Q_+$ (i.e. $l(ws_i)>l(w))$; 

vii) $E_{m\tilde d_i\delta+\alpha_i}=T_{\lambda_i}^{-m}(E_i)$, $F_{m\tilde d_i\delta+\alpha_i}=T_{\lambda_i}^{-m}(F_i)$ $\forall m\in\N, i\in I_0$; 

viii) $\{K_{\alpha}|\alpha\in Q\}$ is a basis of $\U_q^{DJ,0}$, where $K_{\alpha}=\prod_{i\in I}K_i^{m_i}$ if $\alpha=\sum_{i\in I}m_i\alpha_i\in Q$;

ix) $\{E(\gamma)=E_{\gamma_1}\cdot...\cdot E_{\gamma_M}|M\in\N,\ \gamma=(\gamma_1\preceq...\preceq\gamma_{M}),\gamma_{h} \in\hat{\Phi}_+\forall h=1,...,M\}$
is a basis of $\U_q^{DJ,+}$;

x) $\{E(\gamma)K_{\alpha}\Omega(E(\gamma^{\prime}))|\alpha\in Q,\gamma=(\gamma_1\preceq...\preceq\gamma_{M})\in\hat{\Phi}_+^M, \gamma=(\gamma_1\preceq...\preceq\gamma_{M^\prime})\in\hat{\Phi}_+^{M^\prime},M,M^\prime\in\N\}$
is a basis of $\U_q^{DJ}$, called the PBW-basis;

xi) $\forall\alpha\prec\beta\in\hat\Phi_+$ $E_{\beta}E_{\alpha}-q^{(\alpha|\beta)}E_{\alpha}E_{\beta}$ is a linear combination of 
$\{E(\gamma)|\gamma=(\gamma_1\preceq...\preceq\gamma_{M})\in\hat{\Phi}_+^M,\ M\in\N,\ \alpha\prec\gamma_1\}$
(Levendorskii-Soibelman formula).

\vskip .3 truecm 
\rem{\dmpbfn}~
If $\Gamma$ is affine remark \grgr,ix) implies that $dim{\Cal U}_{q,\alpha}^{DJ,+}=dim{\Cal U}_{q,\alpha}^{fin,+}$ $\forall\alpha\in Q_{0,+}$. In particular $\varphi$ is injective.

\vskip .3 truecm 

\ddefi{\embd}
If $\Gamma$ is affine, for $i\in I_0$ let 
$$\varphi_i:\begin{cases} \U_q^{DJ}(A_{1}^{(1)})\rightarrow\U_q^{DJ}(X_{\tilde n}^{(k)})&
{\roman{if}}\ (X_{\tilde n}^{(k)},i)\neq(A_{2n}^{(2)},1)\cr
\U_q^{DJ}(A_{2}^{(2)})\rightarrow\U_q^{DJ}(X_{\tilde n}^{(k)})&
{\roman{if}}\ (X_{\tilde n}^{(k)},i)=(A_{2n}^{(2)},1)\end{cases} $$ be the $\C$-homomorphisms defined on the generators as follows: 
$$q\mapsto q_i,\ K_1^{\pm 1}\mapsto K_i^{\pm 1},\ E_1\mapsto E_i,\ F_1\mapsto F_i$$
and 
$$K_0\mapsto K_{\tilde d_i\delta-\alpha_i},\ E_0\mapsto E_{\tilde d_i\delta-\alpha_i},\ F_0\mapsto F_{\tilde d_i\delta-\alpha_i}\ \ \ {\roman{if}}\ (X_{\tilde n}^{(k)},i)\neq(A_{2n}^{(2)},1)$$
$$K_0\mapsto K_{\delta-2\alpha_1},\ E_0\mapsto E_{\delta-2\alpha_1},\ F_0\mapsto F_{\delta-2\alpha_1}\ \ \ {\roman{if}}\ (X_{\tilde n}^{(k)},i)=(A_{2n}^{(2)},1).$$

\vskip .3 truecm
\rem{\brgrg}~
i) $\varphi_i\Omega=\Omega\varphi_i$, $\varphi_i T_1=T_i\varphi_i$ and
$\varphi_i T_{\lambda_1}=T_{\lambda_i}\varphi_i $ $\forall i\in I_0$;

ii) $\varphi_i$ ($i\in I_0$) is injective (thanks to the PBW-bases).

\vskip .5 truecm

{\bf{\drrl.\ The DRINFELD REALIZATION $\Udr_q$: DEFINITION.}}
\vskip .5 truecm

In this section the definition of the Drinfeld realization ${\Cal U}_q^{Dr}(X_{\tilde n}^{(k)})$ of the affine quantum algebra of type $X_{\tilde n}^{(k)}$ is presented; the definition is discussed and reformulated using the
set $I_0\times\Z$ as index set for the generators 
instead of the set $\tilde I\times\Z$ used in \drld\ and followed in literature (see for instance \cp, \jing, \jzh), because the relations translated  from $\tilde I\times\Z$ to $I_0\times\Z$ seem simpler to handle, even though they lose the immediate connection with the datum 
$(\tilde I,\chi)$. This reformulation, which is useful if one aims to compare the Drinfeld realization with the Drinfeld-Jimbo presentation, is not difficult, but it is presented with some care in order to avoid any ambiguity.

\vskip .3truecm
\nota{\didi}
i) $\omega$ denotes a primitive $k^{\roman{th}}$ root of 1.

ii) Fix the normalization of the $\tilde W$-invariant bilinear form $(\cdot|\cdot)$ on $\tilde Q$ such that 
$min\{\sum_{u\in\Z/k\Z}(\tilde\alpha_{i^{\prime}}|\tilde\alpha_{\chi^u(i^{\prime})})|i^{\prime}\in\tilde I\}=2$.

iii) Denote by $\tilde d$ the number $\tilde d=max\{\tilde d_i|i\in I_0\}$ (in case $A_{2n}^{(2)} $ $\tilde d=1$, otherwise $\tilde d=k$).

iv) Let $Y$ be a function from $\Z^l$ ($\l\in\N$) to any algebra; given $\sigma\in\sy_l$ and $p=(p_1,...,p_l)\in\Z^l$ set
$\sigma.(Y(p))=Y(\sigma.p)=Y(p_{\sigma^{-1}(1)},...,p_{\sigma^{-1}(l)})$. 

v) Analogously if $f\in\C(q)[[u_1^{\pm 1},...,u_l^{\pm 1}]]$ and $u=(u_1,...,u_l)$ define $\sigma.(f(u))$ by $\sigma.(f(u))=f(u_{\sigma^{-1}(1)},...,u_{\sigma^{-1}(l)})$ for all $\sigma\in\sy_l$.

vi) By ``$(R^{\pm})$ is the relation $S^{\pm}=0$'' it is meant
``$(R^{+})$ is the relation $S^{+}=0$ and $(R^{-})$ is the relation $S^{-}=0$''. 

vii) More generally ``$A^{\pm}$ has the property $P^{\pm}$'' means
``$A^{+}$ has the property $P^{+}$ and $A^{-}$ has the property $P^{-}$''.

\vskip .3truecm

For the definition of the Drinfeld realization of affine quantum algebras, that we recall here, see \drld.
\vskip .3truecm

\ddefi{\drcpjdef}
Let $X_{\tilde n}^{(k)}$ be a Dynkin diagram of affine type; the Drinfeld realization of the quantum algebra of type $X_{\tilde n}^{(k)}$ is the $\C(q)$-algebra $\Udr_q(X_{\tilde n}^{(k)})=\Udr_q$ generated by
$${\Cal C}^{\pm 1},\ \ \ 
{\Cal K}_{i^{\prime}}^{\pm 1}\ \ (i^{\prime}\in\tilde I),\ \ \ {\Cal X}_{i^{\prime},r}^{\pm}\ \ ((i^{\prime},r)\in\tilde I\times\Z),\ \ \ 
{\Cal H}_{i^{\prime},r}\ \ ((i^{\prime},r)\in\tilde I\times(\Z\setminus\{0\})),\leqno{({\Cal{G}})}$$
with the following relations $({\Cal{DR}})$:
$${\Cal K}_{\chi(i^{\prime})}={\Cal K}_{i^{\prime}},\ \ {\Cal H}_{\chi(i^{\prime}),r}=\omega^r {\Cal H}_{i^{\prime},r}\ \ (i^{\prime}\in \tilde I,\ r\in\Z\setminus\{0\})
,\leqno{({\Cal Z})}$$
$${\Cal X}_{\chi(i^{\prime}),r}^{\pm}=\omega^r {\Cal X}_{i^{\prime},r}^{\pm}\ \ ((i^{\prime},r)\in \tilde I\times\Z),\leqno{({\Cal {ZX}}^{\pm})}$$
$${\Cal C}{\Cal C}^{-1}=1,\ \ \ [{\Cal C},x]=0\ \ \forall x,\leqno{({\Cal C})}$$
$${\Cal K}_{i^{\prime}}{\Cal K}_{i^{\prime}}^{-1}=1={\Cal K}_{i^{\prime}}^{-1}{\Cal K}_{i^{\prime}}, \ \ \ {\Cal K}_{i^{\prime}}{\Cal K}_{j^{\prime}}={\Cal K}_{j^{\prime}}{\Cal K}_{i^{\prime}}\ \ (i^{\prime},j^{\prime}\in\tilde I)
,\leqno{({\Cal {KK}})}$$
$${\Cal K}_{i^{\prime}}{\Cal X}_{j^{\prime},r}^{\pm}=
q^{\pm\sum_{u\in \Z/k\Z}(\tilde\alpha_{i^{\prime}}|\tilde\alpha_{\chi^{u}(j^{\prime})})}{\Cal X}_{j^{\prime},r}^{\pm}{\Cal K}_{i^{\prime}}\ \ \ (i^{\prime},j^{\prime}\in\tilde I,\ r\in\Z)
,\leqno({\Cal {KX}}^{\pm})$$
$$[{\Cal K}_{i^{\prime}},{\Cal H}_{j^{\prime},r}]=0\ \ (i^{\prime},j^{\prime}\in\tilde I,\ r\in\Z\setminus\{0\})
,\leqno{({\Cal {KH}})}$$
$$[{\Cal X}_{i^{\prime},r}^+,{\Cal X}_{j^{\prime},s}^-]={\sum_{u=0}^{k-1}\delta_{\chi^u(i^{\prime}),j^{\prime}}\omega^{us}\over\sum_{u=0}^{k-1}\delta_{\chi^u(i^{\prime}),i^{\prime}}}\cdot
{{\Cal C}^{-s}{\Cal K}_{i^{\prime}}\tilde {\Cal H}_{i^{\prime},r+s}^+-{\Cal C}^{-r}{\Cal K}
_{i^{\prime}}^{-1}\tilde {\Cal H}_{i^{\prime},r+s}^-\over (q-q^{-1})[{1\over 2}\sum_{u\in\Z/k\Z}(\tilde\alpha_{\chi^u(i^{\prime})}|\tilde\alpha_{i^{\prime}})]_q}\leqno{({\Cal{XX}})}$$
($(i^{\prime},r),(j^{\prime},s)\in\tilde I\times\Z$), 
$$ [{\Cal H}_{i^{\prime},r},{\Cal X}_{j^{\prime},s}^{\pm}]=\pm \tilde b_{i^{\prime}j^{\prime}r}{\Cal C}^{r\mp|r|\over 2}{\Cal X}_{j^{\prime},r+s}^{\pm}\ \ ((i^{\prime},r)\in\tilde I\times(\Z\setminus\{0\}),(j^{\prime},s)\in \tilde I\times\Z),
\leqno{({\Cal {HX}}^{\pm})}$$
$$ [{\Cal H}_{i^{\prime},r},{\Cal H}_{j^{\prime},s}]=\d_{r+s,0} \tilde b_{i^{\prime}j^{\prime}r}{C^r-C^{-r}\over (q-q^{-1})[{1\over 2}\sum_{u\in\Z/k\Z}(\tilde\alpha_{\chi^u(j^{\prime})}|\tilde\alpha_{j^{\prime}})]_q}
\leqno{({\Cal {HH}})}$$
($(i^{\prime},r),(j^{\prime},s)\in \tilde I\times(\Z\setminus\{0\})$),
$$F_{i^{\prime}j^{\prime}}^{\pm}(u_1,u_2){\Cal X}_{i^{\prime}}^{\pm}(u_1){\Cal X}_{j^{\prime}}^{\pm}(u_2)=G_{i^{\prime}j^{\prime}}^{\pm}(u_1,u_2){\Cal X}_{j^{\prime}}^{\pm}(u_2){\Cal X}_{i^{\prime}}^{\pm}(u_1)\ \ \ (i^{\prime},j^{\prime}\in\tilde I),\leqno{(\Cal{XFG}^{\pm})}$$
$$ \sum_{\sigma\in\sy_3}\sigma.((q^{- 3\varepsilon}u_1^{\pm\varepsilon}-(q+q^{-1})u_2^{\pm\varepsilon}+q^{ 3\varepsilon}u_3^{\pm\varepsilon}){\Cal X}_{i^{\prime}}^{\pm}(u_1){\Cal X}_{i^{\prime}}^{\pm}(u_2){\Cal X}_{i^{\prime}}^{\pm}(u_3))=0\leqno{({\Cal X}3^{\varepsilon,\pm})}$$
($i^{\prime}\in\tilde I$, $\tilde a_{\chi(i^{\prime})i^{\prime}}=-1$),
$$\sum_{\sigma\in\sy_{1-a_{ij}}}\!\!\!\!\!\sigma.\sum_{u=0}^{1-a_{ij}}(-1)^u{1-a_{ij}\brack u}_{q_i}{\Cal X}_{i^{\prime},p_{1}}^{\pm}\cdot...\cdot {\Cal X}_{i^{\prime},p_{u}}^{\pm}{\Cal X}_{j^{\prime},v}^{\pm}{\Cal X}_{i^{\prime},p_{u+1}}^{\pm}\cdot...\cdot {\Cal X}_{i^{\prime},p_{1-a_{ij}}}^{\pm}=0\leqno{({\Cal S}^{\pm})}$$
($k=1$, $i^{\prime},j^{\prime}\in \tilde I$, $i^{\prime}\neq j^{\prime}$),
$$\sum_{\sigma\in\sy_2}\sigma.\Big(P_{i^{\prime}j^{\prime}}^{\pm}(u_1,u_2)\big(
{\Cal X}_{j^{\prime}}^{\pm}(v){\Cal X}_{i^{\prime}}^{\pm}(u_1){\Cal X}_{i^{\prime}}^{\pm}(u_2)+\leqno{({\Cal{XP}}^{\pm})}$$
$$-[2]_{q^{\m_{i^{\prime}j^{\prime}}}}
{\Cal X}_{i^{\prime}}^{\pm}(u_1){\Cal X}_{j^{\prime}}^{\pm}(v){\Cal X}_{i^{\prime}}^{\pm}(u_2)+
{\Cal X}_{i^{\prime}}^{\pm}(u_1){\Cal X}_{i^{\prime}}^{\pm}(u_2){\Cal X}_{j^{\prime}}^{\pm}(v)
\big)\Big)=0$$
($k>1$, $i^{\prime},j^{\prime}\in \tilde I$, $\chi(i^{\prime})\neq j^{\prime}$, $\tilde a_{i^{\prime}j^{\prime}}<0$),

where $\tilde{ \Cal H}_{i^{\prime},r}^{\pm}$, $\tilde b_{i^{\prime}j^{\prime}r}$, ${\Cal X}_{i^{\prime}}^{\pm}(u)$, $F_{i^{\prime}j^{\prime}}^{\pm}(u_1,u_2)$, $G_{i^{\prime}j^{\prime}}^{\pm}(u_1,u_2)$, $\varepsilon$, $P_{i^{\prime}j^{\prime}}^{\pm}(u_1,u_2)$ and $\m_{i^{\prime}j^{\prime}}$
are defined as follows:
$$\sum_{r\in\Z}\tilde{\Cal H}_{i^{\prime},\pm r}^{\pm}u^r=exp\left(\pm(q-q^{-1})\Big[{1\over 2}\sum_{u\in\Z/k\Z}(\tilde\alpha_{\chi^u(i^{\prime})}|\tilde\alpha_{i^{\prime}})\Big]_q
\sum_{r>0}{\Cal H}_{i^{\prime},\pm r}u^r\right);$$
$$\tilde b_{i^{\prime}j^{\prime}r}={\sum_{u=0}^{k-1}[r(\tilde\alpha_{i^{\prime}}|\tilde\alpha_{\chi^u(j^{\prime})})]_q\omega^{ru}\over r[{1\over 2}\sum_{u\in\Z/k\Z}(\tilde\alpha_{i^{\prime}}|\tilde\alpha_{\chi^u(i^{\prime})})]_q};$$
$${\Cal X}_{i^{\prime}}^{\pm}(u)=\sum_{r\in\Z}{\Cal X}_{i^{\prime},r}^{\pm}
u^{-r};$$
$$F_{i^{\prime}j^{\prime}}^{\pm}(u_1,u_2)=\prod_{v\in\Z/k\Z\atop \tilde a_{i^{\prime},\chi^v(j^{\prime})}\neq 0}(u_1-\omega^vq^{\pm(\tilde\alpha_{i^{\prime}}|\tilde\alpha_{\chi^v(j^{\prime})})}u_2);$$
$$G_{i^{\prime}j^{\prime}}^{\pm}(u_1,u_2)=\prod_{v\in\Z/k\Z\atop \tilde a_{i^{\prime},\chi^v(j^{\prime})}\neq 0}(q^{\pm(\tilde\alpha_{i^{\prime}}|\tilde\alpha_{\chi^v(j^{\prime})})}u_1-\omega^vu_2);$$
$$\varepsilon=\pm 1;$$
$$P_{i^{\prime}j^{\prime}}^{\pm}(u_1,u_2)=\begin{cases}  1&{\roman {if}}\ \tilde a_{i^{\prime},\chi(i^{\prime})}=0\ {\roman{and}}\ \chi(j^{\prime})\neq j^{\prime},\ {\roman{or}}\ \chi(i^{\prime})=i^{\prime}\cr
{q^{\pm 2k}u_1^k-u_2^k\over q^{\pm 2}u_1-u_2}&{\roman{otherwise};}\end{cases} $$
$$\m_{i^{\prime}j^{\prime}}=\begin{cases}  {k\over\tilde d}\sum_{u\in\Z/k\Z}\delta_{i^{\prime}\chi^u(i^{\prime})}&{\roman {if}}\ \tilde a_{i^{\prime},\chi(i^{\prime})}=0\ {\roman{and}}\ \chi(j^{\prime})\neq j^{\prime},\ {\roman{or}}\ \chi(i^{\prime})=i^{\prime}\cr
k &{\roman{otherwise}.}\end{cases} $$
\vskip .3 truecm
\rem{\ecp}~
In \drld\ not all the relations $({\Cal X}3^{\varepsilon,\pm})$ appear, but just the relations $({\Cal X}3^{1,+})$ and $({\Cal X}3^{-1,-})$; relations $({\Cal X}3^{-1,+})$ and $({\Cal X}3^{1,-})$ are introduced in \cp\ as consequences of relations $({\Cal Z})-({\Cal {XFG}}^{\pm})$, $({\Cal X}3^{1,+})$, $({\Cal X}3^{-1,-})$, $({\Cal S}^{\pm})$, $({\Cal {XP}}^{\pm})$, since their use simplifies some calculations, making evident some symmetries (the stability of the relations under the antiautomorphism $\tilde \Omega$ and the automorphism $\tilde \Theta$). Here we use relations $({\Cal X}3^{\varepsilon,\pm})$ for the same reasons of simplification (see remarks \dautemb.\ombdef\ and \dautemb.\xibdef), proving in proposition \cnat.\rku\ the equivalence stated in \cp.

\vskip .3truecm

\rem{\resp}~
i) For all $r\in\Z$, the algebra generated by $\{{\Cal{Y}}_{i^{\prime}}|i^{\prime}\in \tilde I\}$
with relations
$\{{\Cal{Y}}_{\chi(i^{\prime})}=\omega^r{\Cal{Y}}_{i^{\prime}}|i^{\prime}\in \tilde I\}$ is isomorphic to the algebra generated by $\{Y_{i}|i\in I_0\}$
with relations
$\{Y_{i}=\omega^{r\#\{i^{\prime}\in\tilde I:\bar{i^{\prime}}=i\}}Y_{i}|i\in I_0\}$, where a section $\tilde{}:I_0\to\tilde I$ induces an isomorphism $Y_i\mapsto{\Cal{Y}}_{\tilde i}$. 

ii) Consider $i^{\prime}\in\tilde I$ and let $i\in I_0=\tilde I/\chi$ be the $\chi$-orbit of $i^{\prime}$. Notice that $
\chi(i^{\prime})=i^{\prime}\Leftrightarrow k|\tilde d_i$; more precisely $\sum_{u\in\Z/k\Z}\delta_{i^{\prime},\chi^u(i^{\prime})}=\tilde d_i$ and $\tilde d_i\#\{i^{\prime}\in\tilde I:\bar{i^{\prime}}=i\}=k$.

iii) For all $r\in\Z$, the algebra generated by $\{Y_{i}|i\in I_0\}$
with relations
$\{Y_{i}=\omega^{{kr\over\tilde d_i}}Y_{i}|i\in I_0\}$ is trivially isomorphic to 
the algebra generated by $\{Y_{i}|i\in I_0\}$
with relations $\{Y_i=0|\tilde d_i\not|r\}$, which is trivially isomorphic to 
the free algebra generated by $\{Y_{i}|i\in I_0$ such that $\tilde d_i|r\}$.

iv) Hence, for all $r\in\Z$, the algebra generated by $\{{\Cal{Y}}_{i^{\prime}}|i^{\prime}\in \tilde I,\}$
with relations
$\{{\Cal{Y}}_{\chi(i^{\prime})}=\omega^r{\Cal{Y}}_{i^{\prime}}|i^{\prime}\in \tilde I\}$ is isomorphic to 
the algebra generated by $\{Y_{i}|i\in I_0\}$
with relations
$\{Y_i=0|\tilde d_i\not|r\}$, where a section $\tilde{}:I_0\to\tilde I$ induces an isomorphism $Y_i\mapsto{\Cal{Y}}_{\tilde i}$. 

v) Finally, the algebra generated by $\{{\Cal{Y}}_{i^{\prime},r}|i^{\prime}\in \tilde I,r\in\Z\}$
with relations
$\{{\Cal{Y}}_{\chi(i^{\prime}),r}=\omega^r{\Cal{Y}}_{i^{\prime},r}|i^{\prime}\in \tilde I,r\in\Z\}$ is isomorphic to 
the algebra generated by $\{Y_{i,r}|i\in I_0,r\in\Z\}$
with relations
$\{Y_{i,r}=0|\tilde d_i\not|r\}$, or equivalently to the free algebra generated by $\{Y_{i,r}|i\in I_0,r\in\Z$ such that $\tilde d_i|r\}$.

\vskip .3 truecm

\nota{\iz}
Let us denote by $I_{\Z}$ the set $I_{\Z}=\{(i,r)\in I_0\times\Z|\tilde d_i|r\}$.

\vskip .3 truecm

\cor{\zig}
i) $\Udr_q$ is (isomorphic to) an algebra generated by 
$$C^{\pm1},\ \ \ k_i^{\pm1}\ \ (i\in I_0),\ \ \ X_{i,r}^{\pm}\ \ ((i,r)\in I_0\times{\Z}),\ \ \ H_{i,r}\ \ ((i,r)\in I_0\times({\Z}\setminus \{0\}));\leqno{({G})}$$
the relations $$X_{i,r}^{\pm}=0\ \ \forall (i,r)\in (I_0\times\Z)\setminus I_{\Z}\leqno{(ZX^{\pm})} $$
and
$$H_{i,r}=0\ \ \forall (i,r)\in (I_0\times\Z)\setminus I_{\Z},\leqno{(ZH)} $$
hold in $\Udr_q$.

ii) $\Udr_q$ is generated by 
$$C^{\pm1},\ \ \ k_i^{\pm1}\ \ (i\in I_0),\ \ \ X_{i,r}^{\pm}\ \ ((i,r)\in I_{\Z},\ \ \ H_{i,r}\ \ ((i,r)\in I_{\Z}\setminus( I_0\times\{0\})).\leqno{({G}^{\prime})}$$

\vskip .3 truecm
\rem{\calz}~
The relations $({\Cal{ZX}}^{\pm})$
are equivalent to the condition ${\Cal X}_{\chi(i^{\prime})}^{\pm}(u)={\Cal X}_{i^{\prime}}^{\pm}(\omega^{-1}u)$ for all $i^{\prime}\in\tilde I$.

\vskip .3 truecm

\nota{\dij}
Given $i,j\in I_0$ we set $\tilde d_{ij}=max\{\tilde d_i,\tilde d_j\}$. 

\vskip .3 truecm

\rem{\ipd}~
i) If $\tilde\alpha_{i^{\prime}}$ ($i^{\prime}\in\tilde I$) is a short root then $(\tilde\alpha_{i^{\prime}}|\tilde\alpha_{i^{\prime}})={2k\over\tilde d}$.

ii) $\forall i^{\prime}\in\tilde I$ we have $\sum_{u\in\Z/k\Z}(\tilde\alpha_{i^{\prime}}|\tilde\alpha_{\chi^u(i^{\prime})})=2d_{\bar{i^{\prime}}}$.

\vskip .3truecm

\rem{\coefeq}~
i) Remark that there exists a section \ $\tilde{}:I_0\to\tilde I$ such that given $i,j\in I_0$ we have that $a_{ij}\neq 0\Rightarrow\tilde a_{\tilde i \tilde j}\neq 0$ (of course it is always true that 
$\tilde a_{i^{\prime}j^{\prime}}\neq 0\Rightarrow a_{\bar i^{\prime}\bar j^{\prime}}\neq 0$); 

ii) let  \ $\tilde{}$ \ be a section as in i); then, if $k>1$, $d_ia_{ij}=max\{d_i,d_j\}\tilde a_{\tilde i \tilde j}$.

\vskip .3truecm

\rem{\calv}~
i) The relations $({\Cal{KK}})$, $({\Cal{KX}}^{\pm})$ and
$({\Cal{KH}})$ are compatible with the relations $({\Cal{Z}})$ and 
$({\Cal{ZX}}^{\pm})$, in the sense that for all $i^{\prime},j^{\prime}\in\tilde I$, $r\in\Z$, $s\in\Z\setminus\{0\}$,
$$({\Cal{KK}})_{\chi(i^{\prime}),j^{\prime}}=({\Cal{KK}})_{i^{\prime},j^{\prime}}=({\Cal{KK}})_{i^{\prime},\chi(j^{\prime})},$$ $$({\Cal{KX}}^{\pm})_{\chi(i^{\prime}),j^{\prime},r}=
({\Cal{KX}}^{\pm})_{i^{\prime},j^{\prime},r},\ \ \ ({\Cal{KX}}^{\pm})_{i^{\prime},\chi(j^{\prime}),r}=\omega^r
({\Cal{KX}}^{\pm})_{i^{\prime},j^{\prime},r}$$ and
$$({\Cal{KH}})_{\chi(i^{\prime}),j^{\prime},s}=
({\Cal{KH}})_{i^{\prime},j^{\prime},s},\ \ \ ({\Cal{KH}})_{i^{\prime},\chi(j^{\prime}),s}=\omega^s
({\Cal{KH}})_{i^{\prime},j^{\prime},s};$$
ii)  if $i,j\in I_0$ are such that $\bar{i^{\prime}}=i$, $\bar{j^{\prime}}=j$, 
$({\Cal{KX}}^{\pm})_{i^{\prime},j^{\prime}}$ is equivalent to 
$${\Cal K}_{i^{\prime}}{\Cal X}_{j^{\prime},r}^{\pm}=
q_i^{\pm a_{ij}}{\Cal X}_{j^{\prime},r}^{\pm}{\Cal K}_{i^{\prime}}$$
(see \prel.II), remark \ipd,ii) and notation \drjq.\sbal,i)).
\vskip 2 truecm

\rem{\stbcalv}~
i) If we apply $\chi$ to the expression $\sum_{r>0}{\Cal H}_{i^{\prime},\pm r}u^r$ $(i^{\prime}\in\tilde I)$ we get (see $({\Cal Z})$)
$$\sum_{r>0}{\Cal H}_{\chi(i^{\prime}),\pm r}u^r=\sum_{r>0}{\Cal H}_{i^{\prime},\pm r}(\omega^{\pm 1}u)^r.$$

ii) From i) and from definition \drcpjdef\ we get 
$$\sum_{r\in\Z}\tilde{\Cal H}_{\chi(i^{\prime}),\pm r}^{\pm}u^r=\sum_{r\in\Z}\tilde{\Cal H}_{i^{\prime},\pm r}^{\pm}(\omega^{\pm 1}u)^r\ \ \ \forall i^{\prime}\in\tilde I,$$
that is
$\tilde{\Cal H}_{\chi(i^{\prime}),r}^{\pm}=\omega^{r}\tilde{\Cal H}_{i^{\prime},r}^{\pm}$ $\forall (i^{\prime},r)\in\tilde I\times(\Z\setminus\{0\})$.

iii) The relations $({\Cal{XX}})$ are compatible with relations $({\Cal Z})$ and $({\Cal {ZX}}^{\pm})$:
$$({\Cal{XX}})_{\chi(i^{\prime}),r;j^{\prime},s}=\omega^r({\Cal {XX}})_{i^{\prime},r;j^{\prime},s} ,\ \ \ 
({\Cal{XX}})_{i^{\prime},r;\chi(j^{\prime}),s}=\omega^s({\Cal {XX}})_{i^{\prime},r;j^{\prime},s}.$$

iv)  If \ $\tilde{}$ \ is as in remark \coefeq\ and  $i,j\in I_0$, then $({\Cal{XX}})_{\tilde i,r;\tilde j,s}$ is equivalent to $$ 
[{\Cal X}_{\tilde i,r}^+,{\Cal X}_{\tilde j,s}^-]=\begin{cases} \delta_{ij}
{{\Cal C}^{-s}{\Cal K}_{\tilde i}\tilde {\Cal H}_{\tilde i,r+s}^+-{\Cal C}^{-r}{\Cal K}
_{\tilde i}^{-1}\tilde {\Cal H}_{\tilde i,r+s}^-\over (q_i-q_i^{-1})}&{\roman{if}}\ 
\tilde d_j|s\cr 0&{\roman{otherwise}}\end{cases} 
$$
(see remark \resp,ii)).

\vskip .3truecm
\rem{\cpbijr}~
Let $i^{\prime},j^{\prime}\in\tilde I$, $r\in\Z$; then: 

i) $\tilde b_{i^{\prime}j^{\prime}r}={1\over r}\sum_{u=0}^{k-1}\left[r(\tilde\alpha_{i^{\prime}}|\tilde\alpha_{\chi^u(j^{\prime})})\over d_{\bar{i^{\prime}}}\right]_{q_{\bar{i^{\prime}}}}\omega^{ru}$;

ii) $\tilde b_{\chi(i^{\prime})j^{\prime}r}=\omega^r \tilde b_{i^{\prime}j^{\prime}r}$ and $\tilde b_{i^{\prime}\chi(j^{\prime})r}=\omega^{-r} \tilde b_{i^{\prime}j^{\prime}r}$;

iii) the relations
$({\Cal {HX}}^{\pm})$
and
$({\Cal{HH}})$
are compatible with the relations $({\Cal Z})$ and $({\Cal{ZX}}^{\pm})$:
$$({\Cal{HX}}^{\pm})_{\chi(i^{\prime}),r;j^{\prime},s}=\omega^r({\Cal{HX}}^{\pm})_{i^{\prime},r;j^{\prime},s},\ \ \ ({\Cal{HX}}^{\pm})_{i^{\prime},r;\chi(j^{\prime}),s}=\omega^s({\Cal{HX}}^{\pm})_{i^{\prime},r;j^{\prime},s}$$
and $$({\Cal{HH}})_{\chi(i^{\prime}),r;j^{\prime},s}=\omega^r({\Cal{HH}})_{i^{\prime},r;j^{\prime},s},\ \ \ ({\Cal{HH}})_{i^{\prime},r;\chi(j^{\prime}),s}=\omega^s({\Cal{HH}})_{i^{\prime},r;j^{\prime},s}.$$

\vskip .3 truecm
\nota{\bri}
Let $i,j\in I_0$, $r\in\Z$; $b_{ijr}$ denotes the element of $\C(q)$
$$b_{ijr}=\begin{cases}  0&{\roman{if}}\ \tilde d_{i,j}\not|r\cr
{[2r]_q(q^{2r}+(-1)^{r-1}+q^{-2r})\over r}&{\roman{if}}\ (X_{\tilde n}^{(k)},i,j)=(A_{2n}^{(2)},1,1)\cr
{[\tilde r a_{ij}]_{q_i}\over\tilde r}&{\roman{otherwise,\ with\ }}
\tilde r={r\over \tilde d_{i,j}}.\end{cases} $$

\vskip .3 truecm
\prop{\cntb}
If \ $\tilde{}$ \ is as in remark \coefeq, then
$\forall i,j\in I_0,\ r\in\Z$ $\tilde b_{\tilde i\tilde j r}=b_{ijr}$.
In particular $({\Cal{HX}})_{\tilde i,r;\tilde j,s}$ and $({\Cal{HH}})_{\tilde i,r;\tilde j,s}$ are equivalent to 
$$ [{\Cal H}_{\tilde i,r},{\Cal X}_{\tilde j,s}^{\pm}]=\pm b_{ijr}{\Cal C}^{r\mp|r|\over 2}{\Cal X}_{\tilde j^,r+s}^{\pm}\ \ \ {\roman {and}}\ \ \ 
 [{\Cal H}_{\tilde i,r},{\Cal H}_{\tilde j,s}]=\d_{r+s,0} b_{ijr}{C^r-C^{-r}\over q_j-q_j^{-1}}.
$$
\dim
If $k=1$ the claim is trivial. Suppose now $k>1$ (so that $X_{\tilde n}$ is simply laced) and notice that (see remark \ipd,i))
$$\tilde b_{\tilde i\tilde jr}={1\over r}\sum_{u=0}^{k-1}\left[rk\tilde a_{\tilde i\chi^u(\tilde j)}\over \tilde dd_i\right]_{q_i}\omega^{ru};$$
moreover if $\tilde d_{ij}=k$ then either $\chi(\tilde i)=\tilde i$ or $\chi(\tilde j)=\tilde j$
so that $\tilde a_{\tilde i\chi^u(\tilde j)}=\tilde a_{{\tilde i}\tilde j}$ for all $u$, $\tilde r={r\over k}$,  $\tilde d_h=d_h$ $\forall h\in I_0$, $\tilde d=k$ and, thanks to remark \coefeq,ii), $$\tilde b_{\tilde i\tilde j r}={1\over r}\left[rk\tilde a_{\tilde i\tilde j}\over \tilde dd_i\right]_{q_i}\sum_{u=0}^{k-1}\omega^{ru}=\begin{cases}  0&{\roman {if}}\ k\not|r\cr{1\over \tilde r}\left[r a_{ i j}\over \tilde d\right]_{q_i}={[\tilde r a_{ i j}]_{q_i}\over \tilde r}&{\roman {if}}\ k|r.\end{cases} $$  
If $(X_{\tilde n}^{(k)},i,j)=(A_{2 n}^{(2)},1,1)$ 
$$\tilde b_{\tilde i\tilde j r}={1\over r}([4r]_q+(-1)^r[-2r]_q)={1\over r}[{2r}]_q\left(q^{2r}+(-1)^{r-1}+q^{-2r}\right).$$ 
In the remaining cases $\tilde a_{\tilde i\chi^u(\tilde j)}=0$ when $k\not|u$, $\tilde r=r$ and $max\{d_i,d_j\}={k\over \tilde d}$, hence 
$$\tilde b_{\tilde i\tilde j r}={1\over r}\left[{kr\tilde a_{\tilde i\tilde j}\over \tilde dd_i}\right]_{q_i}=
{1\over r}[r a_{ij}]_{q_i}={1\over \tilde r}[\tilde r a_{ij}]_{q_i}.$$
\vskip .3 truecm

In the next remarks as well as in all the paper the $q$-commutators play a fundamental role in simplifying the description of the elements and in the computations. We recall here their definition and simple properties (see also \jing).
\vskip .3 truecm

\nota{\vcm}
Given $v\in\C(q)\setminus\{0\}$ and $a,b$ elements of a $\C(q)$-algebra, the element $[a,b]_v$ is defined by $[a,b]_v=ab-vba$.
\vskip .3 truecm

\rem{\rulc}~
Let $a,b,c\in\tilde{\Cal U}_q^{Dr}$ and $u,v,w\in\C(q)\setminus\{0\}$. Then:

i) $[a,b]_u=-u[b,a]_{u^{-1}}$;

ii) $[[a,b]_u,b]_v=[[a,b]_v,b]_u=ab^2-(u+v)bab+uvb^2a$;

iii) $[[a,b]_u,c]_v=[a,[b,c]_{v/w}]_{uw}-u[b,[a,c]_{w}]_{v/uw}$. 

If moreover  $a\in{\Cal U}_{q,\alpha}^{Dr}$, $b\in{\Cal U}_{q,\beta}^{Dr}$ and $i\in I_0$ then: 

iv) $[k_ia,b]_u=k_i[a,b]_{q^{-(\alpha_i|\beta)}u}$; 

v) $[a,k_ib]_u=q^{-(\alpha|\alpha_i)}k_i[a,b]_{q^{(\alpha_i|\alpha)}u}$.

\vskip .3truecm
\rem{\rsfgcal}~
Let $i^{\prime},j^{\prime}\in \tilde I$; then:

i) $F_{i^{\prime}j^{\prime}}^{\pm}(u_1,u_2)$ and $G_{i^{\prime}j^{\prime}}^{\pm}(u_1,u_2)$ are homogeneous polynomials of the same degree~ $d$;

ii) $F_{\chi(i^{\prime})j^{\prime}}^{\pm}(u_1,u_2)
=\omega^dF_{i^{\prime}j^{\prime}}^{\pm}(\omega^{-1} u_1,u_2)$, $F_{i^{\prime}\chi(j^{\prime})}^{\pm}(u_1,u_2)=F_{i^{\prime}j^{\prime}}^{\pm}(u_1,\omega^{-1} u_2)$;

iii) $G_{\chi(i^{\prime})j^{\prime}}^{\pm}(u_1,u_2)
=\omega^dG_{i^{\prime}j^{\prime}}^{\pm}(\omega^{-1} u_1,u_2)$, $G_{i^{\prime}\chi(j^{\prime})}^{\pm}(u_1,u_2)=G_{i^{\prime}j^{\prime}}^{\pm}(u_1,\omega^{-1} u_2)$;

iv) the relations
$({\Cal{XFG}}^{\pm})$
are compatible with the relations $({\Cal{ZX}}^{\pm})$:
$$({\Cal{XFG}}^{\pm})_{\chi(i^{\prime}),j^{\prime}}(u_1,u_2)=\omega^d({\Cal{XFG}}^{\pm})_{i^{\prime},j^{\prime}}(\omega^{-1}u_1,u_2),$$
$$({\Cal{XFG}}^{\pm})_{i^{\prime},\chi(j^{\prime})}(u_1,u_2)=({\Cal{XFG}}^{\pm})_{i^{\prime},j^{\prime}}(u_1,\omega^{-1}u_2).$$

\vskip .3 truecm

\rem{\tilfg}~
Let $i^{\prime},j^{\prime}\in \tilde I$ be such that $\tilde a_{i^{\prime},\chi^r(j^{\prime})}= 0$ for all $r\in\Z$; this is equivalent to the condition $a_{\bar{i^{\prime}}\bar{j^{\prime}}}=0$. Then: 

i) $F_{i^{\prime}j^{\prime}}^{\pm}(u_1,u_2)=G_{i^{\prime}j^{\prime}}^{\pm}(u_1,u_2)=1$; 

ii) the relation $({\Cal{XFG}}^{\pm})_{i^{\prime},j^{\prime}}$
is equivalent to
$[{\Cal X}_{i^{\prime}}^{\pm}(u_1),{\Cal X}_{j^{\prime}}^{\pm}(u_2)]=0$,
that is to $$[{\Cal X}_{i^{\prime},r}^{\pm},{\Cal X}_{j^{\prime},s}^{\pm}]=0\ \ \ \forall r,s\in\Z.$$

\vskip .3 truecm
\rem{\bofg}~
Let $i^{\prime},j^{\prime}\in \tilde I$.

i) The condition 
$(\tilde\alpha_{i^{\prime}}|\tilde\alpha_{\chi^r(j^{\prime})})=(\tilde\alpha_{i^{\prime}}|\tilde\alpha_{j^{\prime}})\neq 0\ \ \ \forall r\in\Z$
is equivalent to the con\-di\-tions $a_{\bar{i^{\prime}}\bar{j^{\prime}}}\neq 0$, $\tilde d_{\bar{i^{\prime}}\bar{j^{\prime}}}=k$ and 
implies that $\tilde d_{\bar{i^{\prime}}\bar{j^{\prime}}}(\tilde\alpha_{i^{\prime}}|\tilde\alpha_{j^{\prime}})=d_{\bar{i^{\prime}}}a_{\bar{i^{\prime}}\bar{j^{\prime}}}$;

ii) the condition $\exists! r\in\Z/k\Z$ such that $(\tilde\alpha_{i^{\prime}}|\tilde\alpha_{\chi^r(j^{\prime})})\neq 0$  is equivalent to the conditions  $a_{\bar{i^{\prime}}\bar{j^{\prime}}}\neq 0$, $\tilde d_{\bar{i^{\prime}}\bar{j^{\prime}}}=1$, $(X_{\tilde n}^{(k)},\bar{i^{\prime}},\bar{j^{\prime}})\neq(A_{2n}^{(2)},1,1)$; this condition  
implies that $(\tilde\alpha_{i^{\prime}}|\tilde\alpha_{\chi^r(j^{\prime})})=d_{\bar{i^{\prime}}}a_{\bar{i^{\prime}}\bar{j^{\prime}}}$.

Let $i,j\in I_0$ and choose \ $\tilde{}$ \ as in remark \coefeq. 

iii) If $a_{ij}\neq 0$ and $(X_{\tilde n}^{(k)},i,j)\neq(A_{2n}^{(2)},1,1)$ (that is, $\tilde i$, $\tilde j$ satisfy the conditions of i) or ii) with $r=0$) then
$$F_{\tilde i\tilde j}^{\pm}(u_1,u_2)=u_1^{\tilde d_{ij}}-q_i^{\pm a_{ij}}u_2^{\tilde d_{ij}},\ \ G_{\tilde i\tilde j}^{\pm}(u_1,u_2)=q_i^{\pm a_{ij}}u_1^{\tilde d_{ij}}-u_2^{\tilde d_{ij}}$$
and the relation $({\Cal{XFG}}^{\pm})_{\tilde i,\tilde j}$
is equivalent to
$$[u_1^{\tilde d_{ij}}{\Cal X}_{\tilde i}^{\pm}(u_1),{\Cal X}_{\tilde j}^{\pm}(u_2)]_{q_i^{\pm a_{ij}}}+[u_2^{\tilde d_{ij}}{\Cal X}_{\tilde j}^{\pm}(u_2),{\Cal X}_{\tilde i}^{\pm}(u_1)]_{q_i^{\pm a_{ij}}}=0,$$
that is to $$[{\Cal X}_{\tilde i,r+\tilde d_{ij}}^{\pm},{\Cal X}_{\tilde j,s}^{\pm}]_{q_i^{\pm a_{ij}}}+[{\Cal X}_{\tilde j,s+\tilde d_{ij}}^{\pm},{\Cal X}_{\tilde i,r}^{\pm}]_{q_i^{\pm a_{ij}}}=0\ \ \ \forall r,s\in\Z.$$
Notice that 
$$[{\Cal X}_{\tilde i,r+\tilde d_{ij}}^{-},{\Cal X}_{\tilde j,s}^{-}]_{q_i^{-a_{ij}}}=-q_i^{-a_{ij}}[{\Cal X}_{\tilde j,(s+\tilde d_{ij})-\tilde d_{ij}}^{-},{\Cal X}_{\tilde i,r+\tilde d_{ij}}^{-}]_{q_i^{a_{ij}}}$$
so that 
$({\Cal{XFG}}^{\pm})_{\tilde i,\tilde j}$
is equivalent to
$$[{\Cal X}_{\tilde i,r\pm\tilde d_{ij}}^{\pm},{\Cal X}_{\tilde j,s}^{\pm}]_{q_i^{a_{ij}}}+[{\Cal X}_{\tilde j,s\pm\tilde d_{ij}}^{\pm},{\Cal X}_{\tilde i,r}^{\pm}]_{q_i^{a_{ij}}}=0\ \ \ \forall r,s\in\Z.$$

\vskip .3 truecm
\rem{\fgad}~
Let $(X_{\tilde n}^{(k)},i,j)=(A_{2n}^{(2)},1,1)$; then if \ $\tilde{}$ \ is as in remark \coefeq:

i) $F_{\tilde i\tilde j}^{\pm}(u_1,u_2)=(u_1-q^{\pm 4}u_2)(u_1+q^{\mp 2}u_2)=u_1^2-(q^{\pm 4}-q^{\mp 2})u_1u_2-q^{\pm 2}u_2^2$, $G_{\tilde i\tilde j}^{\pm}(u_1,u_2)=q^{\pm 2}u_1^2-(q^{\pm 4}-q^{\mp 2})u_1u_2-u_2^2$;

ii) the relation $({\Cal{XFG}}^{\pm})_{\tilde i\tilde j}$ is equivalent to
$$[u_1^2{\Cal X}_{\tilde 1}^{\pm}(u_1),{\Cal X}_{\tilde 1}^{\pm}(u_2)]_{q^{\pm 2}}+[u_2^2{\Cal X}_{\tilde 1}^{\pm}(u_2),{\Cal X}_{\tilde 1}^{\pm}(u_1)]_{q^{\pm 2}}+$$
$$-(q^{\pm 4}-q^{\mp 2})(u_1{\Cal X}_{\tilde 1}^{\pm}(u_1)u_2{\Cal X}_{\tilde 1}^{\pm}(u_2)+u_2{\Cal X}_{\tilde 1}^{\pm}(u_2)u_1{\Cal X}_{\tilde 1}^{\pm}(u_1))=0,$$
that is to $$[{\Cal X}_{\tilde 1,r+2}^{\pm},{\Cal X}_{\tilde 1,s}^{\pm}]_{q^{\pm 2}}-q^{\pm 4}[{\Cal X}_{\tilde 1,r+1}^{\pm},{\Cal X}_{\tilde 1,s+1}^{\pm}]_{q^{\mp 6}}+$$
$$+[{\Cal X}_{\tilde 1,s+2}^{\pm},{\Cal X}_{\tilde 1,r}^{\pm}]_{q^{\pm 2}}-q^{\pm 4}[{\Cal X}_{\tilde 1,s+1}^{\pm},{\Cal X}_{\tilde 1,r+1}^{\pm}]_{q^{\mp 6}}=0\ \ \ \forall r,s\in\Z.$$
As in remark \bofg\ notice that in this case $({\Cal{XFG}}^{\pm})_{\tilde i\tilde j}$
is equivalent to
$$[{\Cal X}_{\tilde 1,r\pm 2}^{\pm},{\Cal X}_{\tilde 1,s}^{\pm}]_{q^{2}}-q^{4}[{\Cal X}_{\tilde 1,r\pm 1}^{\pm},{\Cal X}_{\tilde 1,s\pm 1}^{\pm}]_{q^{-6}}+$$
$$+[{\Cal X}_{\tilde 1,s\pm 2}^{\pm},{\Cal X}_{\tilde 1,r}^{\pm}]_{q^{2}}-q^{4}[{\Cal X}_{\tilde 1,s\pm 1}^{\pm},{\Cal X}_{\tilde 1,r\pm 1}^{\pm}]_{q^{- 6}}=0\ \ \ \forall r,s\in\Z.$$

\vskip .3 truecm

\rem{\txe}~
Let $i^{\prime}\in\tilde I$;

i) the relations $({\Cal X}3^{\varepsilon,\pm})$  are compatible with relations $({\Cal {ZX}}^{\pm})$: $$({\Cal X}3^{\varepsilon,\pm})_{\chi(i^{\prime})}(u_1,u_2,u_3)=({\Cal X}3^{\varepsilon,\pm})_{i^{\prime}}(\omega^{-1}u_1,\omega^{-1}u_2,\omega^{-1}u_3);$$

ii) the condition $\tilde a_{\chi(i^{\prime})i^{\prime}}=-1$ is equivalent to the condition
$(X_{\tilde n}^{(k)},\bar{i^{\prime}})=(A_{2n}^{(k)},1)$; 

iii) the relations $({\Cal X}3^{\varepsilon,\pm})$  are equivalent to
$$ \sum_{\sigma\in\sy_3}\sigma.((q^{- 3\varepsilon}u_1^{\pm\varepsilon}{\Cal X}_{i^{\prime}}^{\pm}(u_1){\Cal X}_{i^{\prime}}^{\pm}(u_2){\Cal X}_{i^{\prime}}^{\pm}(u_3)-q^{-\varepsilon}{\Cal X}_{i^{\prime}}^{\pm}(u_2)u_1^{\pm\varepsilon}{\Cal X}_{i^{\prime}}^{\pm}(u_1){\Cal X}_{i^{\prime}}^{\pm}(u_3)+$$
$$-q^{\varepsilon}{\Cal X}_{i^{\prime}}^{\pm}(u_3)u_1^{\pm\varepsilon}{\Cal X}_{i^{\prime}}^{\pm}(u_1){\Cal X}_{i^{\prime}}^{\pm}(u_2)+q^{3\varepsilon}
{\Cal X}_{i^{\prime}}^{\pm}(u_3){\Cal X}_{i^{\prime}}^{\pm}(u_2)u_1^{\pm\varepsilon}{\Cal X}_{i^{\prime}}^{\pm}(u_1))=0,$$
which is 
$$ q^{-3\varepsilon}\sum_{\sigma\in\sy_3}\sigma.[[u_1^{\pm\varepsilon}{\Cal X}_{i^{\prime}}^{\pm}(u_1),{\Cal X}_{i^{\prime}}^{\pm}(u_2)]_{q^{2\varepsilon}},{\Cal X}_{i^{\prime}}^{\pm}(u_3)]_{q^{4\varepsilon}}=0$$
or equivalently
$$\sum_{\sigma\in\sy_3}\sigma.[[{\Cal X}_{i^{\prime}, r_1\pm\varepsilon}^{\pm},{\Cal X}_{i^{\prime},r_2}^{\pm}]_{q^{2\varepsilon}},{\Cal X}_{i^{\prime},r_3}^{\pm}]_{q^{4\varepsilon}}=0\ \ \ \forall r_1,r_2,r_3\in\Z.$$

\vskip .3 truecm
\rem{\scpt}~
The relations $({\Cal S}^{\pm})$ are compatible with $({\Cal {ZX}}^{\pm})$.

\vskip .3 truecm

\rem{\pspol}~
Let $k>1$, $i^{\prime},j^{\prime}\in \tilde I$ be such that $\chi(i^{\prime})\neq j^{\prime}$, $\tilde a_{i^{\prime}j^{\prime}}<0$ (this is equivalent to the condition $k>1$, $a_{\bar i^{\prime}\bar j^{\prime}}<0$). It is immediate to see that: 

i) $P_{\chi(i^{\prime})j^{\prime}}^{\pm}(u_1,u_2)=P_{i^{\prime}j^{\prime}}^{\pm}(u_1,u_2)=P_{i^{\prime}\chi(j^{\prime})}^{\pm}(u_1,u_2)$;

ii) $\m_{\chi(i^{\prime})j^{\prime}}=\m_{i^{\prime}j^{\prime}}=\m_{i^{\prime}\chi(j^{\prime})}$; 

iii) $P_{i^{\prime}j^{\prime}}^{\pm}(u_1,u_2)$ is homogeneous (of some degree $d$);

iv) the relations $({\Cal{XP}}^{\pm})$ are compatible with $({\Cal{ZX}}^{\pm})$:
$$({\Cal{XP}}^{\pm})_{\chi(i^{\prime}),j^{\prime}}(u_1,u_2;v)=\omega^d({\Cal{XP}}^{\pm})_{i^{\prime},j^{\prime}}(\omega^{-1}u_1,\omega^{-1}u_2;v)$$ and $$({\Cal{XP}}^{\pm})_{\chi(i^{\prime}),j^{\prime}}(u_1,u_2;\omega^{-1}v)=({\Cal{XP}}^{\pm})_{i^{\prime},j^{\prime}}(u_1,u_2;v).$$

Moreover if $i,j\in I_0$ are such that $i=\bar{i^{\prime}}$, $j=\bar{j^{\prime}}$ then: 

v) the condition $\ \tilde a_{i^{\prime},\chi(i^{\prime})}=0\ {\roman{and}}\ \chi(j^{\prime})\neq j^{\prime},\ {\roman{or}}\ \chi(i^{\prime})=i^{\prime}$ is equivalent to the condition $a_{ij}=-1$; 

vi) $\m_{i^{\prime}j^{\prime}}=\begin{cases}  d_i&{\roman{if}}\ a_{ij}=-1\cr
k(=-a_{ij})&{\roman{otherwise}};\end{cases} $

vii) the relation $({\Cal{XP}}^{\pm})_{i^{\prime},j^{\prime}}$ is equivalent to

$$\sum_{\sigma\in\sy_2}\sigma.\sum_{r,s\geq 0\atop r+s=-1-a_{ij}}q^{\pm 2s
}\big({\Cal X}_{j^{\prime}}^{\pm}(v)u_1^s{\Cal X}_{i^{\prime}}^{\pm}(u_1)u_2^r{\Cal X}_{i^{\prime}}^{\pm}(u_2)+$$
$$-[2]_{q^{\m_{i^{\prime}j^{\prime}}}}u_1^s
{\Cal X}_{i^{\prime}}^{\pm}(u_1){\Cal X}_{j^{\prime}}^{\pm}(v)u_2^r{\Cal X}_{i^{\prime}}^{\pm}(u_2)+$$
$$+u_1^s
{\Cal X}_{i^{\prime}}^{\pm}(u_1)u_2^r{\Cal X}_{i^{\prime}}^{\pm}(u_2){\Cal X}_{j^{\prime}}^{\pm}(v)\big)=0$$
that is
$$\sum_{\sigma\in\sy_2}\sigma.\sum_{r,s\geq 0\atop r+s=-1-a_{ij}}
q^{\pm 2s
}\big({\Cal X}_{j^{\prime},v}^{\pm}{\Cal X}_{i^{\prime},p_1+s
}^{\pm}{\Cal X}_{i^{\prime},p_2+r}^{\pm}+$$
$$-[2]_{q^{\m_{i^{\prime}j^{\prime}}}}{\Cal X}_{i^{\prime},p_1+s
}^{\pm}{\Cal X}_{j^{\prime},v}^{\pm}{\Cal X}_{i^{\prime},p_2+r}^{\pm}+{\Cal X}_{i^{\prime},p_1+s}^{\pm}{\Cal X}_{i^{\prime},p_2+r}^{\pm}{\Cal X}_{j^{\prime},v}^{\pm}\big)=0,$$
or equivalently
$$\sum_{\sigma\in\sy_2}\sigma.\sum_{r,s\geq 0\atop r+s=-1-a_{ij}}
q^{2s}\big({\Cal X}_{j^{\prime},v}^{\pm}{\Cal X}_{i^{\prime},p_1\pm s}^{\pm}{\Cal X}_{i^{\prime},p_2\pm r}^{\pm}+$$
$$-[2]_{q^{\m_{i^{\prime}j^{\prime}}}}{\Cal X}_{i^{\prime},p_1\pm s}^{\pm}{\Cal X}_{j^{\prime},v}^{\pm}{\Cal X}_{i^{\prime},p_2\pm r}^{\pm}+{\Cal X}_{i^{\prime},p_1\pm s}^{\pm}{\Cal X}_{i^{\prime},p_2\pm r}^{\pm}{\Cal X}_{j^{\prime},v}^{\pm}\big)=0.$$

\vskip .3 truecm

We are now ready to write down an equivalent definition of $\Udr_q(X_{\tilde n}^{(k)})$, using the generators $(G)$. 

\vskip .3 truecm 
\prop{\drdef}
$\Udr_q(X_{\tilde n}^{(k)})$ is (isomorphic to) the $\C(q)$-algebra generated by 
$$C^{\pm1},\ \ \ k_i^{\pm1}\ \ (i\in I_0),\ \ \ X_{i,r}^{\pm}\ \ ((i,r)\in I_0\times{\Z}),\ \ \ H_{i,r}\ \ ((i,r)\in I_0\times({\Z}\setminus \{0\})),\leqno{(G)}$$
with the following relations $(DR)$:
$$X_{i,r}^{\pm}=0\ \ \forall (i,r)\in (I_0\times\Z)\setminus I_{\Z},\leqno{(ZX^{\pm})} $$
$$H_{i,r}=0\ \ \forall (i,r)\in (I_0\times\Z)\setminus I_{\Z},\leqno{(ZH)} $$
$$[C,x]=0\ \ \forall x,\ \ \ k_ik_j=k_jk_i\ \ \ (i,j\in I_0),\leqno{(CUK)}$$
$$CC^{-1}=1,\ \ \ k_ik_i^{-1}=1=k_i^{-1}k_i\ \ \ (i\in I_0),\leqno{(CK)}$$
$$k_iX_{j,r}^{\pm}=q_i^{\pm a_{ij}}X_{j,r}^{\pm}k_i\ \ (i\in I_0,\ (j,r)\in I_0\times\Z),\leqno{(KX^{\pm})}$$
$$[k_i,H_{j,s}]=0\ \ (i\in I_0,\ (j,s)\in I_0\times(\Z\setminus\{0\})),\leqno{(KH)}$$
$$[X_{i,r}^+,X_{j,s}^-]=\begin{cases} \d_{i,j}
{C^{-s}k_i\tilde H_{i,r+s}^+ -C^{-r}k_i^{-1}\tilde H_{i,r+s}^-\over q_i-q_i^{-1}}&{\roman {if}}\ \tilde  d_j|s\cr 0&{\roman {otherwise}}
\end{cases} 
\leqno{(XX)}$$
$((i,r),(j,s)\in I_0\times\Z)$
$$ [H_{i,r},X_{j,s}^{\pm}]=\pm b_{ijr}C^{r\mp|r|\over 2}X_{j,r+s}^{\pm}\ \ ((i,r)\in I_0\times(\Z\setminus\{0\}),(j,s)\in I_0\times\Z),\leqno{(HX^{\pm})}$$
$$ [H_{i,r},H_{j,s}]=\d_{r+s,0}b_{ijr}{C^r-C^{-r}\over q_j-q_j^{-1}}\ \ ((i,r),(j,s)\in I_0\times(\Z\setminus\{0\})),\leqno{(HH)}$$
$$[X_{i,r\pm\tilde d_{ij}}^{\pm},X_{j,s}^{\pm}]_{q_i^{a_{ij}}}
+[X_{j,s\pm\tilde d_{ij}}^{\pm},X_{i,r}^{\pm}]_{q_j^{a_{ji}}}=0\ ((i,r),(j,s)\in I_0\times\Z,\ 
a_{ij}\!\!<\!\!0),\leqno{(XD^{\pm})}$$
$$
\sum_{\sigma\in\sy_2}\sigma.[X_{i,r_1\pm\tilde d_i}^{\pm},X_{i,r_2}^{\pm}]_{q_i^{2}}=0\ \ ((r_1,r_2)\in\Z^2,\ (X_{\tilde n}^{(k)},i)\neq(A_{2n}^{(2)},1)),\leqno{(X1^{\pm})}$$
$$
\sum_{\sigma\in\sy_2}\sigma.([X_{1,r_1\pm 2}^{\pm},X_{1,r_2}^{\pm}]_{q^{2}}-q^{4}[X_{1,r_1\pm 1}^{\pm},X_{1,r_2\pm 1}^{\pm}]_{q^{-6}})=0\leqno{(X2^{\pm})}$$
$\ \ \ \ ((r_1,r_2)\in\Z^2,\ X_{\tilde n}^{(k)}=A_{2n}^{(2)}),$
$$\sum_{\sigma\in\sy_3}\sigma.[[X_{1,r_1\pm\varepsilon}^{\pm},X_{1,r_2}^{\pm}]_{q^{2\varepsilon}},X_{1,r_3}^{\pm}]_{q^{4\varepsilon}}=0\leqno{(X3^{\varepsilon,\pm})}$$
$\ \ \ \ ((r_1,r_2,r_3)\in\Z^3,\ X_{\tilde n}^{(k)}=A_{2n}^{(2)}),$
$$\sum_{{\sigma\in\sy_{1-a_{ij}}}}\!\!\!
\!\sigma.\!
\sum_{u=0}^{1-a_{ij}}(-1)^u{1-a_{ij}\brack u}_{q_i}X_{i,r_{1}}^{\pm}\cdot...\cdot X_{i,r_{u}}^{\pm}X_{j,s}^{\pm}X_{i,r_{u+1}}^{\pm}\cdot...\cdot X_{i,r_{1-a_{ij}}}^{\pm}=0\leqno{(SUL^{\pm})}$$
$$\ \ \ \ (i\neq j\in\! I_0,\ a_{ij}\in\{0,-1\}\ {\roman {if}}\ k\neq 1,\ r=(r_1,...,r_{1-a_{ij}})\in
\Z^{1-a_{ij}},\ s\in\Z),$$
$$\sum_{\sigma\in\sy_2}\!\!\sigma.\big(q(X_{j,s}^{\pm}X_{i,r_1\pm 1}^{\pm}X_{i,r_2}^{\pm}\!-[2]_{q^2}X_{i,r_1\pm 1}^{\pm}X_{j,s}^{\pm}X_{i,r_2}^{\pm}\!+X_{i,r_1\pm 1}^{\pm}X_{i,r_2}^{\pm}X_{j,s}^{\pm})+\leqno(S2^{\pm})$$
$$+q^{-1}(X_{j,s}^{\pm}X_{i,r_1}^{\pm}X_{i,r_2\pm 1}^{\pm}-[2]_{q^2}X_{i,r_1}^{\pm}X_{j,s}^{\pm}X_{i,r_2\pm 1}^{\pm}+X_{i,r_1}^{\pm}X_{i,r_2\pm 1}^{\pm}X_{j,s}^{\pm})\big)=0$$
$$\ \ \ \ (i,j\in I_0,\ a_{ij}=-2,\ k=2,\ (r_1,r_2)\in\Z^{2},\ s\in\Z),$$
$$\sum_{\sigma\in\sy_2}\sigma.\Big(q^2(X_{j,s}^{\pm}X_{i,r_1\pm 2}^{\pm}X_{i,r_2}^{\pm}-[2]_{q^3}X_{i,r_1\pm 2}^{\pm}X_{j,s}^{\pm}X_{i,r_2}^{\pm}+\leqno(S3^{\pm})$$
$$+X_{i,r_1\pm 2}^{\pm}X_{i,r_2}^{\pm}X_{j,s}^{\pm})+(X_{j,s}^{\pm}X_{i,r_1\pm 1}^{\pm}X_{i,r_2\pm 1}^{\pm}-[2]_{q^3}X_{i,r_1\pm 1}^{\pm}X_{j,s}^{\pm}X_{i,r_2\pm 1}^{\pm}+$$
$$+X_{i,r_1\pm 1}^{\pm}X_{i,r_2\pm 1}^{\pm}X_{j,s}^{\pm})+q^{-2}(X_{j,s}^{\pm}X_{i,r_1}^{\pm}X_{i,r_2\pm 2}^{\pm}+$$
$$-[2]_{q^3}X_{i,r_1}^{\pm}X_{j,s}^{\pm}X_{i,r_2\pm 2}^{\pm}+X_{i,r_1}^{\pm}X_{i,r_2\pm 2}^{\pm}X_{j,s}^{\pm})\Big)=0$$
$$\ \ \ \ (i,j\in I_0,\ a_{ij}=-3,\ k=3,\ (r_1,r_2)\in\Z^{2},\ s\in\Z),$$
where $\varepsilon\in\{\pm 1\}$ and 
$\tilde H_{i,r}^{\pm}$ and 
$b_{ijr}$ are defined as follows:
$$\sum_{r\in\Z}\tilde H_{i,\pm r}^{\pm}u^r=exp\left(\pm(q_i-q_i^{-1})
\sum_{r>0}H_{i,\pm r}u^r\right);$$
$$b_{ijr}=\begin{cases}  0&{\roman{if}}\ \tilde d_{i,j}\not|r\cr
{[2r]_q(q^{2r}+(-1)^{r-1}+q^{-2r})\over r}&{\roman{if}}\ (X_{\tilde n}^{(k)},i,j)=(A_{2n}^{(2)},1,1)\cr
{[\tilde r a_{ij}]_{q_i}\over\tilde r}&{\roman{otherwise,\ with\ }}
\tilde r={r\over \tilde d_{i,j}}.\end{cases} $$

An isomorphism is given by
$$C^{\pm 1}\mapsto {\Cal C}^{\pm 1},\ \  k_i^{\pm 1}\mapsto {\Cal K}_{\tilde i}^{\pm 1},\ \ X_{i,r}^{\pm}\mapsto {\Cal X}_{\tilde i,r}^{\pm},\ \ H_{i,s}\mapsto {\Cal H}_{\tilde i,s}$$
where $(i,r),(i,s)\in I_0\times\Z$ ($s\neq 0)$ and \ $\tilde{}:I_0\to\tilde I$ is a section as in remark \coefeq; its inverse is 
$${\Cal C}^{\pm 1}\mapsto C^{\pm 1},\ \  {\Cal K}_{i^{\prime}}^{\pm 1}\mapsto k_{\bar{i^{\prime}}}^{\pm 1},\ \ {\Cal X}_{\chi^u(\tilde i),r}^{\pm}\mapsto \omega^{ur}X_{i,r}^{\pm},\ \ {\Cal H}_{\chi^u(\tilde i),s}\mapsto \omega^{us}H_{i,s}$$
($i^{\prime}\in\tilde I$, $i\in I_0$, $u,r\in\Z$, $s\in\Z\setminus\{0\}$).
\dim
The claim follows from remarks \resp, corollary \zig, remarks \calv-\cpbijr, proposition \cntb\ and remarks \rsfgcal-\pspol.

\vskip .3 truecm

\rem{\greq}~

$\Udr_q(X_{\tilde n}^{(k)})$ is (isomorphic to) the $\C(q)$-algebra generated by 
$$C^{\pm1},\ \ \ k_i^{\pm1}\ \ (i\in I_0),\ \ \ X_{i,r}^{\pm}\ \ ((i,r)\in I_{\Z}),\ \ \ H_{i,r}\ \ ((i,r)\in I_{\Z}\setminus(I_0\times\{0\})),\leqno{(G^{\prime})}$$
with relations $(CUK^{\prime})$-$(S3^{\prime\pm})$, where,
for a relation $(R)$, the relation $(R^{\prime})$ is the set of relations in $(R)$ whose left hand side does not involve indices in $(I_0\times\Z)\setminus I_{\Z}$. 

Remark that the only case where the right hand side of some relation in $(R^{\prime})$ involves indices in $(I_0\times\Z)\setminus I_{\Z}$ is the case $(R)=(HX^{\pm})$: in this situation if $(j,r+s)\not\in I_{\Z}$ then $\tilde d_j\not|r$ and $b_{ijr}=0$, hence $(HX^{\prime\pm})$ is the following relation:
$$ [H_{i,r},X_{j,s}^{\pm}]=\begin{cases}  0&{\roman {if}\ } \tilde d_j\not|r\cr\pm b_{ijr}C^{r\mp|r|\over 2}X_{j,r+s}^{\pm}&{\roman {if}\ } \tilde d_j|r\end{cases} \ \ ((i,r)\in I_{\Z}\setminus(I_0\times\{0\}),(j,s)\in I_{\Z}).$$

\vskip .3 truecm

\rem{\prnpr}~
Since in the $\C(q)$-algebra generated by $(G)$ for any of the relations $(R)$ defining $\Udr_q(X_{\tilde n}^{(k)})$ the relations $(ZX,ZH,R)$ are equivalent to the relations $(ZX,ZH,R^{\prime})$, by abuse of notation we shall denote by $(R)$ also the relation $(R^{\prime})$.

\vskip .3 truecm
It is with the presentation of $\Udr_q$ given in proposition \drdef\ that we shall deal from now on.

\vskip 2 truecm

{\bf{\drmr.\ MORE about the DEFINITION of $\Udr_q$.}}
\vskip .5 truecm

The material of this section is presented in order to simplify the exposition and to handle more easily  the relations defining $\Udr_q$, with the aim of sharply reducing them:
some notations will be fixed; a new formulation will be given, mainly in terms of $q$-commutators, of some of the relations of proposition \drrl.\drdef; and some new relations ($(T2^{\pm})$ and $(T3^{\pm})$) will be introduced and proved to be equivalent, under suitable conditions, to $(S2^{\pm})$ and $(S3^{\pm})$. Also the Serre relations are introduced here, but they will be studied in details in section \serrel.

\vskip .3 truecm
\nota{\aaa}
Let ${\Cal U}$ be an algebra and let $(R)$ denote the relations 
$$S_{\zeta}(r,s)=0\ \ \ (\zeta\in{\Cal Z}, r\in\Z^l,s\in\Z^{\tilde l}),\leqno{(R)}$$
where ${\Cal Z}$ is a set, $l\in\Z_+$, $\tilde l\in\{0,1\}$, $S_{\zeta}(r,s)\in
{\Cal U}$. Then:

i) for all $\zeta\in{\Cal Z}$, denote by $(R_{\zeta})$ the relations
$$S_{\zeta}(r,s)=0\ \ \ (r\in\Z^l,s\in\Z^{\tilde l});\leqno{(R_{\zeta})}$$
of course if $\#{\Cal Z}=1$ and ${\Cal Z}=\{\zeta\}$ then 
$(R)=(R_{\zeta})$;

ii) denote by  ${\Cal I}(R)$ the ideal of 
${\Cal U}$ generated by the $S_{\zeta}(r,s)$'s: 
$${\Cal I}(R)=(S_{\zeta}(r,s)|\zeta\in{\Cal Z}, r\in\Z^l,s\in\Z^{\tilde l});$$
of course ${\Cal I}(R)=({\Cal I}(R_{\zeta})|\zeta\in{\Cal Z})$; 

iii) if $(^{(h)}\!R)$ ($h=1,...,m$) are the relations 
$$^{(h)}\!S_{\zeta}(r,s)=0\ \ \ (\zeta\in^{(h)}\!\!\!{\Cal Z}, r\in\Z^{l_h},s\in\Z^{\tilde l_h}),\leqno{(^{(h)}\!R)}$$
where $^{(h)}\!{\Cal Z}$ is a set, $l_h\in\Z_+$, $\tilde l_h\in\{0,1\}$, $^{(h)}\!S_{\zeta}^{\pm}(r,s)\in{\Cal U}$, define
$${\Cal I}(^{(1)}\!R,...,^{(m)}\!R)=({\Cal I}(^{(1)}\!R),...,{\Cal I}(^{(m)}\!R))$$

iv) if $(R^{\pm})$ denotes the relations 
$$S_{\zeta}^{\pm}(r,s)=0\ \ \ (\zeta\in{\Cal Z}, r\in\Z^l,s\in\Z^{\tilde l}),\leqno{(R^{\pm})}$$
where ${\Cal Z}$ is a set, $l\in\Z_+$, $\tilde l\in\{0,1\}$, $S_{\zeta}^{\pm}(r,s)\in
{\Cal U}$,
denote by 
$(R)$ the relations 
$$S_{\zeta^{\prime}}(r,s)=0\ \ \ (\zeta^{\prime}\in{\Cal Z}\times\{\pm\}, r\in\Z^l,s\in\Z^{\tilde l}),\leqno{(R)}$$
where $S_{(\zeta,\pm)}(r,s)=S_{\zeta}^{\pm}(r,s)$; 
in particular 
$${\Cal I}(R)=({\Cal I}(R^{+}),{\Cal I}(R^{-}));$$
moreover denote by ${\Cal I}^{\pm}(R)$ the ideals
$${\Cal I}^+(R)={\Cal I}(R^+)\ \ \ {\roman{and}}\ \ \ {\Cal I}^-(R)={\Cal I}(R^-).$$
\vskip .3truecm
\nota{\scmp}
For $i,j\in I_0$, $l\in\N$, $a\in\Z$, $r=(r_1,...,r_l)\in\Z^l$, $s\in\Z$ we set $$X_{i,j;l;a}^{\pm}(r;s)=\sum_{u=0}^l(-1)^u{l\brack u}_{q_i^a}
X_{i,r_1}^{\pm}\cdot...\cdot X_{i,r_u}^{\pm}X_{j,s}^{\pm}X_{i,r_{u+1}}^{\pm}\cdot...\cdot X_{i,r_l}^{\pm}.$$
\vskip .3truecm
\rem{\scpm}~
The relations $(SUL^{\pm})$, $(S2^{\pm})$ and $(S3^{\pm})$ can be written in a more compact form as: 
$$\sum_{{\sigma\in\sy_{1-a_{ij}}}}
\sigma.X_{i,j;1-a_{ij};1}^{\pm}(r;s)=0,$$
which is $(SUL^{\pm})$,
and 
$$\sum_{{\sigma\in\sy_2}}
\sigma.\sum_{u,v\geq 0\atop u+v=-1-a_{ij}}
q^{v-u}X_{i,j;2;-a_{ij}}^{\pm}(r_1\pm v,r_2\pm u;s)=0,$$
which is $(S2^{\pm})$, $(S3^{\pm})$ and also $(SUL^{\pm})$ in the case $a_{ij}=-1$.
\vskip .3truecm

In order to express the relations $(SUL^{\pm})$ in terms of $q$-commutators, and for further use and simplifications, we introduce the following notation.
\vskip .3truecm

\nota{\nnmm}
For $i\neq j\in I_0$, $l\in\N$, $a\in\Z$, $r=(r_1,...,r_l)\in\Z^l$, $s\in\Z$ set
$$M_{i,j;l;a}^{\pm}(r_1,...,r_{l};s)=\begin{cases}  X_{j,s}^{\pm}&{\roman{if}}\ l=0\cr
[M_{i,j;l-1;a}^{\pm}(r_1,...,r_{l-1};s),X_{i,r_l}^{\pm}]_{q_i^{-a_{ij}-2a(l-1)}}&{\roman{if}}\ l>0.\end{cases} $$

\vskip .3truecm
\rem{\sercol}~
The relations $(SUL^{\pm})$ can be formulated in terms of $q$-commutators as
$$\sum_{\sigma\in\sy_{1-a_{ij}}}\sigma.M_{i,j;1-a_{ij};1}^{\pm}(r;s)=0$$
$(i\neq j\in I_0,\ a_{ij}\in\{0,-1\}\ {\roman {if}}\ k\neq 1,\ r\in
\Z^{1-a_{ij}},\ s\in\Z)$.

\vskip .3truecm
\rem{\sercod}~
Also the relations $(S2^{\pm})$ and $(S3^{\pm})$ can be formulated in terms of $q$-commuta\-tors: 

i) $(S2^{\pm})$ can be written as:
$$\sum_{\sigma\in\sy_2}\sigma.\big((q^2+q^{-2})[[X_{j,s}^{\pm},X_{i,r_1\pm 1}^{\pm}]_{q^2},X_{i,r_2}^{\pm}]+q^2[[X_{i,r_1\pm 1}^{\pm},X_{i,r_2}^{\pm}]_{q^{ 2}},X_{j,s}^{\pm}]_{q^{- 4}}\big)=0;$$

ii) moreover $(S2^{+})$ can be written also in one of the following equivalent ways:
$$\sum_{\sigma\in\sy_2}\sigma.\big((q^2+q^{-2})[[X_{j,s}^{+},X_{i,r_1}^{+}]_{q^{- 2}},X_{i,r_2+1}^{+}]+[X_{j,s}^{+},[X_{i,r_2+1}^{+},X_{i,r_1}^{+}]_{q^{ 2}}]_{q^{-4}}\big)=0;$$
$$\sum_{\sigma\in\sy_2}\sigma.\big([[X_{j,s}^{+},X_{i,r_1+1}^{+}]_{q^{- 2}},X_{i,r_2}^{+}]-q^{2}[X_{i,r_1+1}^{+},[X_{j,s}^{+},X_{i,r_2}^{+}]_{q^{-2}}]_{q^{-4}}\big)=0;$$

iii) $(S3^{\pm})$ can be formulated in terms of $q$-commutators as follows:
$$\sum_{\sigma\in\sy_2}\!\!\sigma.\big((q^2+q^{-4})[[X_{j,s}^{\pm},\!X_{i,r_1\pm 2}^{\pm}]_{q^3},\!X_{i,r_2}^{\pm}]_{q^{- 1}}+$$
$$+
(1-q^{-2}+q^{-4})[[X_{j,s}^{\pm},\!X_{i,r_1\pm 1}^{\pm}]_{q^3},\!X_{i,r_2\pm 1}^{\pm}]_q+$$
$$+q^2[[X_{i,r_1\pm 2}^{\pm},X_{i,r_2}^{\pm}]_{q^{ 2}}+[X_{i,r_2\pm 1}^{\pm},X_{i,r_1\pm 1}^{\pm}]_{q^{ 2}},X_{j,s}^{\pm}]_{q^{- 6}}\big)=0;$$

\vskip .3truecm

\ddefi{\skq}
Consider the case $k>1$, $X_{\tilde n}^{(k)}\neq A_{2 n}^{(2)}$ and introduce the relations $(Tk^{\pm})$: 
$$\sum_{\sigma\in\sy_2}\sigma.[[X_{j,s}^{\pm},X_{i,r_1\pm 1}^{\pm}]_{q^{2}},X_{i,r_2}^{\pm}]=0\ \ \ (i,j\in I_0,\ a_{ij}=-2,\ r\in\Z^2,\ s\in\Z);\leqno{(T2^{\pm})}$$
$$\sum_{\sigma\in\sy_2}\sigma.((q^{2}+1)[[X_{j,s}^{\pm},X_{i,r_{1}\pm 2}^{\pm}]_{q^{3}},X_{i,r_{2}}^{\pm}]_{q^{-1}}+[[X_{j,s}^{\pm},X_{i,r_{1}\pm 1}^{\pm}]_{q^{3}},X_{i,r_{2}\pm 1}^{\pm}]_{q})=0\leqno{(T3^{\pm})}$$
$\ \ \ (i,j\in I_0,\ a_{ij}=-3,\ r\in\Z^2,\ s\in\Z)$.

\vskip .3truecm

\prop{\xuit}
Let $k>1$, $X_{\tilde n}^{(k)}\neq A_{2n}^{(2)}$.

Then ${\Cal I}(X1^{\pm},Sk^{\pm})={\Cal I}(X1^{\pm},Tk^{\pm})$. 

More precisely if $i,j\in I_{0}$ are such that $a_{ij}<-1$ we have that  $${\Cal I}(X1_i^{\pm},Sk^{\pm})={\Cal I}(X1_i^{\pm},Tk^{\pm})\ \ \ ({\roman {see\ notation}}\ \aaa,i)).$$ 

In particular, $(S2^{\pm})$ and $(S3^{\pm})$ can be replaced respectively by $(T2^{\pm})$ and $(T3^{\pm})$ among the defining relations of ${\Cal U}_q^{Dr}$.
\dim
It is enough to notice that 
$$\Big[\sum_{\sigma\in\sy_2}\sigma.[X_{i,r_1\pm 1}^{\pm},X_{i,r_2}^{\pm}]_{q^{ 2}},X_{j,s}^{\pm}\Big]_{q^{- 4}}$$
and
$$[[X_{i,r_1\pm 2}^{\pm},X_{i,r_2}^{\pm}]_{q^{ 2}}+[X_{i,r_2\pm 1}^{\pm},X_{i,r_1\pm 1}^{\pm}]_{q^{ 2}},X_{j,s}^{\pm}]_{q^{- 6}}$$
belong to ${\Cal I}(X1_i^{\pm})$.

\vskip .3 truecm

\ddefi{\relser}
We recall also the Serre relations
$$\sum_{{\sigma\in\sy_{1-a_{ij}}}}
\sigma.X_{i,j;1-a_{ij};1}^{\pm}(r;s)=0\ \ \ (i\neq j\in I_0,\ r\in\Z^{1-a_{ij}},\ s\in\Z).\leqno{(S^{\pm})}$$
\vskip .3 truecm

\rem{\sercom}~
The Serre relations can be formulated in terms of $q$-commutators as
$$\sum_{{\sigma\in\sy_{1-a_{ij}}}}
\sigma.M_{i,j;1-a_{ij};1}^{\pm}(r;s)=0\ \ \ (i\neq j\in I_0,\ r\in\Z^{1-a_{ij}},\ s\in\Z).$$

\rem{\sercot}~
The right hand sides of relations $(Tk^{\pm})$ and $(S^{\pm})$ are zero, hence remark \drrl.\prnpr\ holds for these relations (see also remark \drrl.\greq).

\vskip .3 truecm
The comparison of the defining relations of ${\Cal U}_q^{Dr}$ with the Serre relations is the matter of section \serrel. 

\vskip .3 truecm
\nota{\mmnn}
Let us introduce also the following notations: 

i) for $i,j\in I_0$, $r,s\in\Z$
$$M_{(2)}^{\pm}((i,r),(j,s))=
[X_{i,r\pm\tilde d_{ij}}^{\pm},X_{j,s}^{\pm}]_{q_i^{a_{ij}}}+[X_{j,s\pm\tilde d_{ij}}^{\pm},X_{i,r}^{\pm}]_{q_j^{a_{ji}}};$$

ii) for $i\in I_0$, $r=(r_1,r_2)\in\Z^2$
$$M_i^{\pm}(r)=
[X_{i,r_1\pm\tilde d_{i}}^{\pm},X_{i,r_2}^{\pm}]_{q_i^{2}};$$

iii) if $X_{\tilde n}^{(k)}=A_{2n}^{(2)}$ and $r=(r_1,r_2)\in\Z^2$
$$M_{(2,2)}^{\pm}(r)=[X_{1,r_1\pm 2}^{\pm},X_{1,r_2}^{\pm}]_{q^{ 2}}-q^{ 4}[X_{1,r_1\pm 1}^{\pm},X_{1,r_2\pm 1}^{\pm}]_{q^{-6}};$$

iv) if $X_{\tilde n}^{(k)}=A_{2n}^{(2)}$ and $r=(r_1,r_2,r_3)\in\Z^3$
$$M_{(3)}^{\varepsilon,\pm}(r)=
[[X_{1,r_1\pm\varepsilon}^{\pm},X_{1,r_2}^{\pm}]_{q^{2\varepsilon}},X_{1,r_3}^{\pm}]_{q^{4\varepsilon}};$$ 

v) if $k>1$ and $r=(r_1,r_2)\in\Z^2$, $s\in\Z$
$$X_{[k]}^{\pm}(r;s)=\sum_{u,v\geq 0\atop u+v=k-1}q^{v-u}
X_{i,j;2;k}^{\pm}(r_1\pm v,r_2\pm u;s)$$
where $i,j\in I_0$ are such that $a_{ij}=-k$;

vi) if $k=2$, $X_{\tilde n}^{(k)}\neq A_{2n}^{(2)}$  and $r=(r_1,r_2)\in\Z^2$, $s\in\Z$
$$M_{[2]}^{\pm}(r;s)=M_{i,j;2;1}^{\pm}(r_1\pm 1,r_2;s)$$
where $i,j\in I_0$ are such that $a_{ij}=-2$;

vii) if $k=3$ and $r=(r_1,r_2)\in\Z^2$, $s\in\Z$
$$M_{[3]}^{\pm}(r;s)=(q^2+1)M_{i,j;2;2}^{\pm}(r_1\pm 2,r_2;s)+M_{i,j;2;1}^{\pm}(r_1\pm 1,r_2\pm 1;s)$$
where $i,j\in I_0$ are such that $a_{ij}=-3$. 

\vskip .3 truecm
\rem{\crnt}~
Of course the following relations depend on $(ZX^{\pm})$: 

i) $M_{(2)}^{\pm}((i,r),(j,s))=0$ if $(r,s)\not\in\tilde d_i\Z\times\tilde d_j\Z$;

ii) $M_i^{\pm}(r)=0$ if $r\not\in(\tilde d_i\Z)^2$;

iii) $M_{i,j;l;a}^{\pm}(r;s)=0$ if $(r,s)\not\in(\tilde d_i\Z)^{l}\times\tilde d_j\Z$;

iv) $X_{i,j;l;a}^{\pm}(r;s)=0$ if $(r,s)\not\in(\tilde d_i\Z)^{l}\times\tilde d_j\Z$;

v) $X_{[k]}^{\pm}(r;s)=0$ if $s\not\in\tilde d\Z$;

vi) $M_{[k]}^{\pm}(r;s)=0$ if $s\not\in k\Z$.

\vskip .3 truecm
\rem{\tcalb}~
Recalling remark \drrl.\prnpr\ (and remark \crnt) we have the following obvious reformulation of the relations $(XD^{\pm})$-$(S3^{\pm})$, $(T2^{\pm})$, $(T3^{\pm})$ and $(S^{\pm})$ in terms of the notations just introduced (notations \scmp, \nnmm, \mmnn):
$$M_{(2)}^{\pm}((i,\tilde d_ir),(j,\tilde d_js))=0\ \ \ (i,j\in I_0,\ a_{ij}<0,\ r,s\in\Z);\leqno(XD^{\pm}):$$
$$\sum_{\sigma\in\sy_2}\sigma.M_i^{\pm}(\tilde d_i r)=0\ \ \ (i\in I_0,\ (X_{\tilde n}^{(k)},i)\neq(A_{2n}^{(2)},i),\ r\in\Z^2);\leqno(X1^{\pm}):$$
$$\sum_{\sigma\in\sy_2}\sigma.M_{(2,2)}^{\pm}(r)=0\ \ \ (r\in\Z^2);\leqno(X2^{\pm}):$$
$$\sum_{\sigma\in\sy_3}\sigma.M_{(3)}^{\varepsilon,\pm}(r)=0\ \ \ (r\in\Z^3);\leqno(X3^{\varepsilon,\pm}):$$
$$\sum_{\sigma\in\sy_{1-a_{ij}}}\sigma.M_{i,j;1-a_{ij};1}^{\pm}(\tilde d_ir;\tilde d_js)=0\leqno(S(UL)^{\pm}):$$ or equivalently
$$\sum_{{\sigma\in\sy_{1-a_{ij}}}}
\sigma.X_{i,j;1-a_{ij};1}^{\pm}(\tilde d_ir;\tilde d_js)=0$$
($i\neq j\in I_0\ (a_{ij}\in\{0,-1\}\ {\roman{if}}\ k\neq 1),\ r\in\Z^{1-a_{ij}},\ s\in\Z$);
$$\sum_{{\sigma\in\sy_2}}\sigma.X_{[k]}^{\pm}(r;\tilde d s)=0\ \ \ (r\in\Z^2,\ s\in\Z);\leqno(Sk^{\pm}): $$
$$\sum_{{\sigma\in\sy_2}}
\sigma.M_{[k]}^{\pm 1}(r;ks)=0\ \ \ (r\in\Z^2,\ s\in\Z).\leqno(Tk^{\pm}):$$

\vskip .5 truecm

{\bf{\drcmp.\ $\bar{\Cal U}_q^{Dr}$ and its STRUCTURES.}}
\vskip .5 truecm

In order to study the relations defining ${\Cal U}_q^{Dr}$ it is convenient to proceed by steps: the algebras $\bar{\Cal U}_q^{Dr}=\bar{\Cal U}_q^{Dr}(X_{\tilde n}^{(k)})$ and $\tilde{\Cal U}_q^{Dr}=\tilde{\Cal U}_q^{Dr}(X_{\tilde n}^{(k)})$ here defined are such that $\tilde{\Cal U}_q^{Dr}$ is a quotient of $\bar{\Cal U}_q^{Dr}$ and ${\Cal U}_q^{Dr}$ is a quotient of $\tilde{\Cal U}_q^{Dr}$.

This section is devoted to introduce some important structures on $\bar{\Cal U}_q^{Dr}$ ($Q$-gradation, homomorphisms between some of these algebras, automorphisms and antihomomorphims of each of them), which will be proved to induce analogous structures on $\tilde{\Cal U}_q^{Dr}$ (see also section \darl) and, which is finally important, on ${\Cal U}_q^{Dr}$ (see also sections \darl\ and \dautemb). 

Some remarks point out the first (trivially) unnecessary relations: $(ZH)$ and $(KH)$ are redundant.

\vskip .3 truecm
\ddefi{\tilu}
We denote by:

i) $\tilde{\Cal U}^{Dr}_q(X_{\tilde n}^{(k)})$  the $\C(q)$-algebra  generated by
$(G)$
with relations $$(ZX^{\pm}),\ (CUK),\ (CK),\ (KX^{\pm}),\ (XX)\ {\roman {and}}$$
$$ [H_{i,r},X_{j,s}^{\pm}]=\pm b_{ijr}C^{r\mp|r|\over 2}X_{j,r+s}^{\pm}\ \ ((i,r),(j,s)\in I_{\Z},\ \tilde d_i\leq|r|\leq\tilde d_{ij}
);\leqno(HXL^{\pm})$$

ii) $\bar{\Cal U}^{Dr}_q(X_{\tilde n}^{(k)})$  the $\C(q)$-algebra  generated by
$$C^{\pm1},\ \ \ k_i^{\pm1}\ \ (i\in I_0),\ \ \ X_{i,r}^{\pm}\ \ ((i,r)\in I_0\times{\Z})
\leqno(\bar{G})$$
with relations $$(ZX^{\pm}),\ (CUK),\ (CK).$$

\rem{\bvqz}~
$\Udr_q(X_{\tilde n}^{(k)})$ is obviously a quotient of  $\tilde{\Cal U}^{Dr}_q(X_{\tilde n}^{(k)})$.
\vskip .3 truecm

We shall prove that also  $\tilde{\Cal U}^{Dr}_q(X_{\tilde n}^{(k)})$ is a quotient of  $\bar{\Cal U}^{Dr}_q(X_{\tilde n}^{(k)})$. 

Since $\bar{\Cal U}^{Dr}_q(X_{\tilde n}^{(k)})\to\tilde{\Cal U}^{Dr}_q(X_{\tilde n}^{(k)})$ is obviously well defined, we just need to prove that this map is surjective, or equivalently that 
$\tilde{\Cal U}^{Dr}_q(X_{\tilde n}^{(k)})$ is generated by $(\bar{G})$.
 To this aim we need some simple remarks.

\vskip.3 truecm
\rem{\greh}~
In $\tilde{\Cal U}^{Dr}_q$ (hence in $\Udr_q$) the following holds: 

i) $\tilde H_{i,0}^{\pm}=1$ $\forall i\in I_0$; 

ii) $\tilde H_{i,\mp r}^{\pm}=0$ $\forall i\in I_0$, $\forall r>0$; 

iii) $\forall r>0$ $\tilde H_{i,\pm r}^{\pm}\mp (q_i-q_i^{-1})H_{i,\pm r}$ belongs to the $\C(q)$-subalgebra 
generated by $\{H_{i,\pm s}|0<s<r\}$; in particular $\{H_{i,\pm s}|(i,s)\in I_{\Z},\ 0<s<r\}$ and $\{\tilde H_{i,\pm s}^{\pm}|(i,s)\in I_{\Z},\ 0<s<r\}$ generate the same $\C(q)$-subalgebra. 
\vskip.3 truecm
\rem{\tiht}~
In $\tilde{\Cal U}^{Dr}_q$ (hence in $\Udr_q$) we have that for all $i\in I_0$ and for all $r\in\Z_+$
$$\tilde H_{i,\pm r}^{\pm}=
(q_i-q_i^{-1})k_i^{\mp 1}[X_{i,\pm r}^{\pm},X_{i,0}^{\mp}].$$

In particular  for all $(i,r)\in I_0\times\Z$ $\tilde H_{i,r}^{\pm}$ lies in the subalgebra of $\tilde{\Cal U}^{Dr}_q(X_{\tilde n}^{(k)})$ generated by $(\bar G)$. 

Consequently (see remark \greh) for all $(i,r)\in I_0\times\Z$ also  $H_{i,r}$ lies in the subalgebra of $\tilde{\Cal U}^{Dr}_q(X_{\tilde n}^{(k)})$ generated by $(\bar{G})$.

\vskip.3 truecm

\cor{\genqz}
i) $\tilde{\Cal U}^{Dr}_q(X_{\tilde n}^{(k)})$ and $\Udr_q(X_{\tilde n}^{(k)})$ are generated by $(\bar{G})$;

ii) $\tilde{\Cal U}^{Dr}_q(X_{\tilde n}^{(k)})$ is a quotient of  $\bar{\Cal U}^{Dr}_q(X_{\tilde n}^{(k)})$.

\vskip .3 truecm
\nota{\bthi}
We denote by $\tilde H_{i,\pm r}^{\pm}$ also the elements in $\bar{\Cal U}^{Dr}_q$ defined by
$$\tilde H_{i,\pm r}^{\pm}=\begin{cases} 
(q_i-q_i^{-1})k_i^{- 1}[X_{i, r}^{+},X_{i,0}^{-}]&{\roman{ if}}\ r,\pm r>0\cr
(q_i-q_i^{-1})[X_{i,-r}^{-},X_{i,0}^{+}]k_i&{\roman{ if}}\ r>0,\pm r<0\cr
1&{\roman{ if}}\ r=0\cr 0&{\roman{ if}}\ r<0,\end{cases} $$
and by $H_{i,r}$ the elements of $\bar{\Cal U}^{Dr}_q$ defined by
$$\sum_{r\in\Z}\tilde H_{i,\pm r}^{\pm}u^r=exp\left(\pm(q_i-q_i^{-1})
\sum_{r>0}H_{i,\pm r}u^r\right);$$

\vskip .3 truecm
\rem{\sepid}~
Relations $(ZH)$  are trivial in $\bar{\Cal U}^{Dr}_q(X_{\tilde n}^{(k)})$.

\vskip .3 truecm
\rem{\idstsep}~
i) $\bar{\Cal U}^{Dr}_q=\bar{\Cal U}^{Dr}_q(X_{\tilde n}^{(k)})$ is $Q$-graded: 
$$\bar{\Cal U}^{Dr}_q=\oplus_{\alpha\in Q}\bar{\Cal U}^{Dr}_{q,\alpha},$$
where $C^{\pm 1}$, $k_i^{\pm 1}\in\bar{\Cal U}^{Dr}_{q,0}$, $X_{i,r}^{\pm}\in\bar{\Cal U}^{Dr}_{q,\pm\alpha_i+r\delta}$ $\forall i\in I_0$, $r\in\Z$, and $\bar{\Cal U}^{Dr}_{q,\alpha}\bar{\Cal U}^{Dr}_{q,\beta}\subseteq\bar{\Cal U}^{Dr}_{q,\alpha+\beta}$..

ii) $\tilde H_{i,r}^{\pm}$ ($r\in\Z$) and $H_{i,r}$ ($r\in\Z\setminus\{0\}$) are homogeneous of degree $r\delta$ for all $i\in I_0$.

iii) Since the relations defining $\tilde{\Cal U}^{Dr}_q$ and $\Udr_q$ are homogeneous, the $Q$-gradation of $\bar{\Cal U}^{Dr}_q$ induces $Q$-gradations on $\tilde{\Cal U}^{Dr}_q=\oplus_{\alpha\in Q}\tilde{\Cal U}^{Dr}_{q,\alpha}$ and on $\Udr_q=\oplus_{\alpha\in Q}\Udr_{q,\alpha}$.

\vskip .3 truecm
\nota{\qgr}
The $\C(q)$-algebra $\C(q)[C^{\pm 1}, k_i^{\pm 1}|i\in I_0]$ is $Q$-graded, with one-dimensional homogeneous components $\C(q)k_{\alpha}$ ($\alpha\in Q$) where we set $$k_{m\delta+\sum_{i\in I_0}m_i\alpha_i}=C^m\prod_{i\in I_0}k_i^{m_i}\ \ \ (m,m_i\in\Z\ \forall i\in I_0).$$

Indeed $\C(q)[C^{\pm 1}, k_i^{\pm 1}|i\in I_0]=\C(q)[Q]$. 

Recall that $\C(q)[C^{\pm 1}, k_i^{\pm 1}|i\in I_0]$ naturally maps in 
$\bar{\Cal U}^{Dr}_{q,0}\subseteq\bar{\Cal U}^{Dr}_q$ (hence in $\tilde{\Cal U}^{Dr}_{q,0}\subseteq\tilde{\Cal U}^{Dr}_q$ and in $\Udr_{q,0}\subseteq\Udr_q$).

\vskip .3 truecm
\rem{\lpdgr}~
i) The relations $(CUK)$, $(CK)$ and $(KX)^{\pm}$ are equivalent to a) and b):

a) the $\C(q)$-subalgebra generated by $\{C^{\pm 1}, k_i^{\pm 1}|i\in I_0\}$  is a quotient of the ring of Laurent polynomials $\C(q)[k_i^{\pm 1}|i\in I]$ ($C=\prod_{i\in I}k_i^{r_i}$);

b) for all $\alpha$, $\beta\in Q$ and forall $x$ of degree $\beta$ we have $k_{\alpha}x=q^{(\alpha|\beta)}xk_{\alpha}$.

ii) Relations $(KH)$ depend on relations $(CUK)$, $(CK)$ and $(KX)^{\pm}$, and in particular  are trivial in $\tilde{\Cal U}^{Dr}_q(X_{\tilde n}^{(k)})$.

\vskip .3 truecm
\ddefi{\tilpm}
We denote by ${\Cal F}_q^{+}={\Cal F}_q^{+}(X_{\tilde n}^{(k)})$ and ${\Cal F}_q^{-}={\Cal F}_q^{-}(X_{\tilde n}^{(k)})$ the $\C(q)$-algebras  generated respectively by
$$X_{i,r}^+\ \ \ ((i,r)\in I_0\times\Z)\leqno{(G^+)}$$ 
and 
$$X_{i,r}^-\ \ \ ((i,r)\in I_0\times\Z)\leqno{(G^-)}$$ 
with relations respectively $(ZX^+)$ and $(ZX^-)$.

\vskip 2 truecm
\rem{\rtilpm}~
${\Cal F}_q^{+}(X_{\tilde n}^{(k)})$ and ${\Cal F}_q^{-}(X_{\tilde n}^{(k)})$ are the free $\C(q)$-algebras  generated respectively by
$$X_{i,r}^+\ \ \ ((i,r)\in I_{\Z})\leqno{({G}^{\prime+})}$$ 
and 
$$X_{i,r}^-\ \ \ ((i,r)\in I_{\Z}).\leqno{({G}^{\prime-})}$$

\vskip .3 truecm
\nota{\eqvdd}
${\Cal F}_q^{+}$ and ${\Cal F}_q^{-}$ naturally embed in $\bar{\Cal U}_q^{Dr}$, hence they map in $\tilde{\Cal U}_q^{Dr}$  and in ${\Cal U}_q^{Dr}$; their images  in $\tilde{\Cal U}_q^{Dr}$ are denoted respectively by $\tilde{\Cal U}_q^{Dr,+}$ and $\tilde{\Cal U}_q^{Dr,-}$,  and their images  in ${\Cal U}_q^{Dr}$ are denoted respectively by ${\Cal U}_q^{Dr,+}$ and ${\Cal U}_q^{Dr,-}$.

\vskip .3 truecm
\rem{\fuf}~
i) As subalgebras of $\bar{\Cal U}_q^{Dr}$,  ${\Cal F}_q^{+}$ inherits a $(Q_{0,+}\oplus\Z\delta)$-gradation and  ${\Cal F}_q^{-}$ inherits a $(-Q_{0,+}\oplus\Z\delta)$-gradation; 

ii) more precisely we have that 
$${\Cal F}_q^{\pm}\subseteq\C(q)\oplus\bigoplus_{\alpha\in Q_{0,+},\alpha\neq 0\atop m\in\Z}\bar{\Cal U}^{Dr}_{q,\pm\alpha+m\delta}$$
and similarly
$$\tilde{\Cal U}_q^{Dr,\pm}\subseteq\C(q)\oplus\bigoplus_{\alpha\in Q_{0,+},\alpha\neq 0\atop m\in\Z}\tilde{\Cal U}^{Dr}_{q,\pm\alpha+m\delta}$$
and $${\Cal U}_q^{Dr,\pm}\subseteq
\C(q)\oplus\bigoplus_{\alpha\in Q_{0,+},\alpha\neq 0\atop m\in\Z}\Udr_{q,\pm\alpha+m\delta}.$$

\vskip .3 truecm
The last part of this section is devoted to the definition of automorphisms and antiautomorphisms of the algebras just introduced, which make evident some symmetries in the generators and relations of $\Udr_q$. Thanks to these structures the study of the apparently very complicated relations defining $\Udr_q$ will be strongly simplified in sections \darl, \dlem\ and following.

The next definitions depend on the choice of an automorphism $\eta$ of $\C$. A short discussion about the choice of $\eta$ is outlined in remark \disceta.

\vskip .3truecm
\ddefi{\gendef}
Let us introduce the following homomorphisms and antihomomorphisms:

i) $\bar\Omega:\bar{\Cal U}^{Dr}_q\to\bar{\Cal U}^{Dr}_q$ is the 
anti-homomorphism defined on the generators by 
$$\bar\Omega\big|_{\C}=\eta,\ \ 
q\mapsto q^{-1},\ \ C^{\pm 1}\mapsto C^{\mp 1},\ \ k_i^{\pm 1}\mapsto k_i^{\mp 1},\ \ X_{i,r}^{\pm}\mapsto X_{i,-r}^{\mp}.$$

ii) $\Theta_{{\Cal F}}^{+}:{\Cal F}_q^{+}\to{\Cal F}_q^{+}$ and $\Theta_{{\Cal F}}^{-}:{\Cal F}_q^{-}\to{\Cal F}_q^{-}$ are the homomorphisms defined on the generators respectively by $$\Theta_{{\Cal F}}^{+}:\ \ \ \Theta_{{\Cal F}}^{+}\big|_{\C}=\eta,\ \ q\mapsto q^{-1},\ \ X_{i,r}^{+}\mapsto X_{i,-r}^{+}$$
and $$\Theta_{{\Cal F}}^{-}:\ \ \ \Theta_{{\Cal F}}^{-}\big|_{\C}=\eta,\ \ q\mapsto q^{-1},\ \ X_{i,r}^{-}\mapsto X_{i,-r}^{-}.$$

iii) 
$\bar\Theta:\bar{\Cal U}^{Dr}_q\to\bar{\Cal U}^{Dr}_q$ is the homomorphism defined on the generators by $$\bar\Theta\big|_{\C}=\eta,\ \ q\mapsto q^{-1},\ \ C^{\pm 1}\mapsto C^{\pm 1},\ \ k_i^{\pm 1}\mapsto k_i^{\mp 1},$$
$$ X_{i,r}^+\mapsto -X_{i,-r}^+k_iC^{-r},\ \ X_{i,r}^-\mapsto -k_i^{-1}C^{-r}X_{i,-r}^-.$$

iv) For all $i\in I_0$ $\bar\ti:\bar{\Cal U}^{Dr}_q\to\bar{\Cal U}^{Dr}_q$ is the $\C(q)$-homomorphism defined on the generators by $$C^{\pm 1}\mapsto C^{\pm 1},\ \ \ k_j^{\pm 1}\mapsto (k_jC^{-\delta_{ij}\tilde d_i})^{\pm 1},\ \ \ X_{j,r}^{\pm}\mapsto X_{j,r\mp\d_{ij}\tilde d_i}^{\pm}.$$

v) For $i\in I_0$ let 
$$\bar\phi_i:\begin{cases} \bar{\Cal U}^{Dr}_q(A_{1}^{(1)})\rightarrow\bar{\Cal U}^{Dr}_q(X_{\tilde n}^{(k)})&
{\roman{if}}\ (X_{\tilde n}^{(k)},i)\neq(A_{2n}^{(2)},1)\cr
\bar{\Cal U}^{Dr}_q(A_{2}^{(2)})\rightarrow\bar{\Cal U}^{Dr}_q(X_{\tilde n}^{(k)})&
{\roman{if}}\ (X_{\tilde n}^{(k)},i)=(A_{2n}^{(2)},1)\end{cases} 
$$ be the $\C$-homomorphisms defined on the generators as follows: 
$$q\mapsto q_i
,\ \ \ C^{\pm 1}\mapsto C^{\pm \tilde d_i},\ \ \ k^{\pm 1}\mapsto k_i^{\pm 1},\ \ \ X_r^{\pm}\mapsto X_{i,\tilde d_ir}^{\pm}.$$

\vskip .3truecm

\rem{\genbdef}~
It is immediate to notice that:

i) $\bar\Omega$,  $\Theta_{{\Cal F}}^{\pm}$, $\bar\Theta$, $\bar\ti$ and $\bar\phi_i$  are all well-defined;

ii) $\bar\Omega({\Cal F}_q^{\pm})={\Cal F}_q^{\mp}$;

iii) $\bar\Omega$ and $\bar\Theta$ are involutions of $\bar{\Cal U}^{Dr}_q$, if $\eta$ is an involution of $\C$;

iv) the $\bar\ti$'s are automorphisms of $\bar{\Cal U}^{Dr}_q$ (of infinite order) for all $i\in I_0$; more precisely $<\bar t_i|i\in I_0>\cong Z^{I_0}$;

v) the following commutation properties hold: 
$$\bar\Theta\bar\Omega=\bar\Omega\bar\Theta,\ \ \bar t_i\bar\Omega=\bar\Omega\bar t_i,\ \ \bar t_i\bar\Theta=\bar\Theta\bar t_i^{-1}\ \ {\roman{and}}\ \ \bar t_i\bar t_j=\bar t_j\bar t_i\ \ \ \forall i,j\in I_0$$ 
as maps of $\bar{\Cal U}^{Dr}_q(X_{\tilde n}^{(k)})$ into itself; moreover, $\forall i\in I_0$,  
$$\bar\Omega\bar\phi_i=\bar\phi_i\bar\Omega,\ \ \bar\Theta\bar\phi_i=\bar\phi_i\bar\Theta,\ \ \bar t_i\bar\phi_i=\bar\phi_i\bar t_1,\ \ \bar t_j\bar\phi_i=\bar\phi_i\ \ \forall j\in I_0\setminus\{i\}$$  as maps from $\bar{\Cal U}^{Dr}_q(A_1^{(1)})$ to $\bar{\Cal U}^{Dr}_q(X_{\tilde n}^{(k)})$ if $(X_{\tilde n}^{(k)},i)\neq(A_{2n}^{(2)},1)$, from $\bar{\Cal U}^{Dr}_q(A_2^{(2)})$ to $\bar{\Cal U}^{Dr}_q(X_{\tilde n}^{(k)})$ if $(X_{\tilde n}^{(k)},i)=(A_{2n}^{(2)},1)$;

vi) for all $\alpha=\beta+m\delta\in Q$ with $\beta\in Q_0$, $m\in\Z$ we have: 
$$\bar\Omega(k_{\alpha})=k_{-\alpha},\ \bar\Theta(k_{\beta+m\delta})=k_{-\beta+m\delta},\ \bar t_i(k_{\alpha})=k_{\lambda_i(\alpha)}$$
and
$$\bar\Omega(\bar{\Cal U}^{Dr}_{q,\alpha})=\bar{\Cal U}^{Dr}_{q,-\alpha},\ \bar\Theta(\bar{\Cal U}^{Dr}_{q,\beta+m\delta})=\bar{\Cal U}^{Dr}_{q,\beta-m\delta},\ \bar t_i(\bar{\Cal U}^{Dr}_{q,\alpha})=\bar{\Cal U}^{Dr}_{q,\lambda_i(\alpha)};$$
moreover for all $m_1,m\in\Z$ and for all $i\in I_0$
$$\bar\phi_i(k_{m_1\alpha_1+m\delta})=k_{m_1\alpha_i+\tilde d_im\delta}\ \ \ {\roman{and}}\ \ \ 
\bar\phi_i(\Udr_{q,m_1\alpha_1+m\delta}(A_*^{(*)}))\subseteq\Udr_{q,m_1\alpha_i+\tilde d_im\delta}(X_{\tilde n}^{(k)});$$

vii) on the elements $H_{i,r}$ and $\tilde H_{i,r}^{\pm}$ we have:
$$\bar\Omega(\tilde H_{i,r}^{\pm})=\tilde H_{i,-r}^{\mp},\ \ \bar\Omega(H_{i,r})=H_{i,-r}$$ 
and 
$$\bar\phi_i(\tilde H_{1,r}^{\pm})=\tilde H_{i,\tilde d_ir}^{\pm},\ \ \bar\phi_i(H_{1,r})=H_{i,\tilde d_ir}\ \ \ \forall i\in I_0.$$ 
\vskip .3truecm
\rem{\disceta}~
For the purpose of the present paper, the definition of $\bar\Omega$, $\Theta_{{\Cal{F}}}^{\pm}$, $\bar\Theta$ given in definition \gendef\ could be simplified by requiring these maps to be $\C$-linear (that is $\eta=id_{\C}$). But the choice of a non trivial automorphism $\eta$ of $\C$ becomes sometimes necessary, as when specializing $q$ at a complex value $\epsilon\neq\pm1$: indeed a homomorphism defined over $\C(q)$ (and mapping $q$ to $q^{-1}$) induces a homomorphism on the specialization at $\epsilon$ if and only if the ideal $(q-\epsilon)$ is stable; if, for example, $\epsilon$ is a root of 1, this could be obtained by choosing $\eta(z)=\bar z$ $\forall z\in\C$, that is by requiring the homomorphism to be $\C$-anti-linear. For this reason, from now on we suppose $\eta$ to be the conjugation on $\C$, that is $\bar\Omega$, $\Theta_{{\Cal{F}}}^{\pm}$, $\bar\Theta$ to be $\C$-anti-linear (see definitions \dautemb.\omdef\ and \dautemb.\xidef,  and compare also with definition \drjq.\TTTUQ).

Of course one needs to pay more attention and eventually to choose a different $\eta$ when interested in specializing at complex values $\epsilon$ such that $|\epsilon|\neq 1$.

\vskip .3truecm
Our goal is of course to show that $\bar\Omega$,  $\bar\Theta$, $\bar\ti$ and $\bar\phi_i$ induce 
$\Omega$,  $\Theta$, $\ti$ and $\phi_i$ on $\Udr_q$. This is indeed very easy to show, but we take this occasion to simplify the relations that we have to handle with, passing through $\tilde{\Cal U}_q^{Dr}$ for two reasons: underlying the first redundances of the relations (see corollary \darl.\dnd); discussing separately the relations $(XD)^{+}$-$(S3)^{+}$ 
whose first simplification can be dealt with simultaneously as examples of a general case (see section \dlem).
\vskip 2 truecm

{\bf{\darl.\ The ALGEBRA $\tilde{\Cal U}_q^{Dr}$.}}

\vskip .5truecm
The algebra $\tilde{\Cal U}_q^{Dr}$ and its structures, to which this section is devoted, play a fundamental role in the study and simplification of relations $(XD^{\pm})$-$(S3^{\pm})$. In particular 
the relations are analyzed underlining their consequences on the (anti)automorphisms $\tilde\Omega$, $\tilde\Theta$ and $\tilde t_i$ ($i\in I_0$); relations $(HX^{\pm})$ and $(HH)$ are proved to be redundant; and much smaller sets of generators are provided. 

\vskip .3 truecm
\rem{\rkhgr}~
Remarks \drcmp.\lpdgr,i) and \drcmp.\genbdef,vi) imply immediately that $\bar\Omega$,  $\bar\Theta$, $\bar\ti$ and $\bar\phi_i$ preserve relations $(KX)^{\pm}$.
\vskip .3 truecm
\rem{\fibdef}~
For all $i\in I_0$ $\bar\phi_i$ obviously induces 
$$\tilde \phi_i:\begin{cases} \tilde{\Cal U}^{Dr}_q(A_1^{(1)})\to\tilde{\Cal U}^{Dr}_q(X_{\tilde n}^{(k)})&{\roman{if}}\ (X_{\tilde n}^{(k)},i)\neq(A_{2n}^{(2)},1)\cr
\tilde{\Cal U}^{Dr}_q(A_2^{(2)})\to\tilde{\Cal U}^{Dr}_q(X_{\tilde n}^{(k)})&{\roman{if}}\ (X_{\tilde n}^{(k)},i)=(A_{2n}^{(2)},1)\end{cases} $$ 
and 
$$\phi_i:\begin{cases} \Udr_q(A_1^{(1)})\to\Udr_q(X_{\tilde n}^{(k)})&{\roman{if}}\ (X_{\tilde n}^{(k)},i)\neq(A_{2n}^{(2)},1)\cr
\Udr_q(A_2^{(2)})\to\Udr_q(X_{\tilde n}^{(k)})&{\roman{if}}\ (X_{\tilde n}^{(k)},i)=(A_{2n}^{(2)},1)\end{cases} .$$

\vskip .3 truecm
\rem{\rhxom}

i) $\bar\Omega({\Cal I}^+(HXL))={\Cal I}^-(HXL)$ and $\bar\Omega({\Cal I}^+(HX))={\Cal I}^-(HX)$;

ii) $\bar\Omega$ preserves relations $(HXL)^{\pm}$ and relations $(HX)^{\pm}$.

\vskip .3 truecm
\nota{\xxxtr}
Define relations $(XXD)$ $(XXE)$, $(XXH^{+})$ and $(XXH^{-})$ by:
$$[X_{i,r}^+,X_{j,s}^-]=0 \ \ ((i,r),(j,s)\in I_{\Z},\ i\neq j),\leqno{(XXD)}$$
$$[X_{i,r}^+,X_{i,-r}^-]=
{C^{r}k_i-C^{-r}k_i^{-1}\over q_i-q_i^{-1}} \ \ ((i,r)\in I_{\Z}),\leqno{(XXE)}$$
$$[X_{i,r}^+,X_{i,s}^-]=
{C^{-s}k_i\tilde H_{i,r+s}^+ \over q_i-q_i^{-1}} \ \ ((i,r),(i,s)\in I_{\Z},\ r+s>0),\leqno{(XXH^+)}$$
$$[X_{i,r}^+,X_{i,s}^-]=-
{C^{-r}\tilde H_{i,r+s}^-k_i^{-1}\over q_i-q_i^{-1}} \ \ ((i,r),(i,s)\in I_{\Z},\ r+s<0),\leqno{(XXH^-)}$$
\vskip 1 truecm
\rem{\rxxom}

i) ${\Cal I}(XX)={\Cal I}(XXD,XXE,XXH)$;

ii) $\bar\Omega({\Cal I}(XXD))={\Cal I}(XXD)$ and $\bar\Omega({\Cal I}(XXE))={\Cal I}(XXE)$;

iii)  $\bar\Omega({\Cal I}(XXH^+))={\Cal I}(XXH^-)$;

iv) $\bar\Omega$ preserves relations $(XX)$.

\vskip .3 truecm
\cor{\boto}
$\bar\Omega$ induces $\tilde\Omega:\tilde{\Cal U}^{Dr}_q\to\tilde{\Cal U}^{Dr}_q$.

\vskip .3 truecm
\rem{\xdost}

i) $\bar t_i({\Cal I}(XXD))={\Cal I}(XXD)$ and $\bar t_i({\Cal I}(XXE))={\Cal I}(XXE)$ $\forall i\in I_0$;

ii) ${\Cal I}(XXD,XXE)$ is the $\bar t_i^{\pm 1}$-stable ideal ($\forall i\in I_0$) generated by 
$$\Big\{[X_{i,0}^+,X_{j,0}^-]-\delta_{ij}{k_i-k_i^{-1}\over q_i-q_i^{-1}}|i,j\in I_0\Big\}.$$

\vskip .3 truecm

We want to show now that for all $i\in I_0$ $\bar t_i$ induces $\tilde t_i:\tilde{\Cal U}^{Dr}_q\to\tilde{\Cal U}^{Dr}_q$.
Since $\bar t_i$ commutes with $\bar\Omega$ remarks \rkhgr, \rhxom,i),
\rxxom,i) and iii) and \xdost,i) imply that it is enough to concentrate on ${\Cal I}(HXL^+)$, ${\Cal I}(XXH^+)$.

\vskip .3 truecm
\rem{\ttih}~
i) Remark that if $r+s>0$
$$(q_i-q_i^{-1})C^{s}k_i^{-1}[X_{i,r}^+,X_{i,s}^-]=\bar t_i^{s/\tilde d_i}((q_i-q_i^{-1})k_i^{-1}[X_{i,r+s}^+,X_{i,0}^-])=\bar t_i^{s/\tilde d_i}(\tilde H_{i,r+s}^+),$$
so that relations $(XXH)^+$ are equivalent to 
$$\bar t_i^s(\tilde H_{i,r}^+)=\tilde H_{i,r}^+\ \ \ \forall i\in I_0,\ r>0,\ s\in\Z;$$

ii) $\bar t_i^{\pm 1}({\Cal I}(XXH^+))={\Cal I}(XXH^+)$ $\forall i\in I_0$;

iii) ${\Cal I}^+(XXH)$ is the $\bar t_i^{\pm 1}$-stable ideal ($\forall i\in I_0$) generated by 
$$\{\bar t_i(\tilde H_{i,r}^+)-\tilde H_{i,r}^+|i\in I_0,\ r>0\}.$$
\vskip .3 truecm
\rem{\tind}~
Remark that for all $(i,r)\in I_{\Z}\setminus(I_0\times\{0\})$, $(j,s)\in I_{\Z}$ and $h\in I_0$
$$\bar t_h^{\pm 1}([H_{i,r},X_{j,s}^{+}]-b_{ijr}C^{r-|r|\over 2}X_{j,r+s}^{+})\!=\!
[\bar t_h^{\pm 1}(H_{i,r}),X_{j,s\mp \delta_{jh}\tilde d_j}^{+}]-b_{ijr}C^{r-|r|\over 2}X_{j,r+s\mp \delta_{jh}\tilde d_j}^{+}.$$
Then, thanks to remark \ttih\ and to the definition of $\bar t_i$ (see definition \drcmp.\gendef), we have that:

i) $\bar t_i^{\pm 1}({\Cal I}^+(HX))\subseteq{\Cal I}^+(XXH,HX)$;

ii) $\bar t_i^{\pm 1}({\Cal I}^+(HXL))\subseteq{\Cal I}^+(XXH,HXL)$;

iii) ${\Cal I}^+(XXH,HXL)$ is the $\bar t_i^{\pm 1}$-stable ideal ($\forall i\in I_0$) generated by 
$$\{\bar t_i(\tilde H_{i,r}^+)-\tilde H_{i,r}^+,[H_{i,s},X_{j,0}^{+}]-b_{ijs}C^{s-|s|\over 2}X_{j,s}^{+}|i\in I_0,\ r>0, \tilde d_i\leq|s|\leq\tilde d_{ij}\}.$$

\vskip .3 truecm
\cor{\ttibd}
i) For all $i\in I_0$ $\bar t_i$ induces $\tilde t_i:\tilde{\Cal U}^{Dr}_q\to\tilde{\Cal U}^{Dr}_q$;

ii) For all $i,j\in I_0$ $\tilde t_i(\tilde H_{j,r}^+)=\tilde H_{j,r}^+$ $\forall r\in\Z$ and $\tilde t_i(H_{j,r})=H_{j,r}$  $\forall r\in\Z\setminus\{0\}$.

\vskip .3 truecm

We come now to $\bar\Theta$ recalling that $\bar\Theta\bar\Omega=\bar\Omega\bar\Theta$ and
$\bar\Theta\bar t_i^{\pm 1}=\bar t_i^{\mp 1}\bar\Theta$ for all $i\in I_0$.

\vskip .3 truecm
\rem{\thxde}~
Notice that $[X_{i,0}^+,X_{j,0}^-]-\delta_{ij}{k_i-k_i^{-1}\over q_i-q_i^{-1}}$ is fixed by $\bar\Theta$; hence, thanks to remark \xdost,ii), ${\Cal I}(XXD,XXE)$ is $\bar\Theta$-stable.

\vskip .3 truecm
\rem{\txh}~
i) For all $i\in I_0$ and for all $r>0$ $$\bar\Theta(\tilde H_{i,r}^+)=\bar t_i^{r\over\tilde d_i}(\tilde H_{i,-r}^-)+(q_i-q_i^{-1})[X_{i,-r}^+X_{i,0}^-,k_i]C^{-r};$$

ii) for all $i\in I_0$ and for all $r> 0$ $\bar\Theta(\tilde H_{i,\pm r}^{\pm})-\tilde H_{i,\mp r}^{\mp}$ and $\bar\Theta(H_{i,\pm r})-H_{i,\mp r}$ lie in ${\Cal I}(KX,XXH^{\mp})$;

iii) for all $i\in I_0$ and for all $r>0$ 
$$\bar\Theta(\bar t_i(\tilde H_{i,r}^+)-\tilde H_{i,r}^+)=
\bar t_i^{-1}(\bar\Theta(\tilde H_{i,r}^+))-\bar\Theta(\tilde H_{i,r}^+)\in{\Cal I}(KX,XXH^{-});$$

iv) for all $i,j\in I_0$, $\tilde d_i\leq|r|\leq\tilde d_{ij}$, $s\in\Z$ 
$$\bar\Theta([H_{i,r},X_{j,s}^{+}]-b_{ijr}C^{r-|r|\over 2}X_{j,r+s}^{+})=$$
$$=-
[\bar\Theta(H_{i,r}),X_{j,-s}^{+}k_jC^{-s}]+b_{ijr}C^{r-|r|\over 2}X_{j,-(r+s)}^{+}k_jC^{-(r+s)}=$$
$$=-
([\bar\Theta(H_{i,r}),X_{j,-s}^{+}k_j]k_j^{-1}-b_{ijr}C^{-r-|r|\over 2}X_{j,-(r+s)}^{+})k_jC^{-s}$$
belongs to ${\Cal I}(KX,XXH,HXL^+)$;

Then:

v) $\bar\Theta({\Cal I}(XXH^+))\subseteq{\Cal I}(KX,XXH^-)$;

vi) $\bar\Theta({\Cal I}(HXL^+))\subseteq{\Cal I}(KX,XXH,HXL^+)$;

vii) ${\Cal I}(KX,XXH)$ and ${\Cal I}(KX,XXH,HXL^{\pm})$ are $\bar\Theta$-stable.

\vskip .3 truecm
\cor{\thibd}
i) $\bar\Theta$ induces $\tilde\Theta:\tilde{\Cal U}^{Dr}_q\to\tilde{\Cal U}^{Dr}_q$;

ii) For all $i\in I_0$ $\tilde\Theta(\tilde H_{i,r}^+)=\tilde H_{i,-r}^-$ $\forall r\in\Z$ and $\tilde\Theta(H_{i,r})=H_{i,-r}$  $\forall r\in\Z\setminus\{0\}$.

\vskip 1 truecm
\rem{\fth}
i) Let $f:Q_{0,+}\to\Z$ be defined by:
$$f(0)=0,\ \ \ f(\alpha+\alpha_i)=f(\alpha)+(\alpha|\alpha_i)\ \ \ \forall\alpha\in Q_{0,+},\ i\in I_0;$$
notice that $f$ is well defined, because $(\alpha|\alpha_i)+(\alpha+\alpha_i|\alpha_j)=(\alpha|\alpha_j)+(\alpha+\alpha_j|\alpha_i)$.

ii) $\forall X^+\in{\Cal F}_{q,\alpha+m\delta}^{+}$ and $\forall X^-\in{\Cal F}_{q,-\alpha+m\delta}^{-}$
(where $\alpha\in Q_{0,+}$, $m\in\Z$) we have that in $\tilde{\Cal U}_{q}^{Dr}$
$$\tilde\Theta\pi^{+}(X^+)=(-1)^hq^{f(\alpha)}\pi^{+}\Theta_{{\Cal F}}^{+}(X^+)k_{\alpha}C^{-m}$$ and 
$$\tilde\Theta\pi^{-}(X^-)=(-1)^hq^{-f(\alpha)}C^{-m}k_{-\alpha}\pi^{-}\Theta_{{\Cal F}}^{-}(X^-),$$
where $\pi^{\pm}\!:{\Cal F}_{q}^{\pm}\!\to\tilde{\Cal U}_{q}^{Dr}$ is the restriction to ${\Cal F}_{q}^{\pm}$ of the natural projection $\bar{\Cal U}_{q}^{Dr}\! \to\!\tilde{\Cal U}_{q}^{Dr}$ and $h=\sum_{i\in I_0}m_i$ if $\alpha=\sum_{i\in I_0}m_i\alpha_i$.

iii) In particular $\tilde\Theta\pi^{\pm}(X^{\pm})$ and $\pi^{\pm}\Theta_{{\Cal F}}^{\pm}(X^{\pm})$ are equal up to invertible elements of $\tilde{\Cal U}_{q}^{Dr}$.

\vskip .3 truecm
We present now some more remarks about generators and relations of $\tilde{\Cal U}_{q}^{Dr}$.

\vskip .3 truecm
For the next proposition see the analogous results for ${\Cal U}_q^{DJ}$, in \beck\ and \damcina. 

\vskip .3 truecm

\prop{\sop}
In $\tilde{\Cal{U}}_q^{Dr}$ we have ${\Cal{I}}(HX^{\pm})\subseteq{\Cal{I}}(XD^{\pm},X1^{\pm},X2^{\pm})$.
\dim
In order to avoid repetitive computations we use the behaviour of relations $(XD^{\pm})$, $(X1^{\pm})$ and $(X2^{\pm})$ under the action of $\tilde\Omega$, $\tilde\Theta$ and $\tilde t_i$ ($i\in I_0$), which is an independent result proved in remarks \dlem.\ggg\ and \dlem.\stabid: here it allows us to reduce to the study of $[H_{i,r},X_{j,0}^{+}]$ with $r>0$. Indeed: 
$$\tilde t_j^s([H_{i,r},X_{j,0}^{+}])=[H_{i,r},X_{j,-\tilde d_js}^{+}]\ \ \ {\roman{and}}\ \ \ 
\tilde t_j^s(X_{j,r}^{+})=X_{j,r-\tilde d_js}^{+}$$
$$({\roman{moreover}}\ \ \ [H_{i,r},X_{j,s}^{+}]=0\ \ \ {\roman{and}}\ \ \ b_{ijr}X_{j,r+s}^{+}=0\ \ \ {\roman{if}}\ \tilde d_j\not|s),$$
$$\tilde\Theta([H_{i,r},X_{j,s}^{+}])=-[H_{i,-r},X_{j,-s}^{+}]k_jC^{-s}\ \ \ {\roman{and}}\ \ \ 
\tilde\Theta(X_{j,r+s}^{+})=-X_{j,-r-s}^{+}k_jC^{-r-s},$$
$$\tilde\Omega([H_{i,r},X_{j,s}^{+}])=-[H_{i,-r},X_{j,-s}^{-}]\ \ \ {\roman{and}}\ \ \ 
\tilde\Omega(C^{{r- |r|\over 2}}X_{j,r+s}^{+})=C^{{-r+|r|\over 2}}X_{j,-r-s}^{-},$$
and of course $b_{ijr}=\tilde t_j(b_{ijr})=\tilde\Theta(b_{ijr})=\tilde\Omega (b_{ijr})=b_{ij-r}$; 

Given an element $x\in\tilde{\Cal{U}}_q^{Dr}$ define the operators ${}_lx$ and ${}_rx$ on $\tilde{\Cal{U}}_q^{Dr}$ respectively as the left and right multiplication by $x$; if we have elements $x_s\in\tilde{\Cal{U}}_q^{Dr}$ ($s\in\N$) set ${}_lx(u)=\sum_{s\in\N}{}_lx_{s}u^s$ and ${}_rx(u)=\sum_{s\in\N}{}_rx_{s}u^s$; notice that if $f:\tilde{\Cal{U}}_q^{Dr}\to\tilde{\Cal{U}}_q^{Dr}$ is such that $f(x_s)=x_s$ for all $s\in \N$ then 
${}_lx(u)$ and ${}_rx(u)$ commute with $f$. 

Let $i,j\in I_0$: we want to study 
$({}_l\tilde H_{i}^+(u)-{}_r\tilde H_{i}^+(u))(X_{j,s}^+)$ and deduce from it $({}_lH_{i}(u)-{}_rH_{i}(u))(X_{j,s}^+)$ (setting $H_{i,0}=0$). 
To this aim remark that
$${}_l\tilde H_{i}^+(u)=exp((q_i-q_i^{-1}){}_lH_{i}(u)),\ \ \ {}_r\tilde H_{i}^+(u)=exp((q_i-q_i^{-1}){}_rH_{i}(u))$$ and both commute with $\tilde t_j$.

The next computations are performed in $\tilde{\Cal{U}}_q^{Dr}/{\Cal I}(XD^+,X1^+,X2^+)$.

Remarking that ($r>\tilde d_{ij}$)
$$[\tilde H_{i,r}^+,X_{j,0}^+]=(q_i-q_i^{-1})k_i^{-1}[[X_{i,r}^+,X_{i,0}^-],X_{j,0}^+]_{q_i^{a_{ij}}}=$$
$$=(q_i-q_i^{-1})k_i^{-1}([X_{i,r}^+,[X_{i,0}^-,X_{j,0}^+]]_{q_i^{a_{ij}}}-[X_{i,0}^-,[X_{i,r}^+,X_{j,0}^+]_{q_i^{a_{ij}}}])=$$
$$=(q_i-q_i^{-1})k_i^{-1}(q_i^{a_{ij}}[\delta_{ij}{k_i-k_i^{-1}\over q_i-q_i^{-1}},X_{i,r}^+]_{q_i^{-a_{ij}}}+[[X_{i,r}^+,X_{j,0}^+]_{q_i^{a_{ij}}},X_{i,0}^-])=$$
$$=(q_i-q_i^{-1})k_i^{-1}(\delta_{ij}[2]_{q_i}k_iX_{i,r}^++[[X_{i,r}^+,X_{j,0}^+]_{q_i^{a_{ij}}},X_{i,0}^-]),$$
let us distinguish two cases: 

i) $(X_{\tilde n}^{(k)},i,j)\neq(A_{2 n}^{(2)},1,1)$: then, thanks to $(XD^{+})$, $(X1^{+})$ and $(HXL^{+})$, we have
$$[\tilde H_{i,r}^+,X_{j,0}^+]=(q_i-q_i^{-1})k_i^{-1}(\delta_{ij}[2]_{q_i}k_iX_{i,r}^+-[[X_{j,\tilde d_{ij}}^+,X_{i,r-\tilde d_{ij}}^+]_{q_i^{a_{ij}}},X_{i,0}^-])=$$
$$=(q_i-q_i^{-1})k_i^{-1}(\delta_{ij}[2]_{q_i}k_iX_{i,r}^+-[X_{j,\tilde d_{ij}}^+,[X_{i,r-\tilde d_{ij}}^+,X_{i,0}^-]]_{q_i^{a_{ij}}}+$$
$$+q_i^{a_{ij}}[X_{i,r-\tilde d_{ij}}^+,[X_{j,\tilde d_{ij}}^+,X_{i,0}^-]]_{q_i^{-a_{ij}}})=$$
$$=(q_i-q_i^{-1})\delta_{ij}([2]_{q_i}X_{i,r}^+-{1\over q_i-q_i^{-1}}[\tilde H_{i,\tilde d_{i}}^+,X_{i,r-\tilde d_{i}}^+])+$$
$$+(q_i-q_i^{-1})q_i^{a_{ij}}k_i^{-1}[[X_{i,r-\tilde d_{ij}}^+,X_{i,0}^-],X_{j,\tilde d_{ij}}^+]_{q_i^{-a_{ij}}}=$$
$$=(q_i-q_i^{-1})\delta_{ij}([2]_{q_i}X_{i,r}^+-[H_{i,\tilde d_{i}},X_{i,r-\tilde d_{i}}^+])+q_i^{a_{ij}}[\tilde H_{i,r-\tilde d_{ij}}^+,X_{j,\tilde d_{ij}}^+]_{q_i^{-2a_{ij}}}=$$
$$=q_i^{a_{ij}}\tilde H_{i,r-\tilde d_{ij}}^+X_{j,\tilde d_{ij}}^+-{q_i^{-a_{ij}}}X_{j,\tilde d_{ij}}^+\tilde H_{i,r-\tilde d_{ij}}^+;$$
hence, using again $(HXL^{+})$,
$$({}_l\tilde H_{i}^+(u)-{}_r\tilde H_{i}^+(u))(X_{j,0}^+)=(q_i^{a_{ij}}{}_l\tilde H_{i}^+(u)-q_i^{-a_{ij}}{}_r\tilde H_{i}^+(u))\tilde t_j^{{\tilde d_{ij}\over\tilde d_j}}u^{\tilde d_{ij}}(X_{j,0}^+),$$
or equivalently
$${}_l\tilde H_{i}^+(u)(1-q_i^{a_{ij}}\tilde t_j^{{\tilde d_{ij}\over\tilde d_j}}u^{\tilde d_{ij}})(X_{j,0}^+)=
{}_r\tilde H_{i}^+(u)(1-q_i^{-a_{ij}}\tilde t_j^{{\tilde d_{ij}\over\tilde d_j}}u^{\tilde d_{ij}})(X_{j,0}^+);$$
from this we get
$$(q_i-q_i^{-1})({}_lH_{i}(u)-{}_r H_{i}(u))(X_{j,0}^+)=$$
$$=log(1-q_i^{-a_{ij}}\tilde t_j^{{\tilde d_{ij}\over\tilde d_j}}u^{\tilde d_{ij}})-log(1-q_i^{a_{ij}}\tilde t_j^{{\tilde d_{ij}\over\tilde d_j}}u^{\tilde d_{ij}}),$$
that is 
$$[H_{i,r},X_{j,0}^+]=\begin{cases}  0&{\roman{if}}\ \tilde d_{ij}\not|r
\cr{q_i^{{ra_{ij}\over\tilde d_{ij}}}-q_i^{{-ra_{ij}\over\tilde d_{ij}}}\over {r\over \tilde d_{ij}}(q_i-q_i^{-1})}X_{j,r}^+=b_{ijr}X_{j,r}^+&{\roman{otherwise}}.\end{cases} $$

ii) $(X_{\tilde n}^{(k)},i,j)=(A_{2 n}^{(2)},1,1)$: the computations are a little more complicated than in case i), but substantially similar; we separate the cases $r=2$ and $r>2$ and, thanks to $(X2^{+})$ and $(HXL^{+})$, we get:
$$[\tilde H_{1,2}^+,X_{1,0}^+]=(q-q^{-1})k_1^{-1}([2]_{q}k_1X_{1,2}^++[[X_{1,2}^+,X_{1,0}^+]_{q^{2}},X_{1,0}^-])=$$
$$=(q-q^{-1})k_1^{-1}([2]_{q}k_1X_{1,2}^++(q^4-q^{-2})[(X_{1,1}^+)^2,X_{1,0}^-])=$$
$$=(q-q^{-1})k_1^{-1}([2]_{q}k_1X_{1,2}^++(q^4-q^{-2})(X_{1,1}^+[X_{1,1}^+,X_{1,0}^-]+[X_{1,1}^+,X_{1,0}^-]X_{1,1}^+))=$$
$$=(q-q^{-1})[2]_qX_{1,2}^++(q^4-q^{-2})(q^{-2}X_{1,1}^+\tilde H_{1,1}^++\tilde H_{1,1}^+X_{1,1}^+)=$$
$$=(q^2-q^{-2})X_{1,2}^++(q^4-q^{-2})\tilde H_{1,1}^+X_{1,1}^++(q^2-q^{-4})X_{1,1}^+\tilde H_{1,1}^+,$$
hence, for all $s\in\Z$,
$$[\tilde H_{1,2}^+,X_{1,s}^+]=(q^2-q^{-2})X_{1,s+2}^++(q^4-q^{-2})[\tilde H_{1,1}^+,X_{1,s+1}^+]_{-q^{-2}};$$
$r>2$:
$$[\tilde H_{1,r}^+,X_{1,0}^+]=(q-q^{-1})k_1^{-1}([2]_{q}k_1X_{1,r}^++[[X_{1,r}^+,X_{1,0}^+]_{q^{2}},X_{1,0}^-])=$$
$$=(q-q^{-1})k_1^{-1}([2]_{q}k_1X_{1,r}^++
[(q^4-q^{-2})[X_{1,r-1}^+,X_{1,1}^+]_{-1}-
[X_{1,2}^+,X_{1,r-2}^+]_{q^{2}},X_{1,0}^-])=$$
$$=(q-q^{-1})k_1^{-1}([2]_{q}k_1X_{1,r}^++$$
$$+(q^4-q^{-2})[X_{1,r-1}^+,[X_{1,1}^+,X_{1,0}^-]]_{-1}+(q^4-q^{-2})[X_{1,1}^+,[X_{1,r-1}^+,X_{1,0}^-]]_{-1}+$$
$$-
[X_{1,2}^+,[X_{1,r-2}^+,X_{1,0}^-]]_{q^{2}}+q^{2}[X_{1,r-2}^+,[X_{1,2}^+,X_{1,0}^-]]_{q^{-2}})=$$
$$=(q^2-q^{-2})X_{1,r}^++$$
$$+(q^4-q^{-2})[\tilde H_{1,1}^+,X_{1,r-1}^+]_{-q^{-2}}+(q^4-q^{-2})[\tilde H_{1,r-1}^+,X_{1,1}^+]_{-q^{-2}}+$$
$$+q^{2}
[\tilde H_{1,r-2}^+,X_{1,2}^+]_{q^{-4}}-[\tilde H_{1,2}^+,X_{1,r-2}^+]=$$
$$=(q^4-q^{-2})\tilde H_{1,r-1}^+X_{1,1}^++(q^2-q^{-4})X_{1,1}^+\tilde H_{1,r-1}^++$$
$$+q^{2}
\tilde H_{1,r-2}^+X_{1,2}^+-q^{-2}X_{1,2}^+\tilde H_{1,r-2}^+;$$
this implies, using again $(HXL^{+})$, that 
$$({}_l\tilde H_{1}^+(u)-{}_r\tilde H_{1}^+(u))(X_{1,0}^+)=$$
$$=((q^4-q^{-2}){}_l\tilde H_{1}^+(u)\tilde t_1u+(q^2-q^{-4}){}_r\tilde H_{1}^+(u))\tilde t_1u+$$
$$+q^{2}{}_l\tilde H_{1}^+(u)(\tilde t_1u)^2-q^{-2}{}_r\tilde H_{1}^+(u))(\tilde t_1u)^2)(X_{1,0}^+),$$
or equivalently
$${}_l\tilde H_{1}^+(u)(1-q^{4}\tilde t_1u)(1+q^{-2}\tilde t_1u)(X_{1,0}^+)=$$
$$=
{}_r\tilde H_{1}^+(u)(1-q^{-4}\tilde t_1u)(1+q^{2}\tilde t_1u)(X_{1,0}^+);$$
from this we get
$$(q-q^{-1})({}_lH_{1}(u)-{}_r H_{1}(u))(X_{1,0}^+)=$$
$$=(log(1-q^{-4}(\tilde t_1u))+log(1+q^{2}(\tilde t_1u))-log(1-q^{4}(\tilde t_1u))-log(1+q^{-2}(\tilde t_1u)))(X_{1,0}^+),$$
that is 
$$[H_{1,r},X_{1,0}^+]={-q^{-4r}+(-1)^{r-1}q^{2r}+q^{4r}-(-1)^{r-1}q^{-2r}\over r(q-q^{-1})}X_{1,r}^+=b_{11r}X_{1,r}^+.$$

\vskip .3truecm

\prop{\boh}
In $\tilde{\Cal U}_{q}^{Dr}$ we have
${\Cal I}(HH)\subseteq{\Cal I}(HX)$. 
\dim
Thanks to remark \greh, to the fact that 
$$[H_{i,r},H_{j,s}]=-[H_{j,s},H_{i,r}]=\tilde\Omega[H_{j,-s},H_{i,-r}]$$ and to the definition of $b_{ijr}$, it is enough to prove that in $\tilde{\Cal U}_q^{Dr}/{\Cal I}(HX)$
$$[H_{i,r},\tilde H_{j,s}^+]=\delta_{r+s,0}b_{ijr}(C^r-C^{-r})\ {\roman{if}}\ |r|\geq s>0.$$ This is an easy computation:
$$[H_{i,r},\tilde H_{j,s}^+]=(q_j-q_j^{-1})k_j^{-1}[H_{i,r},[X_{j,s}^+,X_{j,0}^-]]=$$
$$=(q_j-q_j^{-1})k_j^{-1}([[H_{i,r},X_{j,s}^+],X_{j,0}^-]-[[H_{i,r},X_{j,0}^-],X_{j,s}^+])=$$
$$=(q_j-q_j^{-1})b_{ijr}k_j^{-1}(C^{{r-|r|\over 2}}[X_{j,r+s}^+,X_{j,0}^-]-C^{{r+|r|\over 2}}[X_{j,s}^+,X_{j,r}^-])=$$
$$=b_{ijr}k_j^{-1}(C^{{r-|r|\over 2}}k_j\tilde H_{j,r+s}^+-C^{{-r-|r|\over 2}-s}k_j^{-1}\tilde H_{j,r+s}^-+$$
$$-C^{{-r
+|r|\over 2}}k_j\tilde H_{j,r+s}^++C^{{r+|r|\over 2}-s}k_j^{-1}\tilde H_{j,r+s}^-)=$$
$$=b_{ijr}k_j^{-1}((C^{{r-|r|\over 2}}-C^{{-r
+|r|\over 2}})k_j\tilde H_{j,r+s}^++(C^{{r+|r|\over 2}}-C^{{-r-|r|\over 2}})C^{-s}k_j^{-1}\tilde H_{j,r+s}^-)=$$
$$=b_{ijr}k_j^{-1}(C^{{r-|r|\over 2}}-C^{{-r
+|r|\over 2}})k_j\tilde H_{j,r+s}^+=$$
$$=\delta_{r+s,0}b_{ijr}(C^r-C^{-r}).$$

\vskip .3 truecm

\cor{\dnd}
i) Relations $(ZH)$ and $(KH)$ are redundant.

ii) In $\tilde{\Cal U}^{Dr}_q$ relations $(HX^{\pm})$ depend on $(XD^{\pm})$, $(X1^{\pm})$ and $(X2^{\pm})$ and relations $(HH)$ depend on $(XD)$, $(X1)$ and $(X2)$.

iii) $\Udr_q(X_{\tilde n}^{(k)})$ is the quotient of $\tilde{\Cal U}^{Dr}_q(X_{\tilde n}^{(k)})$  by the ideal generated by the rela\-tions $(XD^{\pm})$-$(S3^{\pm})$.

\vskip .3truecm
\rem{\solocom}~
It is worth remarking that corollary \dnd,ii) allows us to reduce the relations $(HX^{\pm})$ and $(HH)$ to relations involving just the $X_{i,r}^{\pm}$'s, without using the $H_{i,r}$'s whose connection with the $\tilde H_{i,r}^{\pm}$'s (these last can be expressed in terms of commutators between 
the $X_{i,r}^{+}$'s and the $X_{i,r}^{-}$'s, see remark \drcmp.\tiht) is complicated to handle. 
Indeed relations $(HXL^{\pm})$ can be translated as follows: 

i) If $\tilde d_i\leq |r|<\tilde d_{ij}$ then $[H_{i,r},X_{j,s}^{\pm}]=0$, that is $X_{j,s}^{\pm}$ commutes with the subalgebra generated by $\{H_{i,r}|\tilde d_i\leq |r|<\tilde d_{ij}\}$, which is $\{\tilde H_{i,r}^{\pm}|\tilde d_i\leq |r|<\tilde d_{ij}\}$ (see remark \drcmp.\greh); hence these relations can be rewritten as 
$$[\tilde H_{i,r}^{\pm},X_{j,s}^{+}]=0\ \ {\roman {and}}\ \ [\tilde H_{i,r}^{\pm},X_{j,s}^{-}]=0\ \ {\roman {if}}\ |r|<\tilde d_{ij};$$

ii) If $|r|=\tilde d_{ij}$ and $\pm r>0$ then $\tilde H_{i,r}^{\pm}\mp(q_i-q_i^{-1})H_{i,r}$ commutes with $
X_{j,s}^{\pm}$, by i) and remark \drcmp.\greh, hence in the relations 
$[H_{i,r},X_{j,s}^{\pm}]=\pm b_{ij\tilde d_{ij}}C^{{r\mp|r|\over 2}}X_{j,r+s}^{\pm}$ we can replace $H_{i,r}$ with  $\pm{\tilde H_{i,r}^{\pm}\over (q_i-q_i^{-1})}$. 

Then relations $(HXL^{\pm})$ are equivalent to 

$$[[X_{i,\pm r}^{\pm},X_{i,0}^{\mp}],X_{j,s}^{\pm}]_{q_i^{a_{ij}}}=b_{ijr}k_i^{\pm 1}X_{j,s\pm r}^{\pm}$$
and 
$$[[X_{i,\pm r}^{\pm},X_{i,0}^{\mp}],X_{j,s}^{\mp}]_{q_i^{-a_{ij}}}=-b_{ijr}C^{\pm r}k_i^{\pm 1}X_{j,s\pm r}^{\mp}$$
with $0<r\leq\tilde d_{ij}$.

Remark also that among the relations defining $\Udr_q$ there are no other relations involving the $H_{i,r}$'s.

\vskip .3truecm

\rem{\grem}~
Remark that $\forall i\in I_0$ $\{C^{\pm\tilde d_i}, k_i^{\pm 1}, X_{i,r}^{\pm},H_{i,s}|\tilde d_i|r,s;\ s\neq 0\}$ generates $Im(\tilde\phi_i)\subseteq\tilde{\Cal U}_q^{Dr}$ over $\C(q_i)$. Then the following sets generate $Im(\tilde\phi_i)$ (hence $Im(\phi_i)\subseteq{\Cal U}_q^{Dr}$) over $\C(q_i)$: 

i) $\{C^{\pm\tilde d_i}, k_i^{\pm 1}, X_{i,r}^{\pm}|\tilde d_i|r\}$;

ii) $\{C^{\pm\tilde d_i}, k_i^{\pm 1}, X_{i,0}^{\pm}, H_{i,\pm\tilde d_i}\}$;

iii) $\{C^{\pm\tilde d_i}, k_i^{\pm 1}, X_{i,0}^{\pm}, X_{i,\mp\tilde d_i}^{\pm}\}$.

Moreover

iv) $\{C^{\pm 1}, k_i^{\pm 1}, X_{i,0}^{\pm}, X_{i_0,\mp 1}^{\pm}|i\in I_0\}$ (where $i_0$ is any fixed element of $I_0$ with the property $\tilde d_{i_0}=1$) generates $\tilde{\Cal U}_q^{Dr}$ (hence ${\Cal U}_q^{Dr}$) over $\C(q)$.

\dim

i) See remarks \drcmp.\greh\ and \drcmp.\tiht;

ii) follows from i) by induction on $|r|$, using that $$\forall r\in\Z\ \ \ [H_{i,\pm\tilde d_i},X_{i,r}^{+}]=
b_{ii\tilde d_i}C^{{\pm 1-1\over 2}\tilde d_i}X_{i,r\pm\tilde d_i}^{+}$$ and applying $\tilde\Omega$ (the set $\{C^{\pm\tilde d_i}, k_i^{\pm 1}, X_{i,0}^{\pm}, H_{i,\pm\tilde d_i}\}$ is $\tilde\Omega$-stable);

iii) is an immediate consequence of ii) and of the fact that 
$[X_{i,\tilde d_i}^+,X_{i,0}^-]=k_iH_{i,\tilde d_i}$, again applying $\tilde\Omega$;

iv) $\forall i\in I_0$ there exists a sequence of different indices $i_0,i_1,...,i_l=i$ in $I_0$ such that $a_{i_{h-1}i_h}<0$ and $\tilde d_{i_{h-1}}|\tilde d_{i_h}$ $\forall h=1,...,l$. 

We prove by induction on $h$ that $Im(\phi_{i_h})$ is contained in the $\C(q)$-subalgebra of $\tilde{\Cal U}_q^{Dr}$ generated by $\{C^{\pm 1}, k_i^{\pm 1}, X_{i,0}^{\pm}, X_{i_0,\mp 1}^{\pm}|i\in I_0\}$, the claim for $h=0$ being iii). 

For $h>0$ it is again enough to use iii), remarking  that $$[H_{i_{h-1},-\tilde d_{i_h}},X_{i_h,0}^{+}]=
b_{i_{h-1}i_h\tilde d_{i_h}}C^{-\tilde d_{i_h}}X_{i_h,-\tilde d_{i_h}}^+\neq 0$$ and applying $\tilde\Omega$.

\vskip .5truecm
{\bf{\dautemb.\ ${\Cal U}_q$: (anti)AUTOMORPHISMS and RELATIONS.}}
\vskip .5 truecm

The main point of this section is to describe in some details how the (anti)auto\-mor\-phisms $\tilde\Omega$, $\tilde\Theta$ and $\tilde t_i$ ($i\in I_0$) act on the generators of the ideal of $\tilde{\Cal U}_q$ defining ${\Cal U}_q$. As a corollary $\tilde\Omega$, $\tilde\Theta$ and $\tilde t_i$ ($i\in I_0$) induce analogous $\Omega$, $\Theta$ and $t_i$ ($i\in I_0$) on ${\Cal U}_q$. But the important consequence of this analysis (together with the study of the commutation with the elements  $H_{i,r}$'s) is the reduction of the huge amount of relations $(XD^{\pm})$-$(S3^{\pm})$ to relations involving only the positive $X_{i,r}^{+}$ (which is obvious and well known) and, which is new, to the analogous relations with ``constant parameters'' (see section \dlem). Lemmas \dlem.\dvssdz\ and \dlem.\dvss\ are the fundamental tool of this paper, which makes possible and easy the computations of the following sections, leading to theorems \cnat.\tldg\ and \serrel.\ssrr.

\vskip .3 truecm
\nota{\zl}
Let $l\in\N$; then:

i) $\1=\1_l=(1,...,1)\in\Z^l$;

ii) $\{e_1,...,e_l\}$ is the canonical basis of $\Z^l$

iii) for all $r=(r_1,...,r_l)\in\Z^l$ $\bar r\in\Z^l$ denotes $\bar r=(r_l,...,r_1)$.

\vskip .3truecm
\ddefi{\omdef}
$\Omega:\Udr_q\to\Udr_q$ is the $\C$-anti-linear anti-homomorphism induced by $\tilde\Omega$ (and by $\bar\Omega$, see definition \drcmp.\gendef\ and remark \drcmp.\genbdef,vii)
), that is the $\C$-anti-linear anti-homomorphism defined on the generators by $$q\mapsto q^{-1},\ \ C^{\pm 1}\mapsto C^{\mp 1},\ \ k_i^{\pm 1}\mapsto k_i^{\mp 1},\ \ X_{i,r}^{\pm}\mapsto X_{i,-r}^{\mp},\ \ H_{i,r}\mapsto H_{i,-r}.$$
\vskip .3truecm

\rem{\ombdef}~
$\Omega$ is a well-defined involution of $\Udr_q$. Indeed 
$$
\tilde\Omega(M_{(2)}^{\pm}((i,r),(j,s)))=-q_i^{-a_{ij}}M_{(2)}^{\mp}((i,-r),(j,-s)),$$
$$
\tilde\Omega(M_{i}^{\pm}(r))=-q_i^{-2}M_{i}^{\mp}(-r),$$
$$\tilde\Omega(M_{(2,2)}^{\pm}(r))=-q^{-2}M_{(2,2)}^{\mp}(-r),$$
$$\tilde\Omega(M_{(3)}^{\varepsilon,\pm}(r))=q^{-6\varepsilon}M_{(3)}^{\varepsilon,\mp}(-r),$$
$$\tilde\Omega(M_{i,j;l;a}^{\pm}(r;s))=(-1)^lq_i^{l(a_{ij}+a(l-1))}M_{i,j;l;a}^{\mp}(-r;-s),$$
$$\tilde\Omega(X_{i,j;l;a}^{\pm}(r;s))=(-1)^lX_{i,j;l;a}^{\mp}(-\bar r;-s),$$
$$\tilde\Omega(X_{[k]}^{\pm}(r;s))=X_{[k]}^{\mp}(-\bar r;-s),$$
$$\tilde\Omega(M_{[2]}^{\pm}(r;s))=q^{-2}M_{[2]}^{\mp}(-r;-s),$$
$$\tilde\Omega(M_{[3]}^{\pm}(r;s))=q^{-4}M_{[3]}^{\mp}(-r;-s).$$

\vskip .3truecm

\ddefi{\xidef}
$\Theta:\Udr_q\to\Udr_q$ is the $\C$-anti-linear homomorphism induced by $\tilde\Theta$ (and by $\bar\Theta$, see definition \drcmp.\gendef\ and corollary \darl.\thibd,ii)), that is the $\C$-anti-linear homomorphism defined on the generators by $$q\mapsto q^{-1},\ \ C^{\pm 1}\mapsto C^{\pm 1},\ \ k_i^{\pm 1}\mapsto k_i^{\mp 1},$$
$$ X_{i,r}^+\mapsto -X_{i,-r}^+k_iC^{-r},\ \ X_{i,r}^-\mapsto -k_i^{-1}C^{-r}X_{i,-r}^-,\ \ H_{i,r}\mapsto H_{i,-r}.$$
\vskip .3truecm

\rem{\xibdef}~
$\Theta$ is a well-defined involution of $\Udr_q$. Indeed
$$
\Theta_{{\Cal F}}^{\pm}(M_{(2)}^{\pm}((i,r),(j,s)))=-q_i^{-a_{ij}}M_{(2)}^{\pm}((i,-r\mp\tilde d_{ij}),(j,-s\mp\tilde d_{ij})),$$
$$
\Theta_{{\Cal F}}^{\pm}(M_{i}^{\pm}(r))=-q_i^{-2}M_{i}^{\pm}(-\bar r\mp\tilde d_{ij}\1),$$
$$
\Theta_{{\Cal F}}^{\pm}(M_{(2,2)}^{\pm}(r))=-q^{-2}M_{(2,2)}^{\pm}(-\bar r\mp 2\1),$$
$$\Theta_{{\Cal F}}^{\pm}(M_{(3)}^{\varepsilon,\pm}(r))=M_{(3)}^{-\varepsilon,\pm}(-r),$$
$$\Theta_{{\Cal F}}^{\pm}(X_{i,j;l;a}^{\pm}(r;s))=X_{i,j;l;a}^{\pm}(-r;-s),$$
$$\Theta_{{\Cal F}}^{\pm}(X_{[k]}^{\pm}(r;s))=X_{[k]}^{\pm}(- r\mp (k-1)\1;-s).$$

\vskip .3truecm

\ddefi{\tidef}
For all $i\in I_0$ $\ti:\Udr_q\to\Udr_q$ is  the $\C(q)$-homomorphism induced by $\tilde t_i$ (and by $\bar t_i$, see definition \drcmp.\gendef\ and corollary \darl.\ttibd,ii)), that is the $\C(q)$-homomorphism defined on the generators by $$C^{\pm 1}\mapsto C^{\pm 1},\ \ \ k_j^{\pm 1}\mapsto (k_jC^{-\delta_{ij}\tilde d_i})^{\pm 1},\ \ \ X_{j,r}^{\pm}\mapsto X_{j,r\mp\d_{ij}\tilde d_i}^{\pm},\ \ \ H_{j,r}\mapsto H_{j,r}.$$
\vskip .3truecm
\rem{\tibdef}~
It is immediate to check that the $\ti$'s are well-defined automorphisms of $\Udr_q$. Indeed
$$\tilde\ti(M_{(2)}^{\pm}((j,r),(h,s)))=
M_{(2)}^{\pm}((j,r\mp\delta_{ij}\tilde d_i),(h,s\mp\delta_{ih}\tilde d_i)),$$
$$\tilde\ti(M_{j}^{\pm}(r))=M_{j}^{\pm}(r\mp\d_{ij}\tilde d_i\uno),$$
$$\tilde\ti(M_{(2,2)}^{\pm}(r))=M_{(2,2)}^{\pm}(r\mp\d_{i1}\uno),$$
$$\tilde\ti(M_{(3)}^{\varepsilon,\pm}(r))=M_{(3)}^{\varepsilon,\pm}(r\mp\d_{i1}\uno),$$
$$\tilde\ti(M_{j,h;l;a}^{\pm}(r;s))=M_{j,h;l;a}^{\pm}(r\mp\delta_{ij}\tilde d_i\1;s\mp\delta_{ih}\tilde d_i),$$
$$\tilde\ti(X_{j,h;l;a}^{\pm}(r;s))=X_{j,h;l;a}^{\pm}(r\mp\delta_{ij}\tilde d_i\1;s\mp\delta_{ih}\tilde d_i),$$
$$\tilde\ti(X_{[k]}^{\pm}(r;s))=X_{[k]}^{\pm}(r\mp\delta_{ij}\1;s\mp\delta_{ih}\tilde d),$$
$$\tilde\ti(M_{[k]}^{\pm}(r;s))=M_{[k]}^{\pm}(r\mp\delta_{ij}\1;s\mp\delta_{ih}k);$$
in the last two identities $j,h\in I_0$ are such that $a_{jh}=-k$.

\vskip .3truecm

\rem{\aucom}~
Of course (see remark \drcmp.\genbdef, v) and vi)) $$\Theta\Omega=\Omega\Theta,\ \ t_i\Omega=\Omega t_i,\ \ t_i\Theta=\Theta t_i^{-1}\ \ {\roman{and}}\ \  t_i t_j= t_j t_i\ \ \ \forall i,j\in I_0$$ as maps of $\Udr_q(X_{\tilde n}^{(k)})$ into itself, and, $\forall i\in I_0$, $$\Omega\phi_i=\phi_i\Omega,\ \ \Theta\phi_i=\phi_i\Theta,\ \ t_i\phi_i=\phi_i t_1,\ \  t_j\phi_i=\phi_i\ \ \forall j\in I_0\setminus\{i\}$$  as maps from $\Udr_q(A_1^{(1)})$ to $\Udr_q(X_{\tilde n}^{(k)})$ if $(X_{\tilde n}^{(k)},i)\neq(A_{2n}^{(2)},1)$, from $\Udr_q(A_2^{(2)})$ to $\Udr_q(X_{\tilde n}^{(k)})$ if $(X_{\tilde n}^{(k)},i)=(A_{2n}^{(2)},1)$;

for all $\alpha\in Q$, $\beta\in Q_0$, $m\in\Z$, $i\in I_0$ we have
$$\Omega({\Cal U}_{q,\alpha}^{Dr})={\Cal U}_{q,-\alpha}^{Dr},\ \ \ 
\Theta({\Cal U}_{q,\beta+m\delta}^{Dr})={\Cal U}_{q,\beta-m\delta}^{Dr},\ \ \ 
t_i({\Cal U}_{q,\alpha}^{Dr})={\Cal U}_{q,\lambda_i(\alpha)}^{Dr}.$$

\vskip 2 truecm
{\bf{\dlem.\ REDUCTION to 
RELATIONS with CONSTANT PARAMETER.}}
\vskip .5 truecm

We shall now apply the structures introduced until now to the analysis of the relations defining ${\Cal U}_q^{Dr}$.
\vskip .3 truecm
\nota{\relid}
Let $(R)$ be relations as in notation \drmr.\aaa\ and 
define the following ideals of ${\Cal U}$:
$${\Cal I}_{cte}(R)=(S_{\zeta}(r\1_l,s)|\zeta\in{\Cal Z}, r\in\Z,s\in\Z^{\tilde l}),$$
$${\Cal I}_{0}(R)=(S_{\zeta}(\underline 0)|\zeta\in{\Cal Z}),$$
where $\underline 0\in\Z^{l+\tilde l}$;
more precisely, if $\#{\Cal Z}=1$ and ${\Cal Z}=\{\zeta\}$,
given $r\in\Z$ and $s\in\Z^{\tilde l}$, let 
$${\Cal I}_{(r,s)}(R)(={\Cal I}_{(r,s)}(R_{\zeta}))=(S_{\zeta}(r\1_l,s)).$$
If $(^{(h)}\!R)$ ($h=1,...,m$) are as in notation \drmr.\aaa,iii), define
$${\Cal I}_{cte}(^{(1)}\!R,...,^{(m)}\!R)=({\Cal I}_{cte}(^{(1)}\!R),...,{\Cal I}_{cte}(^{(m)}\!R))$$
and
$${\Cal I}_{0}(^{(1)}\!R,...,^{(m)}\!R)=({\Cal I}_{0}(^{(1)}\!R),...,{\Cal I}_{0}(^{(m)}\!R)).$$
Finally, if moreover each $(^{(h)}\!R)=(^{(h)}\!R^+)\cup(^{(h)}\!R^-)$ where $(^{(h)}\!R^{\pm})$ is as in notation \drmr.\aaa,iv), we shall also use the notation 
$${\Cal I}_{*}^{\pm}(^{(1)}\!R,...,^{(m)}\!R)=({\Cal I}_{*}(^{(1)}\!R^{\pm}),...,{\Cal I}_{*}(^{(m)}\!R^{\pm}))$$
where 
$*\in\{\emptyset,cte,0\}$.

\vskip .3 truecm
\rem{\rgen}~
With the notations fixed in notations \drmr.\aaa\ and \relid\ we have: 

i) ${\Cal I}_0(R)\subseteq{\Cal I}_{cte}(R)\subseteq{\Cal I}(R)$;

ii) ${\Cal I}_{*}(R)=({\Cal I}_{*}(R_{\zeta})|\zeta\in{\Cal Z})$ $\forall *\in\{\emptyset,cte,0\}$;

iii) for all $\zeta\in{\Cal Z}$ ${\Cal I}_{cte}(R_{\zeta})=({\Cal I}_{(r,s)}(R_{\zeta})|(r,s)\in\Z^{1+\tilde l})$;

iv) $l= 1\Rightarrow{\Cal I}_{cte}(R)={\Cal I}(R)$.

\vskip .3 truecm
\rem{\bbb}~
Let $(R^{\pm})$ be relations in $\tilde{\Cal U}^{Dr}_q$ as in notation \drmr.\aaa,iv) and suppose that for all $\zeta\in{\Cal Z}$ and for all $r\in\Z^l$, 
$s\in\Z^{\tilde l}$ there exists an invertible element $u_{\zeta,r,s}$ of $\tilde{\Cal U}_q^{Dr}$ such that $\tilde\Omega(S_{\zeta}^{+}(r,s))=u_{\zeta,r,s}S_{\zeta}^{-}(-r,-s)$ (notice that if $(R^{\pm})$ has this property then $(R_{\zeta}^{\pm})$ has the same property). With the notations fixed in notations \drmr.\aaa\ and \relid\ we have: 

i) $\tilde\Omega({\Cal I}_{r,s}^+(R_{\zeta}))={\Cal I}_{-r,-s}^-(R_{\zeta})$ $\forall \zeta\in{\Cal Z}$, $(r,s)\in\Z^{1+\tilde l}$;

ii) $\tilde\Omega({\Cal I}^+(R))={\Cal I}^-(R)$, $\tilde\Omega({\Cal I}_{cte}^+(R))={\Cal I}_{cte}^-(R)$ and  $\tilde\Omega({\Cal I}_{0}^+(R))={\Cal I}_{0}^-(R)$;

iii) ${\Cal I}(R)$, ${\Cal I}_{cte}(R)$ and ${\Cal I}_{0}(R)$ are the $\tilde\Omega$-stable ideals generated respectively by ${\Cal I}^+(R)$, ${\Cal I}_{cte}^+(R)$ and ${\Cal I}_{0}^+(R)$;
\vskip .3 truecm
\rem{\ccc}~
Let $(R)$ be relations in $\tilde{\Cal U}^{Dr}_q$ as in notation \drmr.\aaa\ and suppose that for all $\zeta\in{\Cal Z}$ there exist $i,j\in I_0$ such that for all $r\in\Z^l$, $s\in\Z^{\tilde l}$ we have:
$$\tilde t_i(S_{\zeta}(r,s))=S_{\zeta}(r-\1_l,s),$$
$$\tilde t_j(S_{\zeta}(r,s))=S_{\zeta}(r,s-\1_{\tilde l})$$
$$\tilde t_h(S_{\zeta}(r,s))=S_{\zeta}(r,s)\ \ \forall h\neq i,j.$$
Then:

i) ${\Cal I}(R_{\zeta})$ and  ${\Cal I}_{cte}(R_{\zeta})$ are $\tilde t_i^{\pm 1}$-stable for all $ \zeta\in{\Cal Z}$ and $i\in I_0$;

ii) $\forall \zeta\in{\Cal Z}$ and $\forall (r,s),(\tilde r,\tilde s)\in\Z^{1+\tilde l}$ there exists 
$\tilde t\in<\tilde t_i|i\in I_0>$ such that 
${\Cal I}_{(\tilde r,\tilde s)}(R_{\zeta})=\tilde t({\Cal I}_{ (r,s)}(R_{\zeta}))$;

iii) $\forall \zeta\in{\Cal Z}$ and $\forall (r,s)
\in\Z^{1+\tilde l}$ ${\Cal I}_{cte}(R_{\zeta})$ is the $\tilde t_i^{\pm 1}$-stable (for all $i\in I_0$) ideal generated by ${\Cal I}_{(r,s)}(R_{\zeta})$;

iv) ${\Cal I}_{cte}(R)$ is the $\tilde t_i^{\pm 1}$-stable ideal (for all $i\in I_0$) generated by ${\Cal I}_0(R)$.

\vskip .3 truecm
\cor{\ddd}
i) If $(R)$ satisfies the conditions of remarks \bbb\ and \ccc\ then ${\Cal I}_{cte}(R)$ is the $\tilde \Omega$-stable and $\tilde t_i$-stable (for all $i\in I_0$) ideal generated by ${\Cal I}_0^{+}(R)$. More precisely $\forall\zeta\in{\Cal{Z}}$ ${\Cal I}_{cte}(R_{\zeta})$ is the $\tilde \Omega$-stable and $\tilde t_i$-stable (for all $i\in I_0$) ideal generated by ${\Cal I}_0^{+}(R_{\zeta})$.

ii) Let $(^{(1)}\!R^{\pm})$ and $(^{(2)}\!R^{\pm})$ be as in remarks \bbb\ and \ccc, and suppose that 
${\Cal I}_0^{+}(^{(1)}\!R)\subseteq{\Cal I}_{cte}^{+}(^{(2)}\!R)$; then ${\Cal I}_{cte}(^{(1)}\!R)\subseteq{\Cal I}_{cte}(^{(2)}\!R)$. 

More precisely ${\Cal I}_{cte}^{+}(^{(1)}\!R)\subseteq{\Cal I}_{cte}^{+}(^{(2)}\!R)$ if and only if for all $\zeta\in^{(1)}\!\!\!{\Cal Z}$ there exists $(r,s)\in\Z\times\Z^{\tilde l_h}$ such that
${\Cal I}_{(r,s)}^{+}(^{(1)}\!R_{\zeta})\subseteq{\Cal I}_{cte}^{+}(^{(2)}\!R)$, and if this is the case we have also
${\Cal I}_{cte}^{-}(^{(1)}\!R)\subseteq{\Cal I}_{cte}^{-}(^{(2)}\!R)$.
\vskip .3 truecm
\rem{\eee}~
With the notations fixed in notation \drmr.\aaa\ suppose that 
$$\sigma.S_{\zeta}(r,s)=S_{\zeta}(r,s)\ \ \ \forall\zeta\in{\Cal Z},\ r\in\Z^l,\ s\in\Z^{\tilde l},\ \sigma\in\sy_l,$$
where $\sigma.S_{\zeta}(r,s)=S_{\zeta}(\sigma.r,s)$, see notation \drrl.\didi,iv). 

This condition is equivalent to the existence of elements $N_{\zeta}(r,s)\in\tilde{\Cal U}_q^{Dr}$ such that 
$$S_{\zeta}(r,s)=\sum_{\sigma\in\sy_l}\sigma.N_{\zeta}(r,s).$$
Notice that in general the elements $N_{\zeta}(r,s)$ ($r\in\Z^l,\ s\in\Z^{\tilde l}$) are not uniquely determined by the 
$S_{\zeta}(r,s)$'s. 

But $N_{\zeta}(r\1_l,s)={1\over l!}S_{\zeta}(r\1_l,s)$ for all $(r,s)\in\Z\times\Z^{\tilde l}$.

\vskip .3 truecm
\rem{\fff}~
In the hypotheses of remark \eee\ suppose that:

i) ${\Cal Z}\subseteq\begin{cases}  I_0&{\roman {if}}\ \tilde l=0\cr
\{(i,j)\in I_0\times I_0|i\neq j\}&{\roman {if}}\ \tilde l=1;\end{cases} $

ii) if $\tilde l=0$ and $i\in{\Cal Z}\subseteq I_0$ there exists $c_{p,\sigma}\in\C(q)$ ($p\in\Z^l$, $\sigma\in\sy_l$) such that for all $r\in\Z^l$
$$N_i(r)=\sum_{{p\in\Z^l\atop\sigma\in\sy_l}}c_{p,\sigma}X_{i,\tilde d_i(r_{\sigma(1)}+p_1)}^+\cdot...\cdot X_{i,\tilde d_i(r_{\sigma(l)}+p_l)}^+;$$

iii) if $\tilde l=1$ and $(i,j)\in{\Cal Z}\subseteq I_0\times I_0$ there exists $\tilde c_{p,\sigma,u}\in\C(q)$ ($p\in\Z^l$, $\sigma\in\sy_l$, $u\in\{0,...,l\}$) such that for all $(r,s)\in\Z^l\times\Z^{\tilde l}$
$$N_{(i,j)}(r,s)=$$
$$=\!\!\!\!\sum_{{p\in\Z^l,\sigma\in\sy_l\atop u=0,...,l}}\!\!\!\!\tilde c_{p,\sigma,u}X_{i,\tilde d_i(r_{\sigma(1)}+p_1)}^{+}...X_{i,\tilde d_i(r_{\sigma(u)}+p_u)}^{+}X_{j,\tilde d_js}^{+}X_{i,\tilde d_i(r_{\sigma(u+1)}+p_{u+1})}^{+}...X_{i,\tilde d_i(r_{\sigma(l)}+p_l)}^{+}.$$
Then the conditions of remark \ccc\ are satisfied.

\vskip .3 truecm
\rem{\ggg}~
The relations $(XD^{\pm})$-$(S3^{\pm})$, as well as $(Tk^{\pm})$ 
and $(S^{\pm})$, are of the form described in remark \fff\ and satisfy the hypotheses of remark \bbb, so that they all satisfy the conditions of remarks \bbb\ and \ccc\ and in particular the properties stated in colollary \ddd,i).
\vskip .3 truecm
\rem{\cmh}~
If the relations $(R)$ are of the form described in remarks \eee\ and \fff\ we have that $\forall h\in I_0, p\in \Z_+$:
$$[H_{h,p},S_i(r)]=b_{hip}\sum_{u=1}^{l}S_i(r+{p\over\tilde d_{i}}e_u)\ \ \ {\roman {if}}\ \tilde l=0,$$
$$[H_{h,p},S_{(i,j)}(r,s)]=b_{hjp}S_{(i,j)}(r,s+{p\over\tilde d_{j}})+b_{hip}\sum_{u=1}^{l}S_{(i,j)}(r+{p\over\tilde d_{i}}e_u,s)\ \ {\roman {if}}\ \tilde l=1,$$
where $S_{\zeta}(r,s)=0$ if $(r,s)\not\in\Z^{l+\tilde l}$.

\vskip .3 truecm
Our next goal is studying the ideals ${\Cal I}^{\pm}(R)$ and ${\Cal I}(R)$ (see notations \drmr.\aaa\ and \relid), providing a set of generators smaller and simpler than all of $\{S_{\zeta}^{\pm}(r,s)|\zeta\in{\Cal Z},r\in\Z^l,s\in\Z^{\tilde l}\}$. 
More precisely we shall show that under suitable hypotheses (fulfilled by the relations defining ${\Cal U}_q^{Dr}$ over  $\tilde{\Cal U}_q^{Dr}$) we have ${\Cal I}^{\pm}(R)={\Cal I}_{cte}^{\pm}(R)$.
\vskip .3truecm

\rem{\ctxd}~
The relations $(XD^+)$ satisfy the conditions of remark \rgen,iv), hence in particular ${\Cal I}(XD^+)={\Cal I}_{cte}
(XD^+)$ and ${\Cal I}(XD)={\Cal I}_{cte}(XD)$ (see remarks \bbb\ and \ggg).

\vskip .3truecm

We shall generalize in two steps this result for $(XD^{\pm})$ to relations $(R)$ satisfying the properties described in remarks \bbb\ and \fff: the cases $l=2$, $\tilde l=0$ (in particular $(X1^{\pm})$ and $(X2^{\pm})$) shall follow from lemma \dvssdz, while the general case will be an application of lemma \dvss.

Remark that if we considered $\tilde{\Cal U}_q^{Dr}/{\Cal I}(HX)$ instead of $\tilde{\Cal U}_q^{Dr}$ we would not need to deal with the two cases separately, but the result would follow in both cases from lemma \dvss. 

For the next remark recall notation \dautemb.\zl.
\vskip .3truecm

\rem{\sztmdz}
Consider an algebra $\U$ over a field of characteristic 0, an automorphism $t$ of $\U$, and elements $z, N(r)\in\U$ ($r\in\Z^2$), such that: 

i) $t(N(r))=N(r+\1)$ $\forall r\in\Z^2$;

ii) $[z,N(r)]=N(r+e_1)+N(r+e_2)=N(r_1+1,r_2)+N(r_1,r_2+1)$  $\forall r=(r_1,r_2)\in\Z^2$. 

If we put $S(r)=\sum_{\sigma\in\sy_2}N(\sigma(r))$ then of course: 

a) $S(r)=S(\bar r)$ $\forall r\in\Z^2$; 

b) $S(r)$ satisfies i) and ii);

c) $S(0)=2N(0)$.

\vskip .3truecm

\lem{\dvssdz}
Let $\U$, $t$, $z$, $N(r)$, $S(r)$ be as in remark \sztmdz.

If $N(0)=0$ then $S(r)=0$ $\forall r\in\Z^l$.
\dim
First of all remark \sztmdz\ implies that it is enough to prove that
$S(0,r)=0$ $\forall r\in\N$: indeed a) of remark \sztmdz\ implies that one can suppose $r_1\leq r_2$; moreover applying $t^{-r_1}$ one reduces to the case $r_1=0$. 

Let us proceed by induction on $r$: 
if $r=0$ the claim is true by hypothesis; 
let $r>0$; then 
by the inductive hypothesis
$S(0,r-1)=0$ and $0=[z,S(0,r-1)]=
S(1,r-1)+S(0,r)$; if $r=1$ $S(1,r-1)+S(0,r)=S(1,0)+S(0,1)=2S(0,1)$, so that $S(0,1)=0$; if $r>1$ then 
$S(1,r-1)=t(S(0,r-2))=0$ by the inductive hypothesis, so that also $S(0,r)=0$.

\vskip .3truecm

\rem{\sztm}~

Consider an algebra $\U$ over a field of characteristic 0, an automorphism $t$ of $\U$, elements $z_m, N_y(r)\in\U$ ($m\in\Z_+$, $y\in\U$, $r\in\Z^l$ with $l\in\Z_+$ fixed), such that: 

i) $t(N_y(r))=N_{y}(r+\1)$ $\forall y\in\U, r\in\Z^l$;

ii) $[z_m,N_y(r)]=N_{[z_m,y]}(r)+\sum_{u=1}^lN_y(r+me_u)$. 

If we put $S_y(r)=\sum_{\sigma\in\sy_l}N_y(\sigma(r))$ then of course: 

a) $S_y(r)=S_y(\sigma(r))$ $\forall\sigma\in\sy_l$;

b) $S_y(r)$ satisfies ii) and iii);

c) $S_y(0)=l!N_y(0)$. 
\vskip .3truecm

\lem{\dvss}
Let $\U$, $t$, $z_m$, $N_y(r)$, $S_y(r)$ be as in remark \sztm\ and let 
$Y\subseteq\U$ be a subset such that $[z_m,Y]\subseteq Y$ $\forall m\in\Z_+$.

If $N_y(0)=0$ $\forall y\in Y$ then $S_y(r)=0$ $\forall y\in Y$ $\forall r\in\Z^l$.
\dim
First of all remark \sztm\ implies that it is enough to prove that
$S_y(r)=0$ $\forall y\in Y$ $\forall r=(r_1,...,r_l)\in\Z^l$ such that $0=r_1\leq...\leq r_l$: indeed a) of remark \sztm\ implies that one can suppose $r_1\leq...\leq r_l$; moreover applying $t^{-r_1}$ one reduces to the case $r_1=0$. 

Let $v=max\{u=1,...,l|r_u=0\}$ and proceed by induction on $v$: 
if $v=l$ then $r=~\!\!0$ and the claim is true by hypothesis; 
let $v<l$ and choose $m=r_{v+1}$; then 
$$max\{u=1,...,l|(r-me_{v+1})_u=0\}=v+1$$ and $$max\{u=1,...,l|(r-me_{v+1}+me_u)_u=0\}=v+1\ \ \forall u>v+1,$$ hence by the inductive hypothesis
$S_y(r-me_{v+1})=0$ and $S_y(r-me_{v+1}+me_u)=0$ $\forall y\in Y$ $\forall u>v+1$; it follows that
$$0=[z_m,S_y(r-me_{v+1})]=
S_{[z_m,y]}(r-me_{v+1})+\sum_{u=1}^lS_y(r-me_{v+1}+me_u)=$$
$$=\sum_{u=1}^{v+1}S_y(r-me_{v+1}+me_u)=(v+1)S_y(r).$$

\vskip .3truecm
\prop{\ancid}
Consider the notations fixed in notation \drmr.\aaa\ and suppose that $(R^+)$ satisfies the hypotheses of remark \fff. Then: 

i) if $l=2$ and $\tilde l=0$ we have ${\Cal I}^+(R)={\Cal I}_{cte}^+(R)$ in $\tilde{\Cal U}^{Dr}_q(X_{\tilde n}^{(k)})$; 

ii) in any case ${\Cal I}^+(R)={\Cal I}_{cte}^+(R)$ in $\tilde{\Cal U}^{Dr}_q(X_{\tilde n}^{(k)})/\Cal I(HX^+)$;

iii) if moreover the hypotheses of remark \bbb\ are satisfied then we have also ${\Cal I}^-(R)={\Cal I}_{cte}^-(R)$ (in $\tilde{\Cal U}^{Dr}_q(X_{\tilde n}^{(k)})$ and in $\tilde{\Cal U}^{Dr}_q(X_{\tilde n}^{(k)})/{\Cal I}(HX^-)$ respectively).
\dim
Let $\zeta\in{\Cal Z}$; thanks to remarks \ccc,iv) and \darl.\tind, for all $i\in I_0$ $\tilde t_i$ induces an automorphism $t_i^{\prime}$ of $\tilde{\Cal U}^{Dr}_q(X_{\tilde n}^{(k)})/{\Cal I}_{cte}^+(R_{\zeta})$ and of $\tilde{\Cal U}^{Dr}_q(X_{\tilde n}^{(k)})/{\Cal I}_{cte}^+(HX,R_{\zeta})$.

i) Fix $i\in{\Cal Z}\subseteq I_0$, and notice that the data $$\U=\tilde{\Cal U}^{Dr}_q/\I_{cte}^+(R_i),\ \ t=t_i^{\prime-1},\ \ z={1\over b_{ii \tilde d_i}}H_{(i,\tilde d_i )},\ \ N(r)=N_i(r)\ \ (r\in\Z^2)$$ 
satisfy the conditions of remark \sztmdz\ and lemma \dvssdz. 

Since we have that 
${\Cal I}_0^+(R_i)=0$ in $\tilde{\Cal U}_q^{Dr}/{\Cal I}_{cte}^+(R_i)$, lemma \dvssdz\ implies that
${\Cal I}^+(R_{i})=0$ in $\tilde{\Cal U}_q^{Dr}/{\Cal I}_{cte}^+(R_{i})$, or equivalently that in $\tilde{\Cal U}_q^{Dr}$
${\Cal I}^+(R_{i})\subseteq{\Cal I}_{cte}^+(R_{i})$, and this for all $i\in I_0\subseteq{\Cal Z}$, so that ${\Cal I}^+(R)={\Cal I}_{cte}^+(R)$
thanks to remark \rgen,i) and ii). 

ii) Fix $\zeta=\begin{cases}  i\in{\Cal Z}\subseteq I_0&{\roman {if}}\ \tilde l=0\cr
(i,j)\in{\Cal Z}\subseteq I_0\times I_0&{\roman {if}}\ \tilde l=1\end{cases} $, and notice that the data $$\U=\tilde{\Cal U}^{Dr}_q/\I_{cte}^+(HX,R_{\zeta}),\ \ t=t_i^{\prime-1},\ \ z_m={1\over b_{ii \tilde d_i\!m}}H_{(i,\tilde d_i m)}\ \forall m\in\Z_+,$$ 
$$Y=\begin{cases} \{0,1\}&{\roman{if\ }}
\tilde l=0\cr
\{aX_{j,\tilde d_js}^+|s\in\Z,\ a\in\C(q)\}&{\roman{if\ }}\tilde l=1,\end{cases} $$
and for $r\in\Z^l$
$$N_y(r)=\begin{cases}  yN_i^+(r)&{\roman{if\ }}\tilde l=0\cr
aN_{(i,j)}^+(r,s)&{\roman{if\ }}\tilde l=1, y=aX_{j,\tilde d_js}^+\end{cases} $$
satisfy the conditions of remark \sztm\ and lemma \dvss. 

Since we have that 
${\Cal I}_0^+(R_i)=0$ in $\tilde{\Cal U}_q^{Dr}/{\Cal I}_{cte}^+(HX,R_i)$ if $\tilde l=0$ and  
${\Cal I}_{0,s}^+(R_{(i,j)})=0$ in $\tilde{\Cal U}_q^{Dr}/{\Cal I}_{cte}^+(HX,R_{(i,j)})$ for all $s\in\Z$ if 
$\tilde l=1$, lemma \dvss\ implies that
${\Cal I}^+(R_{\zeta})=0$ in $\tilde{\Cal U}_q^{Dr}/{\Cal I}_{cte}^+(HX,R_{\zeta})$, or equivalently that in $\tilde{\Cal U}_q^{Dr}/{\Cal I}_{cte}^+(HX)$
${\Cal I}^+(R_{\zeta})\subseteq{\Cal I}_{cte}^+(R_{\zeta})$, so that ${\Cal I}^+(R)={\Cal I}_{cte}^+(R)$
thanks to remark \rgen,i) and ii). 

It follows that ${\Cal I}^+(R)\subseteq({\Cal I}^+(HX),{\Cal I}_{cte}^+(R))$.

iii) follows from i) and ii) thanks to remarks \bbb\ and \darl.\rhxom.

\vskip .3 truecm
\rem{\iii}~
In proposition \ancid,i) the hypothesis $\tilde l=0$ is not necessary: the claim would hold also in case $\tilde l=1$. But this case is omitted here because it is not really needed in this paper and its proof, very  similar, would just require a little more complicated, and repetitive, exposition (see the proof of proposition \ancid,ii)).

\vskip 1 truecm
\cor{\hhh}
i) ${\Cal I}^{\pm}(X1)={\Cal I}_{cte}^{\pm}(X1)$ and ${\Cal I}^{\pm}(X2)={\Cal I}_{cte}^{\pm}(X2)$. 

ii) If $(R^{\pm})$ is one of $(X3^{\varepsilon,\pm})$-$(S3^{\pm})$, $(Tk^{\pm})$ and $(S^{\pm})$ then 
$${\Cal I}^{\pm}(R)\subseteq({\Cal I}^{\pm}(HX),{\Cal I}_{cte}^{\pm}(R)).$$

iii) If $(R^{\pm})$ is one of $(X3^{\varepsilon,\pm})$-$(S3^{\pm})$ 
and $(S^{\pm})$ then 
$${\Cal I}^{\pm}(R,XD,X1,X2)={\Cal I}_{cte}^{\pm}(R,XD,X1,X2).$$
\dim
The claims follow from proposition \ancid, remarks \ggg\ and \ctxd\ and proposition \darl.\sop.

\vskip .3 truecm
\rem{\ctepg}~
Remark \iii\ would imply that furthermore ${\Cal I}^{\pm}(R)={\Cal I}_{cte}^{\pm}(R)$ even in the case when
$(R^{\pm})$  is one of $(Sk^{\pm})$, $(Tk^{\pm})$ and, if $a_{ij}=-1$, also $(S(UL)_{(i,j)}^{\pm})$.

\vskip .3 truecm
\cor{\nrelm}
Proposition \ancid\ implies that $${\Cal U}_q^{Dr}\!\!=\tilde{\Cal U}_q^{Dr}/{\Cal I}_{cte}(XD,X1,X2,X3^{\varepsilon},SUL,S2,S3).$$

\vskip .3 truecm
The final remark of this section is a straightforward consequence of remark \dautemb.\xibdef.
\vskip .3 truecm
\rem{\stabid}~
i) If $(R^{\pm})$ is one of $(XD^{\pm})$-$(X2^{\pm})$, $(SUL^{\pm})$-$(S3^{\pm})$, $(Tk^{\pm})$, $(S^{\pm})$ then ${\Cal I}^+(R)$ and ${\Cal I}_{cte}^+(R)$ 
are $\tilde\Theta$-stable.

ii) ${\Cal I}^+(X3^1,X3^{-1})$ and ${\Cal I}_{cte}^+(X3^1,X3^{-1})$ are the $\tilde\Theta$-stable ideal generated respectively by ${\Cal I}^+(X3^1)$ and 
${\Cal I}_{cte}^+(X3^1)$;

iii) ${\Cal I}_{cte}^+(X3^1,X3^{-1})$ is the $\tilde\Theta$-stable and $t_i$-stable ideal (for all $i\in I_0$) generated by 
${\Cal I}_0^+(X3^1)$;

iv) ${\Cal I}(X3^1,X3^{-1},X2)={\Cal I}_{cte}(X3^1,X3^{-1},X2)$ is the $\tilde\Omega$-stable, $\tilde\Theta$-stable and $\tilde t_i$-stable ideal (for all $i\in I_0$) generated by 
${\Cal I}_0^+(X3^1,X2)$.

\vskip .5 truecm
{\bf{\cnat.\ MORE about REDUNDANT RELATIONS. }}
\vskip .5 truecm

In this section we prove some dependences among the relations $(XD^{\pm})$-$(S3^{\pm})$ making  systematic recourse to the properties of $q$-commutators (remark \drrl.\rulc) and to corollary \dlem.\ddd,ii). 

\vskip .3 truecm

\prop{\rku}
With the notations of remark \dlem.\relid, 
${\Cal I}_{cte}^{\pm}(X2)$, ${\Cal I}_{cte}^{\pm}(X3^{-1})\subseteq{\Cal I}_{cte}^{\pm}(X3^{1})$. 
\dim
${\Cal I}_{0}^+(X2)\subseteq{\Cal I}_{cte}^+(X3^{1})$:
$$[[[X_{1,1}^+,X_{1,0}^+]_{q^2},X_{1,0}^+]_{q^4},X_{1,-1}^-]\in{\Cal I}_{cte}^+(X3^1);$$
but $$[[[X_{1,1}^+,X_{1,0}^+]_{q^2},X_{1,0}^+]_{q^4},X_{1,-1}^-]=
\Big[\Big[{Ck_1-C^{-1}k_1^{-1}\over q-q^{-1}},X_{1,0}^+\Big]_{q^2},X_{1,0}^+\Big]_{q^4}+$$
$$+[[X_{1,1}^+,k_1^{-1}H_{1,-1}]_{q^2},X_{1,0}^+]_{q^4}+
[[X_{1,1}^+,X_{1,0}^+]_{q^2},k_1^{-1}H_{1,-1}]_{q^4}=$$
$$=k_1^{-1}(q^2[2]_qC^{-1}[X_{1,0}^+,X_{1,0}^+]_{q^6}-
q^2[[H_{1,-1},X_{1,1}^+],X_{1,0}^+]_{q^6}-
q^4[H_{1,-1},[X_{1,1}^+,X_{1,0}^+]_{q^2}])\!\!=$$
$$=-q^4[3]_q!C^{-1}k_1^{-1}([X_{1,1}^+,X_{1,-1}^+]_{q^2}-(q^4-q^{-2})(X_{1,0}^+)^2)=-q^4[3]_q!C^{-1}k_1^{-1}M_{(2,2)}^+(0);$$ 
then
${\Cal I}_{0}^+(X2)\subseteq{\Cal I}_{cte}^+(X3^{1})$.

${\Cal I}_{0}^+(X3^{-1})\subseteq{\Cal I}_{cte}^+(X3^{1})$: 
notice that
$${\Cal I}_{cte}^+(X3^1)\ni 
{1\over b_{ii\tilde d_i}}\tilde t_1([H_{1,1},M_{(3)}^{1,+}(0)])=$$
$$=[[X_{1,1}^+,X_{1,-1}^+]_{q^2},X_{1,-1}^+]_{q^4}+[[X_{1,0}^+,X_{1,0}^+]_{q^2},X_{1,-1}^+]_{q^4}+[[X_{1,0}^+,X_{1,-1}^+]_{q^2},X_{1,0}^+]_{q^4}=$$
$$=[[X_{1,1}^+,X_{1,-1}^+]_{q^2}-(q^4-q^{-2})(X_{1,0}^+)^2,X_{1,-1}^+]_{q^4}+(q^4-q^{-2}+1-q^2-q^4)(X_{1,0}^+)^2X_{1,-1}^++$$
$$+(1+q^6)X_{1,0}^+X_{1,-1}^+X_{1,0}^++(-q^8+q^{2}-q^4+q^6-q^2)X_{1,-1}^+(X_{1,0}^+)^2,$$
so that 
$$-(q^2-1+q^{-2})((X_{1,0}^+)^2X_{1,-1}^+-(q^4+q^2)X_{1,0}^+X_{1,-1}^+X_{1,0}^++q^6X_{1,-1}^+(X_{1,0}^+)^2)=$$
$$=-q^6(q^2-1+q^{-2})[[X_{1,-1}^+,X_{1,0}^+]_{q^{-2}},X_{1,0}^+]_{q^{-4}}=-q^6(q^2-1+q^{-2})M_{(3)}^{-1,+}(0)$$
is an element of ${\Cal I}_{cte}^+(X3^1)$. The claims follow again from corollary \dlem.\ddd.

\vskip .3 truecm
\cor{\tsmd}
i) ${\Cal I}^{\pm}(X3^1)={\Cal I}_{cte}^{\pm}(X3^{1})$ (see corollary \dlem.\hhh,iii)).

ii) ${\Cal I}_{cte}^+(X3^1)={\Cal I}_{cte}^+(X3^{-1})$ is $\tilde\Theta$-stable. 

iii) Moreover $({\Cal I}_{cte}^+(X3^1),{\Cal I}_{cte}^-(X3^{-1}))={\Cal I}_{cte}(X3^{1})={\Cal I}_{cte}(X3^{-1})$  is $\tilde\Omega$-stable.

\vskip .3 truecm
\prop{\prku}
i) ${\Cal U}_q^{Dr}(A_1^{(1)})=\tilde{\Cal U}_q^{Dr}(A_1^{(1)})/{\Cal I}_{cte}(X1)$;

ii) ${\Cal U}_q^{Dr}(A_2^{(2)})=\tilde{\Cal U}_q^{Dr}(A_2^{(2)})/{\Cal I}_{cte}(X3^{1})=\tilde{\Cal U}_q^{Dr}(A_2^{(2)})/({\Cal I}_{cte}^+(X3^{1}),{\Cal I}_{cte}^-(X3^{-1}))$.

\vskip .3 truecm

\prop{\rdxd}
${\Cal I}_{cte}^{\pm}(XD)\subseteq{\Cal I}_{cte}^{\pm}(SUL)$. 
\dim
Let $i,j\in I_0$ be such that $a_{ij}<0$; since $-1\in\{a_{ij},a_{ji}\}$, in the study of $[X_{i,\tilde d_{ij}}^+,X_{j,0}^+]_{q_i^{a_{ij}}}+
[X_{j,\tilde d_{ij}}^+,X_{i,0}^+]_{q_i^{a_{ij}}}$ we can suppose that 
$a_{ij}=-1$, and in particu\-lar $\tilde d_j\leq\tilde d_i=\tilde d_{ij}$ and, if $X_{\tilde n}^{(k)}= A_{2n}^{(2)}$, $i\neq 1$. 
Then $[[X_{j,0}^+,X_{i,0}^+]_q,X_{i,0}^+]_{q^{-1}}$ is an element of ${\Cal I}_{cte}^+(SUL)$, and so is $[[[X_{j,0}^+,X_{i,0}^+]_q,X_{i,0}^+]_{q^{-1}},X_{i,\tilde d_{i}}^-]$. 
But 
$$[[[X_{j,0}^+,X_{i,0}^+]_q,X_{i,0}^+]_{q^{-1}},X_{i,\tilde d_{i}}^-]=$$
$$=[[X_{j,0}^+,X_{i,0}^+]_q,C^{-\tilde d_i}k_iH_{i,\tilde d_i}]_{q^{-1}}+[[X_{j,0}^+,C^{-\tilde d_i}k_iH_{i,\tilde d_i}]_q,X_{i,0}^+]_{q^{-1}}=$$
$$=C^{-\tilde d_i}k_i(q^{-1}[[X_{j,0}^+,X_{i,0}^+]_q,H_{i,\tilde d_i}]+q[[X_{j,0}^+,H_{i,\tilde d_i}],X_{i,0}^+]_{q^{-3}})=$$
$$=C^{-\tilde d_i}k_i(-b_{ii\tilde d_i}q^{-1}[X_{j,0}^+,X_{i,\tilde d_i}^+]_q-b_{ij\tilde d_i}q^{-1}[X_{j,\tilde d_i}^+,X_{i,0}^+]_q-b_{ij\tilde d_i}q[X_{j,\tilde d_i}^+,X_{i,0}^+]_{q^{-3}})=$$
$$=[2]_{q_i}C^{-\tilde d_i}k_i([X_{i,\tilde d_i}^+,X_{j,0}^+]_{q^{-1}}+
[X_{j,\tilde d_i}^+,X_{i,0}^+]_{q^{-1}}).$$
Then ${\Cal I}_0^+(XD)\subseteq{\Cal I}_{cte}^+(SUL)$, and the claim follows using corollary \dlem.\ddd.

\vskip .3 truecm
\lem{\cyi}
For $i\in I_0$, $a\in\N$ define $Y_{i,a}\in\tilde{\Cal U}^{Dr}_q(X_{\tilde n}^{(k)})$ as follows: 
$$Y_{i,0}=X_{i,\tilde d_i}^+;\ \ \ Y_{i,a+1}=[Y_{i,a},X_{i,0}^+]_{q_i^{2(a+1)}}.$$
(Notice that $Y_{i,1}=M_i^+(0)$).

Then:

i) $[Y_{i,a},X_{i,0}^-]=(b_{ii\tilde d_i}-[a]_{q_i}[a+1]_{q_i})k_iY_{i,a-1}$ $\forall a>0$;

ii) $[X_{i,j;1-a_{ij};1}^+(0;0),X_{j,\tilde d_{ij}}^-]=C^{-\tilde d_{ij}}k_jb_{ji\tilde d_{ij}}Y_{i,-a_{ij}}$.
\dim
i) $[Y_{i,a},X_{i,0}^-]=$
$$=[[...[...[X_{i,\tilde d_i}^+,X_{i,0}^+]_{q_i^2},...X_{i,0}^+]_{q_i^{2u}},...X_{i,0}^+]_{q_i^{2a}},X_{i,0}^-]=$$
$$=[...[...[k_iH_{i,\tilde d_i},X_{i,0}^+]_{q_i^2},...X_{i,0}^+]_{q_i^{2u}},...X_{i,0}^+]_{q_i^{2a}}+$$
$$+\sum_{u=1}^a[...[\Big[[...[X_{i,\tilde d_i}^+,X_{i,0}^+]_{q_i^2},...X_{i,0}^+]_{q_i^{2(u-1)}},{k_i-k_i^{-1}\over q_i-q_i^{-1}}\Big]_{q_i^{2u}},X_{i,0}^+]_{q_i^{2(u+1)}},...X_{i,0}^+]_{q_i^{2a}}=$$
$$=k_i[...[...[H_{i,\tilde d_i},X_{i,0}^+],...X_{i,0}^+]_{q_i^{2(u-1)}},...X_{i,0}^+]_{q_i^{2(a-1)}}+$$
$$+\sum_{u=1}^a{q_i^{-2u}-q_i^{2u}\over q_i-q_i^{-1}}k_i[...[[...[X_{i,\tilde d_i}^+,X_{i,0}^+]_{q_i^2},...X_{i,0}^+]_{q_i^{2(u-1)}},X_{i,0}^+]_{q_i^{2u}},...X_{i,0}^+]_{q_i^{2(a-1)}}=$$
$$=b_{ii\tilde d_i}k_i[...[...[X_{i,\tilde d_i}^+,X_{i,0}^+]_{q_i^{2}},...X_{i,0}^+]_{q_i^{2(u-1)}},...X_{i,0}^+]_{q_i^{2(a-1)}}+$$
$$-\sum_{u=1}^a[2u]_{q_i}k_i[...[[...[X_{i,\tilde d_i}^+,X_{i,0}^+]_{q_i^2},...X_{i,0}^+]_{q_i^{2(u-1)}},X_{i,0}^+]_{q_i^{2u}},...X_{i,0}^+]_{q_i^{2(a-1)}}=$$
$$=(b_{ii\tilde d_i}-[a]_{q_i}[a+1]_{q_i})k_iY_{i,a-1};$$.

ii) $[X_{i,j;1-a_{ij};1}^+(0;0),X_{j,\tilde d_{ij}}^-]=$
$$=[[...[...[X_{j,0}^+,X_{i,0}^+]_{q_i^{a_{ij}}},...X_{i,0}^+]_{q_i^{a_{ij}+2(u-1)}},...X_{i,0}^+]_{q_i^{-a_{ij}}},X_{j,\tilde d_{ij}}^-]=$$
$$={1\over q_j-q_j^{-1}}[...[...[C^{-\tilde d_{ij}}k_j\tilde H_{j,\tilde d_{ij}}^+,X_{i,0}^+]_{q_i^{a_{ij}}},...X_{i,0}^+]_{q_i^{a_{ij}+2(u-1)}},...X_{i,0}^+]_{q_i^{-a_{ij}}}=$$
$$={C^{-\tilde d_{ij}}k_j\over q_j-q_j^{-1}}[...[...[\tilde H_{j,\tilde d_{ij}}^+,X_{i,0}^+],...X_{i,0}^+]_{q_i^{2(u-1)}},...X_{i,0}^+]_{q_i^{-2a_{ij}}}.$$
Recalling remark \darl.\solocom\ we get
$$[X_{i,j;1-a_{ij};1}^+(0;0),X_{j,\tilde d_{ij}}^-]=$$
$$=C^{-\tilde d_{ij}}k_jb_{ji\tilde d_{ij}}[...[...[X_{i,\tilde d_{ij}}^+,X_{i,0}^+]_{q_i^2},...X_{i,0}^+]_{q_i^{2(u-1)}},...X_{i,0}^+]_{q_i^{-2a_{ij}}}=$$
$$=C^{-\tilde d_{ij}}k_jb_{ji\tilde d_{ij}}Y_{i,-a_{ij}}.$$

\vskip .3truecm     

\cor{\misp}
Let $i,j\in I_0$ be such that $a_{ij}<0$ with the condition that $a_{ij}=-1$ if $k>1$; then $M_i^+(0)\in{\Cal I}_{cte}^+(SUL)$.

In particular:

i) in the cases of rank higher than 1 (that is $X_{\tilde n}^{(k)}\neq A_{1}^{(1)},A_{2}^{(2)}$) and different from $D_{n+1}^{(2)}$ and $D_{4}^{(3)}$ we have  ${\Cal I}_{cte}^{\pm}(X1)\subseteq{\Cal I}_{cte}^{\pm}(SUL)$;

ii) in the cases $D_{n+1}^{(2)}$ and $D_{4}^{(3)}$ we have ${\Cal I}_{cte}^{\pm}(X1)\subseteq({\Cal I}_{cte}^{\pm}(SUL),{\Cal I}_{cte}^{\pm}(X1_1))$.

\dim
This is an immediate consequence of lemma \cyi\ (and of corollary \dlem.\ddd) once one notices that the hypotheses imply that 
$X_{i,j;1-a_{ij};1}^+(0;0)\in{\Cal I}_{cte}^+(SUL)$, $b_{ji\tilde d_{ij}}\neq 0$
and $b_{ii\tilde d_i}=[2]_{q_i}$.

\vskip .3 truecm

\rem{\xitd}~
If $k=2$, $X_{\tilde n}^{(k)}\neq A_{2n}^{(2)}$ and $i,j\in I_{0}$ are such that $a_{ij}=-2$ then
$${\Cal I}_{cte}^{\pm}(X1_{i})\subseteq{\Cal I}_{cte}^{\pm}(T2).$$ In particular ${\Cal I}_{cte}^{\pm}(T2)={\Cal I}_{cte}^{\pm}(X1_{i},S2)$.
\dim
$[[X_{j,0}^{+},X_{i,1}^{+}]_{q^2},X_{i,0}^{+}]$ lies in ${\Cal I}_{cte}^{+}(T2)$
and so does
$$[[[X_{j,0}^{+},X_{i,1}^{+}]_{q^2},X_{i,0}^{+}],X_{j,0}^{-}]=$$
$$=\Big[\Big[{k_j-k_j^{-1}\over q^2-q^{-2}},X_{i,1}^{+}\Big]_{q^2},X_{i,0}^{+}\Big]=$$
$$=-q^2k_j[X_{i,1}^{+},X_{i,0}^{+}]_{q^2}.$$
\vskip .3 truecm 
\teo{\tldg}
i) ${\Cal U}_q^{Dr}(X_{\tilde n}^{(1)})=\begin{cases} \tilde{\Cal U}_q^{Dr}(A_1^{(1)})/{\Cal I}_{cte}(X1)&{\roman{if}}\ X_{\tilde n}=A_1
\cr
\tilde{\Cal U}_q^{Dr}(X_{\tilde n}^{(1)})/{\Cal I}_{cte}(SUL)&{\roman{otherwise;}}\end{cases} $

ii) ${\Cal U}_q^{Dr}(X_{\tilde n}^{(2)})=\begin{cases} \tilde{\Cal U}_q^{Dr}(A_{2n}^{(2)})/{\Cal I}_{cte}(X3^1,SUL,S2)&{\roman{if}}\ X_{\tilde n}=A_{2n}\cr
\tilde{\Cal U}_q^{Dr}(X_{\tilde n}^{(2)})/{\Cal I}_{cte}(SUL,T2)&{\roman{otherwise;}}\end{cases} $

iii) ${\Cal U}_q^{Dr}(D_{4}^{(3)})\!=\!\tilde{\Cal U}_q^{Dr}(D_{4}^{(3)})/{\Cal I}_{cte}(X1_1,\!SUL,\!S3)\!=\!\tilde{\Cal U}_q^{Dr}(D_{4}^{(3)})/{\Cal I}_{cte}(X1_1,\!SUL,\!T3)$.

\vskip .3 truecm

\cor{\hmbdefnt}
Let: $\alg$ be a $\C(q)$-algebra, $t_i^{(\alg)}$ ($i\in I_0$) be $\C(q)$-automorphisms of $\alg$, $\Omega^{(\alg)}$ be a $\C$-anti-linear anti-automorphism of $\alg$, $\tilde f:\tilde{\Cal U}^{Dr}_q(X_{\tilde n}^{(k)})\rightarrow\alg$ be a homomorphism of $\C(q)$-algebras such that $\tilde f\compo \tilde t_i=t_i^{(\alg)}\compo\tilde f$ $\forall i\in I_0$ and $\tilde f\compo \tilde\Omega=\Omega^{(\alg)}\compo\tilde f$. 

If:

i) $\tilde f({\Cal I}_0^+(X1))=0$ in case $X_{\tilde n}^{(k)}=A_1^{(1)}$;

ii) $\tilde f({\Cal I}_0^+(SUL))=0$ in case $k=1$, $X_{\tilde n}^{(k)}\neq A_1^{(1)}$;

iii) $\tilde f({\Cal I}_0^+(X3^1,SUL,S2))=0$ in case $X_{\tilde n}^{(k)}=A_{2n}^{(2)}$; 

iv) $\tilde f({\Cal I}_0^+(SUL,T2))=0$ in case $k=2$, $X_{\tilde n}^{(k)}\neq A_{2n}^{(2)}$;

v) $\tilde f({\Cal I}_0^+(X1_1,SUL,T3))=0$ in case $D_{4}^{(3)}$;

then:

$\tilde f$ induces $f:\Udr_q(X_{\tilde n}^{(k)})\rightarrow\alg$ and $f\compo t_i=t_i^{(\alg)}\compo f$ $\forall i\in I_0$, $f\compo \Omega=\Omega^{(\alg)}\compo f$. 
\dim 
Since the hypotheses imply that $ker(\tilde f)$ is a $t_i$-stable ($\forall i\in I_0$) and $\tilde\Omega$-stable ideal of $\Udr_q(X_{\tilde n}^{(k)})$, the claim is an immediate consequence of theorem \hmbdefnt\ and of corollary \dlem.\ddd,i).

\vskip .5 truecm

{\bf{\serrel.\  The SERRE RELATIONS}}.
\vskip .5 truecm

This section is dovoted to the study of the Serre relations (see definition \drmr.\relser). In particular we prove that the Serre relations hold in ${\Cal U}_q^{Dr}$, and that in the case of rank higher than 1 the Serre relations alone are indeed equivalent  to $(XD^{\pm})$-$(S3^{\pm})$ (in $\tilde{\Cal U}_q^{Dr}$), that is ${\Cal U}_q^{Dr}=\tilde{\Cal U}_q^{Dr}/{\Cal I}_{cte}(S^{\pm})$.
We use the notations fixed in notation \drmr.\aaa\ and \dlem.\relid.

\vskip .3 truecm

\rem{\mij}~
i) If $k=1$ 
$(S^{\pm})=(SUL^{\pm})$;

ii) if $k>1$ and $i,j\in I_0$ are such that $a_{ij}<-1$ then $(S^{\pm})=(SUL^{\pm})\cup(S_{(i,j)}^{\pm})$.

\vskip .3 truecm
Before passing to prove that the Serre relations hold in $\Udr_q$, we state the following remark on $q$-commutators, which simplifies many computations in the next propositions. 

\vskip .3 truecm
\rem{\cid}~
Let $a\in\tilde{\Cal U}_q^{Dr}$, $i\in I_0$ such that $(X_{\tilde n}^{(k)},i)\neq(A_{2n}^{(2)},1)$, $u,v\in\C(q)$. Then in $\tilde{\Cal U}_q^{Dr}/{\Cal I}_{cte}^+(X1_i)$ we have, for all $r\in \Z$:

i) $[[a,X_{i,r}^+]_u,X_{i,r+\tilde d_i}^+]_v=q_i^{-2}[[a,X_{i,r+\tilde d_i}^+]_{q_i^2v},X_{i,r}^+]_{q_i^2u}$;

ii) $[[a,\!X_{i,r+2\tilde d_i}^+]_u,\!X_{i,r}^+]_v\!=\!q_i^{2}[[a,\!X_{i,r}^+]_{q_i^{-2}\!v},
\!X_{i,r+2\tilde d_i}^+]_{q_i^{-2}\!u}\!+(q^2-1)[a,\!(X_{i,r+\tilde d_i}^+)^2]_{q_i^{-2}\!uv}$.
\dim
It is a simple computation using remark \drrl.\rulc,iii). 

\vskip .3 truecm

\prop{\tds}
If $k=2$, $X_{\tilde n}^{(2)}\neq A_{2 n}^{(2)}$ and $i,j\in I_0$ are such that $a_{ij}=-2$, then 
$${\Cal I}_{cte}^{\pm}(S_{(i,j)})\subseteq{\Cal I}_{cte}^{\pm}(T2).$$
\dim
In this proof we use that ${\Cal I}_{cte}^{\pm}(X1_i)\subseteq{\Cal I}_{cte}^{\pm}(T2)$ (see proposition \cnat.\xitd) 
and make the following computations in $\tilde{\Cal U}^{Dr}_q/{\Cal I}_{cte}^+(T2)$. 

Since $[H_{i,1},[[X_{j,0}^+,X_{i,0}^+]_{q^2},X_{i,-1}^+]]$ lies in ${\Cal I}_{cte}^+(T2)$ (see definition \drmr.\skq) we see that 
$$[[X_{j,0}^+,X_{i,1}^+]_{q^{2}},X_{i,-1}^+]+[[X_{j,0}^+,X_{i,0}^+]_{q^{2}},X_{i,0}^+]=0$$ in $\tilde{\Cal U}^{Dr}_q/{\Cal I}_{cte}^+(T2)$; but, thanks to lemma \cid,ii), we have that 
$$[[X_{j,0}^+,X_{i,1}^+]_{q^{2}},X_{i,-1}^+]=q^2[[X_{j,0}^+,X_{i,-1}^+]_{q^{-2}},X_{i,1}^+]+(q^2-1)[X_{j,0}^+,(X_{i,0}^+)^2],$$
so that $$
0=q^2[[X_{j,0}^+,X_{i,-1}^+]_{q^{-2}},X_{i,1}^+]+[X_{j,0}^+,(q^2-1)(X_{i,0}^+)^2]+[[X_{j,0}^+,X_{i,0}^+]_{q^{2}},X_{i,0}^+]=$$
$$=
q^2[[X_{j,0}^+,X_{i,-1}^+]_{q^{-2}},X_{i,1}^+]+q^2[[X_{j,0}^+,X_{i,0}^+]_{q^{-2}},X_{i,0}^+]$$
and also
$$[[[X_{j,0}^+,X_{i,-1}^+]_{q^{-2}},X_{i,1}^+]+[[X_{j,0}^+,X_{i,0}^+]_{q^{-2}},X_{i,0}^+],X_{i,0}^+]_{q^{2}}=0.$$ Now, thanks to remark \cid,i) and to relations $(T2^+)$, $$[[[X_{j,0}^+,X_{i,-1}^+]_{q^{-2}},X_{i,1}^+],X_{i,0}^+]_{q^{2}}=[[[X_{j,0}^+,X_{i,0}^+]_{q^{2}},X_{i,-1}^+],X_{i,1}^+]_{q^{-2}}=0,$$
so that
$$[[[X_{j,0}^+,X_{i,0}^+]_{q^{-2}},X_{i,0}^+],X_{i,0}^+]_{q^{2}}=0,$$ which implies ${\Cal I}_{cte}^{\pm}(S_{(i,j)})\subseteq{\Cal I}_{cte}^{\pm}(T2)$, thanks to corollary \dlem.\ddd.

\vskip .3 truecm
Let us concentrate now on the case $A_{2 n}^{(2)}$.
\vskip .3 truecm
\lem{\tat}
Let $X_{\tilde n}^{(k)}= A_{2 n}^{(2)}$; then 
$$[[X_{1,2}^+,X_{1,1}^+]_{q^2},X_{1,0}^+]_{q^4}-(q^2-1)(q^4-1)(q^2+q^{-2})(X_{1,1}^+)^3\in{\Cal I}_{cte}^+(X3^1).$$
\dim
By corollary \cnat.\tsmd 
$$[[X_{1,2}^+,X_{1,1}^+]_{q^2},X_{1,0}^+]_{q^4}+[[X_{1,2}^+,X_{1,0}^+]_{q^2},X_{1,1}^+]_{q^4}+[X_{1,1}^+,X_{1,1}^+]_{q^2},X_{1,1}^+]_{q^4}$$ belongs to ${\Cal I}^+(X3^1)={\Cal I}_{cte}^+(X3^1)$. 

But $[X_{1,2}^+,X_{1,0}^+]_{q^2}-(q^4-q^{-2})(X_{1,1}^+)^2\in{\Cal I}_{cte}^+(X3^1)$ (see proposition \cnat.\rku), so that 
$$[[X_{1,2}^+,X_{1,1}^+]_{q^2},X_{1,0}^+]_{q^4}+(1-q^4)(q^4-q^{-2}+1-q^2)(X_{1,1}^+)^3$$
lies  in ${\Cal I}_{cte}^+(X3^1)$.

\vskip .3 truecm

\prop{\sds}
If $X_{\tilde n}^{(2)}= A_{2 n}^{(2)}$ and $i,j\in I_0$ are such that $a_{ij}=-2$ ($i=1$, $j=2$), then 
${\Cal I}_{cte}^{\pm}(S_{(i,j)})\subseteq{\Cal I}_{cte}^{\pm}(XD,X3^1,S2)$.
\dim
Recall that by the very definition of ${\Cal I}_{cte}^+(S2)$ we have (see remark \drmr.\sercod)
$$(q^{2}+q^{-2})[[X_{j,0}^+,X_{i,1}^+]_{q^2},X_{i,0}^+]+q^{2}[[X_{i,1}^+,X_{i,0}^+]_{q^2},X_{j,0}^+]_{q^{-4}}\in{\Cal I}_{cte}^+(S2)$$
so that also
$$(q^{2}+q^{-2})[X_{i,-1}^+,[[X_{j,0}^+,X_{i,1}^+]_{q^2},X_{i,0}^+]]_{q^{-2}}+q^{2}[X_{i,-1}^+,[[X_{i,1}^+,X_{i,0}^+]_{q^2},X_{j,0}^+]_{q^{-4}}]_{q^{-2}}$$ belongs to ${\Cal I}_{cte}^+(S2)$,
and let us compute the two summands separately in the algebra $\tilde{\Cal U}_{q}^{Dr}/{\Cal I}_{cte}^+(XD,X3^{1},S2)$:
$$[X_{i,-1}^+,[[X_{j,0}^+,X_{i,1}^+]_{q^2},X_{i,0}^+]]_{q^{-2}}=$$
$$=-q^2[X_{i,-1}^+,[[X_{i,1}^+,X_{j,0}^+]_{q^{-2}},X_{i,0}^+]]_{q^{-2}}{\buildrel  XD\over =}q^2[X_{i,-1}^+,[[X_{j,1}^+,X_{i,0}^+]_{q^{-2}},X_{i,0}^+]]_{q^{-2}}
=$$
$$=q^2[[X_{i,-1}^+,[X_{j,1}^+,X_{i,0}^+]_{q^{-2}}],X_{i,0}^+]_{q^{-2}}+q^2
[[X_{j,1}^+,X_{i,0}^+]_{q^{-2}},[X_{i,-1}^+,X_{i,0}^+]_{q^{-2}}]{\buildrel  S2\over =}$$
$${\buildrel  S2\over =}
-q^4[[X_{i,0}^+,[X_{j,1}^+,X_{i,-1}^+]_{q^{-2}}]_{q^{-4}},X_{i,0}^+]_{q^{-2}}-
[[X_{j,1}^+,X_{i,0}^+]_{q^{-2}},[X_{i,0}^+,X_{i,-1}^+]_{q^{2}}]=$$
$$=-q^4[[X_{i,0}^+,[X_{j,1}^+,X_{i,-1}^+]_{q^{-2}}]_{q^{-4}},X_{i,0}^+]_{q^{-2}}+$$
$$-
[X_{j,1}^+,[X_{i,0}^+,[X_{i,0}^+,X_{i,-1}^+]_{q^{2}}]_{q^{4}}]_{q^{-6}}+
q^{-2}[X_{i,0}^+,[X_{j,1}^+,[X_{i,0}^+,X_{i,-1}^+]_{q^{2}}]_{q^{-4}}]_{q^{6}}{\buildrel  XD,X3^{-1}\over =}$$
$${\buildrel  XD,X3^{-1}\over =}\!
q^{4}[[X_{i,0}^+,[X_{i,0}^+,X_{j,0}^+]_{q^{-2}}]_{q^{-4}},X_{i,0}^+]_{q^{-2}}+
q^{-2}[X_{i,0}^+,[X_{j,1}^+,[X_{i,0}^+,X_{i,-1}^+]_{q^{2}}]_{q^{-4}}]_{q^{6}}{\buildrel  S2\over =}$$
$${\buildrel  S2\over =}
q^{-2}[[[X_{j,0}^+,X_{i,0}^+]_{q^{2}},X_{i,0}^+]_{q^{4}},X_{i,0}^+]_{q^{-2}}-
(1+q^{-4})[X_{i,0}^+,[[X_{j,1}^+,X_{i,-1}^+]_{q^{-2}},X_{i,0}^+]]_{q^{6}}{\buildrel  XD\over =}$$
$${\buildrel  XD\over =}q^{-2}[[[X_{j,0}^+,X_{i,0}^+]_{q^{2}},X_{i,0}^+]_{q^{4}},X_{i,0}^+]_{q^{-2}}+
(q^{4}+1)[[[X_{j,0}^+,X_{i,0}^+]_{q^{2}},X_{i,0}^+],X_{i,0}^+]_{q^{-6}}$$
and
$$[X_{i,-1}^+,[[X_{i,1}^+,X_{i,0}^+]_{q^2},X_{j,0}^+]_{q^{-4}}]_{q^{-2}}=$$
$$=[[X_{i,-1}^+,[X_{i,1}^+,X_{i,0}^+]_{q^2}]_{q^{-4}},X_{j,0}^+]_{q^{-2}}+
q^{-4}[[X_{i,1}^+,X_{i,0}^+]_{q^2},[X_{i,-1}^+,X_{j,0}^+]_{q^{2}}]{\buildrel  XD\over =}$$
$${\buildrel  XD\over =}
q^{-6}[X_{j,0}^+,[[X_{i,1}^+,X_{i,0}^+]_{q^2},X_{i,-1}^+]_{q^{4}}]_{q^{2}}+
q^{-2}[[X_{i,1}^+,X_{i,0}^+]_{q^2},[X_{i,0}^+,X_{j,-1}^+]_{q^{-2}}]{\buildrel  X3^{1}\over =}$$
$${\buildrel  X3^{1}\over =}(1-q^{-2})(1-q^{-4})(q^{2}+q^{-2})[X_{j,0}^+,(X_{i,0}^+)^{3}]_{q^{2}}-
q^{-2}[[X_{i,0}^+,X_{j,-1}^+]_{q^{-2}},[X_{i,1}^+,X_{i,0}^+]_{q^2}]=$$
$$=(1-q^{-2})(1-q^{-4})(q^{2}+q^{-2})[X_{j,0}^+,(X_{i,0}^+)^{3}]_{q^{2}}+$$
$$-
q^{-2}[X_{i,0}^+,[X_{j,-1}^+,[X_{i,1}^+,X_{i,0}^+]_{q^2}]_{q^{4}}]_{q^{-6}}+
q^{-4}[X_{j,-1}^+,[X_{i,0}^+,[X_{i,1}^+,X_{i,0}^+]_{q^2}]_{q^{-4}}]_{q^{6}}{\buildrel  X3^{1}\over =}$$
$${\buildrel  X3^{1}\over =}\!\!(1-q^{-2})(1-q^{-4})(q^{2}+q^{-2})[X_{j,0}^+,\!(X_{i,0}^+)^{3}]_{q^{2}}+
q^{2}[X_{i,0}^+,[[X_{i,1}^+,X_{i,0}^+]_{q^2},\!X_{j,-1}^+]_{q^{-4}}]_{q^{-6}}{\buildrel  S2\over =}$$
$${\buildrel  S2\over =}
(q^{2}+q^{-2})
((1-q^{-2})(1-q^{-4})[X_{j,0}^+,(X_{i,0}^+)^{3}]_{q^{2}}-
[X_{i,0}^+,[[X_{j,-1}^+,X_{i,1}^+]_{q^2},X_{i,0}^+]]_{q^{-6}}){\buildrel  XD\over =}$$
$${\buildrel  XD\over =}(q^{2}+q^{-2})
((1-q^{-2})(1-q^{-4})[X_{j,0}^+,(X_{i,0}^+)^{3}]_{q^{2}}+q^{-4}
[[[X_{j,0}^+,X_{i,0}^+]_{q^{-2}},X_{i,0}^+],X_{i,0}^+]_{q^{6}}).
$$
It follows that 
$$q^{-2}[[[X_{j,0}^+,X_{i,0}^+]_{q^{2}},X_{i,0}^+]_{q^{4}},X_{i,0}^+]_{q^{-2}}+
(q^{4}+1)[[[X_{j,0}^+,X_{i,0}^+]_{q^{2}},X_{i,0}^+],X_{i,0}^+]_{q^{-6}}+$$
$$+(q^{2}-1)(1-q^{-4})[X_{j,0}^+,(X_{i,0}^+)^{3}]_{q^{2}}+q^{-2}
[[[X_{j,0}^+,X_{i,0}^+]_{q^{-2}},X_{i,0}^+],X_{i,0}^+]_{q^{6}}=$$
$$=
(q+q^{-1})(q^{3}+q^{-3})\big(X_{j,0}^+(X_{i,0}^+)^3-(q^2+1+q^{-2})X_{i,0}^+X_{j,0}^+(X_{i,0}^+)^2+$$
$$+(q^2+1+q^{-2})(X_{i,0}^+)^2X_{j,0}^+X_{i,0}^+-(X_{i,0}^+)^3X_{j,0}^+\big)$$
is an element of ${\Cal I}_{cte}^+(XD,X3^{1},S2)$, hence
${\Cal I}_{0}^{+}(S_{(i,j)})\subseteq{\Cal I}_{cte}^{+}(XD,X3^1,S2)$.

Thanks to corollary \dlem.\ddd\ we obtain ${\Cal I}_{cte}^{\pm}(S_{(i,j)})\subseteq{\Cal I}_{cte}^{\pm}(XD,X3^1,S2)$.

\vskip .3 truecm

\prop{\tts}
If $k=3$ and $i,j\in I_0$ are such that $a_{ij}=-3$, then 
${\Cal I}_{cte}^{\pm}(S_{(i,j)})\subseteq{\Cal I}_{cte}^{\pm}(X1_i,T3)$.
\dim
Let us start from 
$$(q^2+1)[[X_{j,0}^+,X_{i,2}^+]_{q^{3}},X_{i,0}^+]_{q^{-1}}+[[X_{j,0}^+,X_{i,1}^+]_{q^{3}},X_{i,1}^+]_q,$$
which is an element of ${\Cal I}_{cte}^{+}(T3)$, and remark that 
$$[[X_{j,0}^+,X_{i,2}^+]_{q^{3}},X_{i,0}^+]_{q^{-1}}-q^2[[X_{j,0}^+,X_{i,0}^+]_{q^{-3}},X_{i,2}^+]_{q}-(q^2-1)[X_{j,0}^+,(X_{i,1}^+)^2]$$
belongs to ${\Cal I}_{cte}^{+}(X1_i)$;
but
$$(q^2+1)(q^2-1)[X_{j,0}^+,(X_{i,1}^+)^2]+[[X_{j,0}^+,X_{i,1}^+]_{q^{3}},X_{i,1}^+]_q=$$
$$=q^4[[X_{j,0}^+,X_{i,1}^+]_{q^{-3}},X_{i,1}^+]_{q^{-1}},$$
so that
$$[[X_{j,0}^+,X_{i,1}^+]_{q^{-3}},X_{i,1}^+]_{q^{-1}}+(1+q^{-2})[[X_{j,0}^+,X_{i,0}^+]_{q^{-3}},X_{i,2}^+]_{q}$$ lies in ${\Cal I}_{cte}^{+}(X1_i,T3)$, hence 
$$[[X_{j,-3}^+,X_{i,4}^+]_{q^{3}},X_{i,1}^+]_{q^{-1}}+(1+q^{-2})[[X_{j,-3}^+,X_{i,3}^+]_{q^{3}},X_{i,2}^+]_{q}$$ 
and (applying $\tilde t_j^{-1}\tilde t_i^2$ and $q$-commuting by $X_{i,0}^+$)
$$[[[X_{j,0}^+,X_{i,2}^+]_{q^{3}},X_{i,-1}^+]_{q^{-1}}+(1+q^{-2})[[X_{j,0}^+,X_{i,1}^+]_{q^{3}},X_{i,0}^+]_{q},X_{i,0}^+]_{q^{-3}}$$
lie in ${\Cal I}_{cte}^{+}(XD,X1_i,T3)$;
but
$$[[[X_{j,0}^+,X_{i,2}^+]_{q^{3}},X_{i,-1}^+]_{q^{-1}},X_{i,0}^+]_{q^{-3}}-
q^{-2}[[[X_{j,0}^+,X_{i,2}^+]_{q^{3}},X_{i,0}^+]_{q^{-1}},X_{i,-1}^+]_{q}$$
belongs to ${\Cal I}_{cte}^{+}(X1_i)$ by lemma \cid,
$$q^{-2}[[[X_{j,0}^+,X_{i,2}^+]_{q^{3}},X_{i,0}^+]_{q^{-1}},X_{i,-1}^+]_{q}+{q^{-2}\over q^2+1}
[[[X_{j,0}^+,X_{i,1}^+]_{q^{3}},X_{i,1}^+]_{q},X_{i,-1}^+]_{q}$$
belongs to ${\Cal I}_{cte}^{+}(T3)$, 
$$[[[X_{j,0}^+,X_{i,1}^+]_{q^{3}},X_{i,1}^+]_{q},X_{i,-1}^+]_{q}+$$
$$-q^2[[[X_{j,0}^+,X_{i,1}^+]_{q^{3}},X_{i,-1}^+]_{q^{-1}},X_{i,1}^+]_{q^{-1}}-
(q^2-1)[[X_{j,0}^+,X_{i,1}^+]_{q^{3}},(X_{i,0}^+)^2]$$
belongs to ${\Cal I}_{cte}^{+}(X1_i)$ again  by lemma \cid,
$$q^2[[[X_{j,0}^+,X_{i,1}^+]_{q^{3}},X_{i,-1}^+]_{q^{-1}},X_{i,1}^+]_{q^{-1}}+{q^{2}\over q^2+1}
[[[X_{j,0}^+,X_{i,0}^+]_{q^{3}},X_{i,0}^+]_{q},X_{i,1}^+]_{q^{-1}}$$
belongs to ${\Cal I}_{cte}^{+}(T3)$ and 
$$[[[X_{j,0}^+,X_{i,0}^+]_{q^{3}},X_{i,0}^+]_{q},X_{i,1}^+]_{q^{-1}}-q^{-4}
[[[X_{j,0}^+,X_{i,1}^+]_{q^{3}},X_{i,0}^+]_{q^{5}},X_{i,0}^+]_{q^{3}}$$
belongs to ${\Cal I}_{cte}^{+}(X1_i)$ (by lemma \cid).
So we can conclude that 
$${q^{-4}\over (q^2+1)^2}[[[X_{j,0}^+,X_{i,1}^+]_{q^{3}},X_{i,0}^+]_{q^{5}},X_{i,0}^+]_{q^{3}}-{1-q^{-2}\over q^2+1}[[X_{j,0}^+,X_{i,1}^+]_{q^{3}},(X_{i,0}^+)^2]+$$
$$+(1+q^{-2})[[[X_{j,0}^+,X_{i,1}^+]_{q^{3}},X_{i,0}^+]_{q},X_{i,0}^+]_{q^{-3}}=$$
$$={(q^{2}+1+q^{-2})^2\over (q^2+1)^2}[[[X_{j,0}^+,X_{i,1}^+]_{q^{3}},X_{i,0}^+]_{q},X_{i,0}^+]_{q^{-1}}$$
lies in ${\Cal I}_{cte}^{+}(XD,X1_i,T3)$.

Then 
$[[[[X_{j,0}^+,X_{i,1}^+]_{q^{3}},X_{i,0}^+]_q,X_{i,0}^+]_{q^{-1}},X_{i,-1}^+]_{q^{3}}\in{\Cal I}_{cte}^+(X1_i,T3)$;
but
$$[[[[X_{j,0}^+,X_{i,1}^+]_{q^{3}},X_{i,0}^+]_q,X_{i,0}^+]_{q^{-1}},X_{i,-1}^+]_{q^{3}
}+$$
$$-q^4[[[[X_{j,0}^+,X_{i,1}^+]_{q^{3}},X_{i,-1}^+]_{q^{-1}},X_{i,0}^+]_{q^{-1}},X_{i,0}^+]_{q^{-3}}$$ 
belongs to ${\Cal I}_{cte}^{\pm}(X1_i)$ and since
$$(q^2+1)[[X_{j,0}^+,X_{i,1}^+]_{q^{3}},X_{i,-1}^+]_{q^{-1}}+[[X_{j,0}^+,X_{i,0}^+]_{q^{3}},X_{i,0}^+]_{q}\in
{\Cal I}_{cte}^+(T3)$$
it follows that 
$$[[[[X_{j,0}^+,X_{i,0}^+]_{q^{3}},X_{i,0}^+]_{q},X_{i,0}^+]_{q^{-1}},X_{i,0}^+]_{q^{-3}}\in{\Cal I}_{cte}^+(X1,T3),$$
so that 
${\Cal I}_{0}^{+}(S_{(i,j)})\subseteq{\Cal I}_{cte}^+(X1,T3)$.
Then
${\Cal I}_{cte}^{\pm}(S_{(i,j)})\subseteq{\Cal I}_{cte}^{\pm}(X1,T3)$.

\vskip .3 truecm
\cor{\sernul}
${\Cal I}_{cte}(S)=0$ in $\Udr_q$.
\dim 
The claim is a straightforward consequence of remark \mij\ and propositions \tds, \sds, \tts.

\vskip .3 truecm
We are now able to prove that the quantum algebra of finite type ${\Cal U}_q^{fin}$ is mapped in $\Udr_q$, which was not otherwise clear. 
\vskip .3 truecm
\ddefi{\fimbd}
Let $\phi:\U_q^{fin}\to\Udr_q$ be the $\C(q)$-homomorphism given by
$$K_i^{\pm 1}\mapsto k_i^{\pm 1}, \ \ \ E_i\mapsto X_{i,0}^+,\ \ \ 
F_i\mapsto X_{i,0}^-\ \ \ (i\in I_0).$$

\vskip .3 truecm

\rem{\finbdef}~
$\phi$ is well defined. 
\dim 
It is a straightforward consequence of corollary \sernul\ (and of relations $(CUK)$, $(CK)$, $(KX^{\pm})$, $(XXE)$).

\vskip .3 truecm
We shall now complete the study of the ideal generated by the Serre relations.

\vskip .3 truecm

\rem{\srsp}~
${\Cal I}_{cte}^{\pm}(XD)\subseteq{\Cal I}_{cte}^{\pm}(S)$. 

If $n>1$
${\Cal I}_{cte}^{\pm}(X1),{\Cal I}_{cte}^{\pm}(X2),{\Cal I}_{cte}^{\pm}(X3^1),{\Cal I}_{cte}^{\pm}(X3^{-1})\subseteq{\Cal I}_{cte}^{\pm}(S)$.
\dim
That ${\Cal I}_{cte}^{\pm}(XD)\subseteq{\Cal I}_{cte}^{\pm}(S)$ follows from proposition \cnat.\rdxd\ and from remark \mij.

That ${\Cal I}_{cte}^{\pm}(X1)\subseteq{\Cal I}_{cte}^{\pm}(S)$ is a consequence of lemma \cnat.\cyi\ and of corollary \dlem.\ddd\ (see also corollary \cnat.\misp).

Finally that ${\Cal I}_{cte}^{\pm}(X3^1)\subseteq{\Cal I}_{cte}^{\pm}(S)$ follows again from lemma \cnat.\cyi\ and from corollary \dlem.\ddd, once one notices that $(X_{\tilde n}^{(k)},i,j)=(A_{2n}^{(2)},1,2)$ implies $q_i=q$, $b_{ji\tilde d_{ij}}\neq 0$ and $b_{ii\tilde d_{i}}=[2]_q[3]_q$. 

From this it follows that ${\Cal I}_{cte}^{\pm}(X2),{\Cal I}_{cte}^{\pm}(X3^{-1})\subseteq{\Cal I}_{cte}^{\pm}(S)$ (see proposition \cnat.\rku).
\vskip .3 truecm

\cor{\ctes}
i) ${\Cal I}^{\pm}(S)={\Cal I}_{cte}^{\pm}(S)$ (see corollary \dlem.\hhh,iii) and remark \srsp);

ii) ${\Cal I}(S)=0$ in $\Udr_q$ (see corollary \sernul).

\vskip .3 truecm

\rem{\sdkt}~
If $k>1$ and $X_{\tilde n}\neq A_{2 n}$ we have that 
${\Cal I}_{cte}^{\pm}(Sk)\subseteq{\Cal I}_{cte}^{\pm}(S)
\Leftrightarrow{\Cal I}_{cte}^{\pm}(Tk)\subseteq{\Cal I}_{cte}^{\pm}(S)$.
\dim
Of course we can suppose $n>1$; then the claim depends on the fact that $({\Cal I}_{cte}^{\pm}(X1),{\Cal I}_{cte}^{\pm}(Sk))=
({\Cal I}_{cte}^{\pm}(X1),{\Cal I}_{cte}^{\pm}(Tk))$ (see remark \drmr.\xuit) and that
${\Cal I}_{cte}^{\pm}(X1)\subseteq{\Cal I}_{cte}^{\pm}(S)$ (see remark \srsp).

\vskip .3 truecm

\prop{\sdt}
${\Cal I}_{cte}^{\pm}(S2)\subseteq{\Cal I}_{cte}^{\pm}(S)$.
\dim
Let $k=2$ and $i,j\in I_0$ be such that $a_{ij}=-2$. 

Then 
$[[[X_{j,0}^+,X_{i,0}^+]_{q^{2}},X_{i,0}^+],X_{i,0}^+]_{q^{-2}}$ is an element of ${\Cal I}_{cte}^+(S)$, so that 
$${\Cal I}_{cte}^+(S)\ni[[[[X_{j,0}^+,X_{i,0}^+]_{q^{2}},X_{i,0}^+],X_{i,0}^+]_{q^{-2}},X_{i,1}^-]=$$
$$=[[[X_{j,0}^+,C^{-1}k_iH_{i,1}]_{q^{2}},X_{i,0}^+],X_{i,0}^+]_{q^{-2}}+
[[[X_{j,0}^+,X_{i,0}^+]_{q^{2}},C^{-1}k_iH_{i,1}],X_{i,0}^+]_{q^{-2}}+$$
$$+[[[X_{j,0}^+,X_{i,0}^+]_{q^{2}},X_{i,0}^+],C^{-1}k_iH_{i,1}]_{q^{-2}}=$$
$$=-C^{-1}k_i(q^{2}[[[H_{i,1},X_{j,0}^+],X_{i,0}^+]_{q^{-2}},X_{i,0}^+]_{q^{-4}}+
[[[H_{i,1},X_{j,0}^+],X_{i,0}^+]_{q^{2}},X_{i,0}^+]_{q^{-4}}+$$
$$+[[X_{j,0}^+,[H_{i,1},X_{i,0}^+]]_{q^{2}},X_{i,0}^+]_{q^{-4}}+
q^{-2}[[[H_{i,1},X_{j,0}^+],X_{i,0}^+]_{q^{2}},X_{i,0}^+]+$$
$$+q^{-2}[[X_{j,0}^+,[H_{i,1},X_{i,0}^+]]_{q^{2}},X_{i,0}^+]+
q^{-2}[[X_{j,0}^+,X_{i,0}^+]_{q^{2}},[H_{i,1},X_{i,0}^+]])=$$
$$=-C^{-1}k_i(b_{ij1}(q^2[[X_{j,1}^+,X_{i,0}^+]_{q^{-2}},X_{i,0}^+]_{q^{-4}}+
[[X_{j,1}^+,X_{i,0}^+]_{q^{2}},X_{i,0}^+]_{q^{-4}}+$$
$$+q^{-2}[[X_{j,1}^+,X_{i,0}^+]_{q^{2}},X_{i,0}^+])+ 
b_{ii1}([[X_{j,0}^+,X_{i,1}^+]_{q^{2}},X_{i,0}^+]_{q^{-4}}+$$
$$+q^{-2}[[X_{j,0}^+,X_{i,1}^+]_{q^{2}},X_{i,0}^+]+
q^{-2}[[X_{j,0}^+,X_{i,0}^+]_{q^{2}},X_{i,1}^+]))=$$
$$=-C^{-1}k_i([3]_qb_{ij1}(X_{j,1}^+(X_{i,0}^+)^2-(q^{-2}+1)X_{i,0}^+X_{j,1}^+X_{i,0}^++
q^{-2}(X_{i,0}^+)^2X_{j,1}^+)+$$
$$+[3]_qb_{ii1}(q^{-2}X_{j,0}^+X_{i,1}^+X_{i,0}^+-X_{i,1}^+X_{j,0}^+X_{i,0}^+
-q^{-2}X_{i,0}^+X_{j,0}^+X_{i,1}^++X_{i,0}^+X_{i,1}^+X_{j,0}^+)+$$
$$+b_{ii1}(-q^{-4}X_{j,0}^+[X_{i,1}^+,X_{i,0}^+]_{q^2}+[X_{i,1}^+,X_{i,0}^+]_{q^2}X_{j,0}^+))=$$
$$=-C^{-1}k_i([3]_qb_{ij1}[[X_{j,1}^+,X_{i,0}^+]_{q^{-2}},X_{i,0}^+]+$$
$$+q^{-2}[3]_qb_{ii1}[[X_{j,0}^+,X_{i,1}^+]_{q^{2}},X_{i,0}^+]+b_{ii1}[[X_{i,1}^+,X_{i,0}^+]_{q^2},X_{j,0}^+]_{q^{-4}}).$$

Now let us notice that if $X_{\tilde n}\neq A_{2 n}$ we have $b_{ii1}\neq 0$, $[X_{i,1}^+,X_{i,0}^+]_{q^2}\in{\Cal I}_{cte}^+(S)$ (see remark \srsp) and $\tilde d_j=2$, hence $b_{ij1}=0$: we can conclude that 
$[[X_{j,0}^+,X_{i,1}^+]_{q^{2}},X_{i,0}^+]$ is an element of $
{\Cal I}_{cte}^+(S)$, so that  
${\Cal I}_{cte}^{+}(T2)\subseteq{\Cal I}_{cte}^{+}(S)$ (see corollary \dlem.\ddd), which, thanks to remark \sdkt, is equivalent to ${\Cal I}_{cte}^+(S2)\subseteq{\Cal I}_{cte}^+(S)$.

On the other hand, if $X_{\tilde n}= A_{2 n}$we have $b_{ii1}=[2]_q[3]_q$, $\tilde d_j=1$, $b_{ij1}=-[2]_q$ and 
$[X_{j,1}^+,X_{i,0}^+]_{q^{-2}}+[X_{i,1}^+,X_{j,0}^+]_{q^{-2}}\in{\Cal I}_{cte}^{+}(S)$ (see remark \srsp); 
then we have that 
$$(q^2+q^{-2})[[X_{j,0}^+,X_{i,1}^+]_{q^{2}},X_{i,0}^+]+q^2[[X_{i,1}^+,X_{i,0}^+]_{q^2},X_{j,0}^+]_{q^{-4}}$$ is an element of ${\Cal I}_{cte}^{+}(S)$, that is
${\Cal I}_{cte}^+(S2)\subseteq{\Cal I}_{cte}^+(S)$.

In both cases using corollary \dlem.\ddd\ we get ${\Cal I}_{cte}^{\pm}(S2)\subseteq{\Cal I}_{cte}^{\pm}(S)$.

\vskip .3 truecm

\prop{\std}
${\Cal I}_{cte}^{\pm}(T3)\subseteq{\Cal I}_{cte}^{\pm}(S)$.
\dim
Let $k=3$ ($X_{\tilde n}^{(k)}=D_{4}^{(3)}$) and $i,j\in I_0$ be such that $a_{ij}=-3$ ($i=1$, $j=2$). 

Then 
$[[[[X_{j,0}^+,X_{i,0}^+]_{q^{3}},X_{i,0}^+]_{q},X_{i,0}^+]_{q^{-1}},X_{i,0}^+]_{q^{-3}}$ is an element of ${\Cal I}_{cte}^+(S)$, so that, recalling that $b_{ij1}=0$ and $b_{ii1}=[2]_q$,
$${\Cal I}_{cte}^+(S)\ni[[[[[X_{j,0}^+,X_{i,0}^+]_{q^{3}},X_{i,0}^+]_{q},X_{i,0}^+]_{q^{-1}},X_{i,0}^+]_{q^{-3}},X_{i,1}^-]=$$
$$=[[[[X_{j,0}^+,C^{-1}k_iH_{i,1}]_{q^{3}},X_{i,0}^+]_{q},X_{i,0}^+]_{q^{-1}},X_{i,0}^+]_{q^{-3}}+$$
$$+
[[[[X_{j,0}^+,X_{i,0}^+]_{q^{3}},C^{-1}k_iH_{i,1}]_{q},X_{i,0}^+]_{q^{-1}},X_{i,0}^+]_{q^{-3}}+$$
$$+[[[[X_{j,0}^+,X_{i,0}^+]_{q^{3}},X_{i,0}^+]_{q},C^{-1}k_iH_{i,1}]_{q^{-1}},X_{i,0}^+]_{q^{-3}}+$$
$$+[[[[X_{j,0}^+,X_{i,0}^+]_{q^{3}},X_{i,0}^+]_{q},X_{i,0}^+]_{q^{-1}},C^{-1}k_iH_{i,1}]_{q^{-3}}=$$
$$=-C^{-1}k_i(q^{3}[[[[H_{i,1},X_{j,0}^+],X_{i,0}^+]_{q^{-1}},X_{i,0}^+]_{q^{-3}},X_{i,0}^+]_{q^{-5}}+$$
$$+q
[[[H_{i,1},[X_{j,0}^+,X_{i,0}^+]_{q^{3}}],X_{i,0}^+]_{q^{-3}},X_{i,0}^+]_{q^{-5}}+$$
$$+q^{-1}[[H_{i,1},[[X_{j,0}^+,X_{i,0}^+]_{q^{3}},X_{i,0}^+]_{q}],X_{i,0}^+]_{q^{-5}}+$$
$$+q^{-3}[H_{i,1},[[[X_{j,0}^+,X_{i,0}^+]_{q^{3}},X_{i,0}^+]_{q},X_{i,0}^+]_{q^{-1}}])=$$
$$=\!\!-[2]_qC^{-1}\!k_i(q
[[[X_{j,0}^+,X_{i,1}^+]_{q^{3}},\!X_{i,0}^+]_{q^{-3}},\!X_{i,0}^+]_{q^{-5}}+
q^{-1}[[[X_{j,0}^+,X_{i,1}^+]_{q^{3}},X_{i,0}^+]_{q},X_{i,0}^+]_{q^{-5}}+$$
$$+q^{-1}[[[X_{j,0}^+,X_{i,0}^+]_{q^{3}},X_{i,1}^+]_{q},X_{i,0}^+]_{q^{-5}}+
q^{-3}[[[X_{j,0}^+,X_{i,1}^+]_{q^{3}},X_{i,0}^+]_{q},X_{i,0}^+]_{q^{-1}})+$$
$$+q^{-3}[[[X_{j,0}^+,X_{i,0}^+]_{q^{3}},X_{i,1}^+]_{q},X_{i,0}^+]_{q^{-1}})+
q^{-3}[[[X_{j,0}^+,X_{i,0}^+]_{q^{3}},X_{i,0}^+]_{q},X_{i,1}^+]_{q^{-1}});$$
then, thanks to remarks \srsp\ and \cid, we have that 
$$q
[[[X_{j,0}^+,X_{i,1}^+]_{q^{3}},X_{i,0}^+]_{q^{-3}},X_{i,0}^+]_{q^{-5}}+
q^{-1}[[[X_{j,0}^+,X_{i,1}^+]_{q^{3}},X_{i,0}^+]_{q},X_{i,0}^+]_{q^{-5}}+$$
$$+q^{-3}[[[X_{j,0}^+,X_{i,1}^+]_{q^{3}},X_{i,0}^+]_{q^5},X_{i,0}^+]_{q^{-5}}+
q^{-3}[[[X_{j,0}^+,X_{i,1}^+]_{q^{3}},X_{i,0}^+]_{q},X_{i,0}^+]_{q^{-1}})+$$
$$+q^{-5}[[[X_{j,0}^+,X_{i,1}^+]_{q^{3}},X_{i,0}^+]_{q^5},X_{i,0}^+]_{q^{-1}})+
q^{-7}[[[X_{j,0}^+,X_{i,0}^+]_{q^{3}},X_{i,0}^+]_{q^5},X_{i,1}^+]_{q^{3}}=$$
$$=q^{-3}(q^2+q^{-2})[3]_q[[[X_{j,0}^+,X_{i,1}^+]_{q^{3}},X_{i,0}^+]_{q},X_{i,0}^+]_{q^{-1}}$$
belongs to ${\Cal I}_{cte}^+(S)$; then so does 
$$[[[[X_{j,0}^+,X_{i,1}^+]_{q^{3}},X_{i,0}^+]_{q^{-1}},X_{i,0}^+]_{q},X_{i,-1}^-]=$$
$$=\Big[\Big[\Big[X_{j,0}^+,{Ck_i-C^{-1}k_i^{-1}\over q-q^{-1}}\Big]_{q^{3}},X_{i,0}^+\Big]_{q^{-1}},X_{i,0}^+\Big]_{q}+$$
$$+[[[X_{j,0}^+,X_{i,1}^+]_{q^{3}},k_i^{-1}H_{i,-1}]_{q^{-1}},X_{i,0}^+]_{q}+
[[[X_{j,0}^+,X_{i,1}^+]_{q^{3}},X_{i,0}^+]_{q^{-1}},k_i^{-1}H_{i,-1}]_{q}=$$
$$=k_i^{-1}([3]_{q}C^{-1}[[X_{j,0}^+,X_{i,0}^+]_{q},X_{i,0}^+]_{q^{3}}+$$
$$-q^{-1}[[[H_{i,-1},[X_{j,0}^+,X_{i,1}^+]_{q^{3}}],X_{i,0}^+]_{q^{3}}-
q[H_{i,-1},[[X_{j,0}^+,X_{i,1}^+]_{q^{3}},X_{i,0}^+]_{q^{-1}}])=$$
$$=C^{-1}k_i^{-1}([3]_{q}[[X_{j,0}^+,X_{i,0}^+]_{q^{3}},X_{i,0}^+]_{q}-
q^{-1}[2]_{q}[[[X_{j,0}^+,X_{i,0}^+]_{q^{3}},X_{i,0}^+]_{q^{3}}+$$
$$-q[2]_{q}[[X_{j,0}^+,X_{i,0}^+]_{q^{3}},X_{i,0}^+]_{q^{-1}}-q[2]_{q}[[X_{j,0}^+,X_{i,1}^+]_{q^{3}},X_{i,-1}^+]_{q^{-1}})=$$
$$=-C^{-1}k_i^{-1}([[X_{j,0}^+,X_{i,0}^+]_{q^{3}},X_{i,0}^+]_{q}+
q[2]_{q}[[X_{j,0}^+,X_{i,1}^+]_{q^{3}},X_{i,-1}^+]_{q^{-1}});$$
hence 
${\Cal I}_{0}^+(T3)\subseteq{\Cal I}_{cte}^+(S)$ and, using corollary \dlem.\ddd, ${\Cal I}_{cte}^{\pm}(T3)\subseteq{\Cal I}_{cte}^{\pm}(S)$.

\vskip .3 truecm

\cor{\cstd}
${\Cal I}_{cte}^{\pm}(S3)\subseteq{\Cal I}_{cte}^{\pm}(S)$.
\dim
It follows from remark \sdkt\ and from proposition \std.
\vskip .3 truecm

\cor{\sug}
If $n>1$ ${\Cal I}_{cte}^{\pm}(XD,X1,X2,X3^{\pm1},SUL,S2,S3)={\Cal I}_{cte}^{\pm}(S)$.
\dim
It follows from corollary \sernul, remarks \mij\ and \srsp, proposition \sdt\ and corollary \cstd.
\vskip .3 truecm

\rem{\inttr}~
In $\Udr_q(D_4^{(3)})$ we have $[[[X_{j,0}^+,X_{i,1}^+]_{q^{3}},X_{i,0}^+]_{q},X_{i,0}^+]_{q^{-1}}=0$.
\dim
See the proof of proposition \std.
\vskip .3 truecm 
\teo{\ssrr}
i) ${\Cal U}_q^{Dr}(A_1^{(1)})=\tilde{\Cal U}_q^{Dr}(A_1^{(1)})/{\Cal I}_{cte}(X1)$;

ii) ${\Cal U}_q^{Dr}(A_{2}^{(2)})=\tilde{\Cal U}_q^{Dr}(A_{2}^{(2)})/{\Cal I}_{cte}(X3^1)$;

iii) ${\Cal U}_q^{Dr}(X_{\tilde n}^{(k)})=\tilde{\Cal U}_q^{Dr}(X_{\tilde n}^{(k)})/{\Cal I}_{cte}(S)$ if 
$n>1$ (that is $X_{\tilde n}^{(k)}\neq A_1^{(1)},A_{2}^{(2)}$).
\dim
The claims follow from theorem \cnat.\tldg\ and corollary \sug.

\vskip .3 truecm

\cor{\sspz}
Let: $\alg$ be a $\C(q)$-algebra, $t_i^{(\alg)}$ ($i\in I_0$) be $\C(q)$-automorphisms of $\alg$, $\Omega^{(\alg)}$ be a $\C$-anti-linear anti-automorphism of $\alg$, $\tilde f:\tilde{\Cal U}^{Dr}_q(X_{\tilde n}^{(k)})\rightarrow\alg$ be a homomorphism of $\C(q)$-algebras such that $\tilde f\compo \tilde t_i=t_i^{(\alg)}\compo\tilde f$ $\forall i\in I_0$ and $\tilde f\compo \tilde\Omega=\Omega^{(\alg)}\compo\tilde f$. 

If:

i) $\tilde f({\Cal I}_0^+(X1))=0$ in case $X_{\tilde n}^{(k)}=A_1^{(1)}$;

ii) $\tilde f({\Cal I}_0^+(X3^1))=0$ in case $X_{\tilde n}^{(k)}=A_{2}^{(2)}$; 

iii) $\tilde f({\Cal I}_0^+(S))=0$ in case $X_{\tilde n}^{(k)}\neq A_1^{(1)},A_{2}^{(2)}$;

then:

$\tilde f$ induces $f:\Udr_q(X_{\tilde n}^{(k)})\rightarrow\alg$ and $f\compo t_i=t_i^{(\alg)}\compo f$ $\forall i\in I_0$, $f\compo \Omega=\Omega^{(\alg)}\compo f$. 
\dim 
Since the hypotheses imply that $ker(\tilde f)$ is a $t_i$-stable ($\forall i\in I_0$), 
$\tilde\Omega$-stable ideal of $\Udr_q(X_{\tilde n}^{(k)})$, the claim is an immediate consequence of theorem \ssrr\ and of corollary \dlem.\ddd.

\vskip .3 truecm
\rem{\thsdg}~
It is useful to compare the results of this section with those of section \cnat. The simplification of the relations given in section \cnat\ (theorem \cnat.\tldg\ and corollary \cnat.\hmbdefnt) provides a minimal set of relations of lowest ``degree'' (where the degree of $X_{i_1,r_1}\cdot...\cdot X_{i_h,r_h}$ is meant to be $h$); this minimality in degree can be often useful, in spite of the appearence more complicated of relations like $(S2^{\pm})$ with respect to the simple and familiar Serre relations. On the other hand the advantage of the Serre relations is evident in all the cases, like the application of theorem \ssrr\ and corollary  \sspz\ given in section \hmfr, when the Serre relations play a central role: to this aim recall that the Serre relations are the minimal degree relations defining the positive part of $\U_q^{fin}$ (see definition \fimbd\ and remark \finbdef, and recall \lusz).

\vskip .5 truecm

{\bf{\hmfr.\  The\ HOMOMORPHISM\ $\psi $\ from}}\ $\Udr_q$\ {\bf to}\ $\U_q^{DJ}$.
\vskip .5 truecm

This section is devoted to exhibit a homomorphism $\psi:\Udr_q\to\U_q^{DJ}$ and to prove that it is surjective. 

\vskip .3 truecm

\nota{\defo}
In the following $o:I_0\to\{\pm1\}$ will be a map such that: 

a) $a_{ij}\neq 0\Rightarrow o(i)o(j)=-1$ (see \beck\ for the untwisted case);

b) in the twisted case different from $A_{2n}^{(2)}$ $ a_{ij}=-2\Rightarrow o(i)=1$.
\vskip .3truecm
\rem{\defx}~
A map $o$ as in notation \defo\ exists and is: 

i) determined up to a sign, in the untwisted case and in cases $A_{2n}^{(2)}$ and $D_4^{(3)}$; 

ii) uniquely determined, in cases $A_{2n-1}^{(2)}$ and $E_6^{(2)}$.

\vskip .3truecm

\ddefi{\defpsi}
Let
$\tilde\psi=\tilde\psi_{X_{\tilde n}^{(k)}}:\tilde{\Cal U}^{Dr}_q(X_{\tilde n}^{(k)})\to\U_q^{DJ}(X_{\tilde n}^{(k)})$ be the $\C(q)$-algebra homomorphism 
defined on the generators as follows: 
$$C^{\pm1}\mapsto K_{\d}^{\pm1},\ \ \ 
k_i^{\pm1}\mapsto K_i^{\pm1}\ \ (i\in I_0),$$
$$X_{i,\tilde d_i r}^+\mapsto o(i)^rT_{\lambda_i}^{-r}(E_i),\ \ \ 
X_{i,\tilde d_i r}^-\mapsto o(i)^rT_{\lambda_i}^r(F_i)\ \ (i\in I_0, r\in\Z),$$
$$H_{i,\tilde d_i r}\mapsto \begin{cases}  o(i)^rE_{(\tilde d_i r\d,i)}&\roman{if\ } r>0\cr
o(i)^rF_{(-\tilde d_i r\d,i)}&\roman{if\ } r<0\end{cases} \ \ (i\in I_0, r\in\Z\setminus\{0\}).$$

\vskip .3truecm

\prop{\semi}
i) $\tilde\psi$ is well defined; 

ii) $\tilde\psi\compo\tilde\Omega=\Omega\compo\tilde\psi$;

iii) $\tilde\psi\compo\tilde t_i=T_{\lambda_i}\compo\tilde\psi$ $\forall i\in I_0$; 

iv) $\tilde\psi\compo\tilde\phi_i=\varphi_i\compo\tilde\psi$ $\forall i\in I_0$;

v) $\tilde\psi\compo\tilde\Theta=(\Omega\Xi T_1)\compo\tilde\psi$ in cases $A_1^{(1)}$ and $A_2^{(2)}$.
\dim

i) The relations $(ZX^{\pm})$,  $(CUK)$, $(CK)$ and $(KX^{\pm})$ are obviuosly preserved by $\tilde\psi$; also $(XX)$ (see \beck\ and \damcina) and $(HXL^{\pm})$ hold in $\U_q$: it is enough to notice that $\forall i,j\in I_0$, $\forall r\in\Z$ such that $max\{\tilde d_i, \tilde d_j\}|r$ we have $b_{ijr}=(o(i)o(j))^{\tilde r}x_{ijr}$ where $\tilde r={r\over max\{\tilde d_i,\tilde d_j\}}$ and
$$x_{ijr}=\begin{cases} (o(i)o(j))^r{[ra_{ij}]_{q_i}\over r}& {\roman {if}\ }k=1,\ {\roman {or}}\ X_{\tilde n}^{(k)}=A_{2n}^{(2)} \ {\roman {and}}\ (i,j)\neq(1,1)\\
{[2r]_q\over r}(q^{2r}+(-1)^{r-1}+q^{-2r})& {\roman {if}\ } (X_{\tilde n}^{(k)},i,j)=(A_{2n}^{(2)},1,1)\\
(o(i)o(j))^{\tilde r}{[ra_{ij}^s]_q\over\tilde r[d_i]_q}& {\roman {otherwise}\ }
\end{cases} $$
with $a_{ij}^s=max\{a_{ij},a_{ji}\}$ (see \damcina). 

ii), iii), iv) and v) are trivial. 
\vskip .3 truecm
\teo{\bdf}
Let $X_{\tilde n}^{(k)}$ be different from $A_1^{(1)}$ and $A_2^{(2)}$. Then $\tilde\psi$ induces $$\psi=\psi_{X_{\tilde n}^{(k)}}:\Udr_q(X_{\tilde n}^{(k)})\to\U_q^{DJ}(X_{\tilde n}^{(k)}).$$ 
\dim
Thanks to corollary  \serrel.\sspz, iii) and to proposition \semi, i)-iii) it is enough to prove 
that 
$\tilde \psi({\Cal I}_0^+(S))=0$; but this is obvious since $\tilde \psi({\Cal I}_0^+(S))$ is the ideal generated by the (``positive'') Serre relations. 

\vskip .3 truecm

In order to prove that $\psi$ is well defined also in the remaining cases  
we propose two different arguments: a direct one, requiring just some simple commutation relations in $\U_q(A_1^{(1)})$ and $\U_q(A_2^{(2)})$ (see lemma \scr); and an argument using the injections $\varphi_i$ (see \drjq.\brgrg).
\vskip .3 truecm

\lem{\scr}
In $\U_q(A_1^{(1)})$ we have that:

i) $E_{\delta+\alpha_1}E_1=q^2E_1E_{\delta+\alpha_1}$.

In $\U_q(A_2^{(2)})$ we have that:

ii) $E_{\delta+\alpha_1}E_1-q^2E_1E_{\delta+\alpha_1}=-[4]_qE_{\delta+2\alpha_1}$;

iii) $E_{\delta+2\alpha_1}E_1=q^4E_1E_{\delta+2\alpha_1}$;

iv) $q^{-3}E_{\delta+\alpha-1}E_1^2-(q+q^{-1})E_1E_{\delta+\alpha-1}E_1+
q^3E_1^2E_{\delta+\alpha-1}=0$.
\dim
i) is an immediate application of the Levendorskii-Soibelman formula (see \levsb\ and \damcina);

ii) see \damcina; 

iii) is an immediate application of the Levendorskii-Soibelman formula (see \levsb\ and \damcina);

iv) follows from ii) and iii).
\vskip .3truecm

\teo{\bendef}
$\tilde\psi$ induces $\psi=\psi_{X_{\tilde n}^{(k)}}:\Udr_q(X_{\tilde n}^{(k)})\to\U_q^{DJ}(X_{\tilde n}^{(k)})$.
\dim
Thanks to corollary  \serrel.\sspz, i) and ii), to proposition \semi\ and to theorem \bdf\ it is enough to notice 
that $\tilde\psi({\Cal I}_0^+(X1))=0$ in case $A_1^{(1)}$ and $\tilde\psi({\Cal I}_0^+(X3^1))=0$ in case $A_2^{(2)}$; but this is an immediate consequence of lemma \scr, i) and iv).

Otherwise: 

Let $h=1,2$, $X_{\tilde n}^{(k)}=A_4^{(2)}$,  $i=\begin{cases}  2&{\roman{if}}\ h=1\cr
1&{\roman{if}}\ h=2\end{cases} $ and consider the following well defined diagram:
$$\xymatrix{&\tilde{\Cal U}^{Dr}_q(A_h^{(h)})\ar[rr]^{\tilde\psi_{A_h^{(h)}}}\ar[d] &{} &{\Cal{U}}_q^{DJ}(A_h^{(h)})\ar[d]^{\varphi_i}\\
&\Udr_q(A_h^{(h)})\ar[r]_{\phi_i} &\Udr_q(X_{\tilde n}^{(k)})\ar[r]_{\psi_{X_{\tilde n}^{(k)}}} &{\Cal{U}}_q^{DJ}(X_{\tilde n}^{(k)})}$$

Without loss of genarality we can suppose this diagram to be commutative, by choosing 
$o_{A_h^{(h)}}:1\mapsto o_{X_{\tilde n}^{(k)}}(i)$.

Then $\tilde\psi_{A_h^{(h)}}$ factors through $\Udr_q(A_h^{(h)})$ (that is $\psi_{A_h^{(h)}}$ is well defined) because $\varphi_i$ is injective (see remark \drjq.\brgrg).

\vskip .3truecm

\rem{\ttst}~
i) $\psi\compo\Omega=\Omega\compo\psi$;

ii) $\psi\compo t_i=T_{\lambda_i}\compo\psi$ $\forall i\in I_0$; 

iii) $\psi\compo\phi_i=\varphi_i\compo\psi$ $\forall i\in I_0$;

iv) $\psi\compo\Theta=(\Omega\Xi T_1)\compo\psi$ in cases $A_1^{(1)}$ and $A_2^{(2)}$;

v) $\psi\compo\phi=\varphi$.

\dim

i)-iv) follow from proposition \semi\ and from theorem \bendef.

v) follows from remark \serrel.\finbdef\ and theorem \bendef.
\vskip 1truecm

\cor{\injfin}
$\phi$ is injective.
\dim
It follows from remark \ttst,v), since $\varphi$ is injective (see remark \drjq.\brgrg).

\vskip .3 truecm
We shall now give a proof of the surjectivity of $\psi$.

\vskip .3truecm
\rem{\ezo}~
By the definition of $\psi$ it is obvious that $E_i,F_i,K_i^{\pm1}$ are in the image of $\psi$ for all $i\in I_0$. Moreover, since $K_{\delta}^{\pm 1}$ is in the image of $\psi$, also $K_0^{\pm 1}$ is. 
But by remark \hmfr.\ttst\ $Im(\psi)$ is $\Omega$-stable, so it contains $E_0$ if and only if it contains $F_0$. Then it is enough to prove that $E_0\in Im(\psi)$.
\vskip .3truecm

In the next theorem 
it will be used that for ($i\in I_0$) $Im(\psi)$ is a $T_{\lambda_i}$-stable subalgebra of $\U_q^{DJ}$ containing $E_j,F_j,K_{\tilde j}^{\pm1}$ ($j\in I_0$, $\tilde j\in I$) (see remark \hmfr.\ttst); in particular $\U_{q,\pm\alpha}^{DJ,\pm}\subseteq Im(\psi)$ for all $\alpha=\sum_{i\in I_0}m_i\alpha_i$.
\vskip .3truecm

\teo{\eso}
$\psi:\Udr_q\to\U_q^{DJ}$ is surjective.
\dim
Let $\theta=\delta-\alpha_0=\sum_{i\in I_0}r_i\alpha_i$. Remark that there exists $i\in I_0$ such that 
either
$\tilde d_i=r_i=1$ (recall that $\theta$ is a root) or $\tilde d_i=1$, $r_i=2$, $a_{i0}\neq 0$ (so that in particular 
$\alpha_0+\alpha_i$ and $\theta-\alpha_i$ are roots). Choose such an $i\in I_0$ and let $\tilde\theta=\theta-(r_i-1)\alpha_i$: $\tilde\theta$ is a root.

Let $\lambda_i=\tau_i
s_{i_1}\cdot...\cdot s_{i_N}$ (with ($l(\lambda_i)=N$, $\tau_i\in{\Cal T}$): then $\lambda_i(\tilde\theta)=\tilde\theta-\delta<0$, so that there exists $h$ such that $s_{i_N}\cdot...\cdot s_{i_{h+1}}(\alpha_{i_h})=\tilde\theta$, and we have that $f=T_{i_N}^{-1}\cdot...\cdot T_{i_{r+1}}^{-1}(F_{i_r})\in\U_{q,-\tilde\theta}^{DJ,-}\subseteq Im(\psi)$. 

Since $Im(\psi)$ is $T_{\lambda_i}$-stable we have that $T_{\lambda_i}(f)=-\tau_iT_{i_1}\cdot...\cdot T_{i_{r-1}}(K_{i_r}^{-1}E_{i_r})=-K_{\delta-\tilde\theta}^{-1}e$ (hence $e$) belongs to $Im(\psi)$, with $e\in\U_{q,\delta-\tilde\theta}^{DJ,+}$.

If $r_i=1$, $\delta-\tilde\theta=\alpha_0$ and the claim follows ($e=E_0\in Im(\psi)$).

If $r_i=2$ then $e\in\U_{q,\alpha_0+\alpha_i}^{DJ,+}$: remark that if we are not in case $A_{2n}^{(2)}$, since $l(s_i\lambda_i)=l(\lambda_i)+1$ and $a_{i0}=-1$, $T_i(e)\in\U_{q,\alpha_0}^{DJ,+}$, so that $e=T_i^{-1}(E_0)$; on the other hand in case $A_{2n}^{(2)}$, since $l(s_0\lambda_i)=l(\lambda_i)-1$ and $a_{0i}=-1$, $T_0^{-1}(e)\in\U_{q,\alpha_i}^{DJ,+}$, so that $e=T_0(E_i)$. In both cases $e=-[E_0,E_i]_{q_i^{a_{i0}}}$. 
Commuting $e$ with $F_i(\in Im(\psi))$ we get $Im(\psi)\ni[-[E_0,E_i]_{q_i^{a_{i0}}},F_i]=q_i^{a_{i0}}[[E_i,F_i],E_0]_{q_i^{-a_{i0}}}=[a_{i0}]_{q_i}K_iE_0$, which concludes the proof. 

\vskip .3 truecm
Proving that $\psi$ is injective requires further analysis.

\vskip 2 truecm
{\bf{\ntt.\  APPENDIX: NOTATIONS}}.
\vskip .5 truecm
In this appendix, in order to make it easier for the reader to follow the exposition, most of the notations defined in the paper are collected, with the indication of  the point where they are introduced and eventually characterized.

The present list includes neither the notations related to the definition and the structure of the Drinfeld-Jimbo presentation of the quantum algebras, since they are all given synthetically in section \drjq\ where they can be easily consulted, nor the notations introduced in definition \drrl.\drcpjdef, beacuse there is no reference to them outside section \drrl.

Also the relations listed in proposition \drrl.\drdef\ are not redefined in this appendix, but for some of them other descriptions proposed and used throughout  the paper are here recalled. 

Analogously it seems useful to gather the remaining notations, which are spread out in the paper.

\vskip .3 truecm

{\bf\underbar{Dynkin diagrams, root and weight lattices}}:
\begin{align*}
\tilde\Gamma = {\phantom{:}} & {\roman{(indecomposable)\ Dynkin\ diagram\ of\ finite\ type}}&\prel\cr
\tilde I = {\phantom{:}} & {\roman{set\ of\ vertices\ of}}\ \tilde\Gamma&\prel\cr
\tilde n  = {\phantom{:}}  & \#\tilde I&\prel\cr
\tilde A  = {\phantom{:}}  &{\roman{Cartan\ matrix\ of}}\ \tilde\Gamma&\prel \cr
\chi  = {\phantom{:}}  & {\roman{automorphism\ of}}\ \tilde\Gamma&\prel \cr
k  = {\phantom{:}} & 
o(\chi)&\prel\cr
I_0  = {\phantom{:}}  & \tilde I/\chi&\prel\cr
n  = {\phantom{:}} & \# I_0&\prel\cr
\bar{}  : {\phantom{=}} & \tilde I\to I_0\ {\roman{natural\ projection}}&\prel\cr
\tilde{} : {\phantom{=}}& I_0\to \tilde I\ {\roman{section}}\ (a_{ij}\neq 0\Rightarrow\tilde a_{\tilde i,\tilde j}\neq 0)&\drrl.\resp,\ (\drrl.\coefeq)\cr
I  = {\phantom{:}} & I_0\cup\{0\}&\prel\cr
\Gamma  = {\phantom{:}} & {\roman{Dynkin\ diagram\ of\ affine\ type\ with\ set\ of\ vertices\ }}I&\prel\cr
\Gamma_0  = {\phantom{:}} & {\roman{Dynkin\ subdiagram\ of\ }}\Gamma\ {\roman{with\ set\ of\ vertices\ }}I_0&\prel\cr 
A = {\phantom{:}}  &(a_{ij})_{i,j\in I}\ {\roman{Cartan\ matrix\ of}}\ \Gamma & \prel\cr
A_0  = {\phantom{:}}  &(a_{ij})_{i,j\in I_0}\ {\roman{Cartan\ matrix\ of}}\ \Gamma_0 & \prel\cr
d_i : {\phantom{=}} &   min\{d_i|i\in I\}=1,\ diag(d_i|i\in I)A\ {\roman{symmetric}}&\prle\cr 
\tilde d_i = {\phantom{:}} &  \begin{cases} 
1&{\roman{if}}\ k=1\ {\roman{or}}\ X_{\tilde n}^{(k)}=A_{2n}^{(2)}\cr d_i&{\roman{otherwise;}}\end{cases} &\prle\cr
\tilde d_{i,j}  = {\phantom{:}} & max\{\tilde d_i,\tilde d_j\}  &\drrl.\dij\cr
\tilde d = {\phantom{:}} & max\{\tilde d_i|i\in I_0\} &\drrl.\didi\cr
\tilde Q= {\phantom{:}} & \Z^{\tilde I}=\oplus_{i^{\prime}\in\tilde I}\Z\tilde\alpha_{i^{\prime}}&\prel\cr
(\tilde\alpha_{i^{\prime}}|\tilde\alpha_{j^{\prime}}): {\phantom{=}} &W-{\roman{invariant}},(\tilde\alpha_{i^{\prime}}|\tilde\alpha_{i^{\prime}})={2k\over\tilde d}\ {\roman{if}}\ \tilde\alpha_{i^{\prime}}\ {\roman{is\ short}}&\drrl.\didi\cr
Q= {\phantom{:}} & \Z^I=\oplus_{i\in I}\Z\alpha_i&\prel\cr
Q_0= {\phantom{:}} & \oplus_{i\in I_0}\Z\alpha_i&\prel\cr
(\alpha_i|\alpha_j)= {\phantom{:}}  &d_ia_{ij}&\prle\cr
\delta: {\phantom{=}} &\delta\in Q, \delta-\alpha_0\in Q_0, (\delta|Q)=0&\prle\cr
r_i: {\phantom{=}} &\delta=\sum_{i\in I}r_i\alpha_i&\prel,\prle\cr
\theta= {\phantom{:}}  &\delta-\alpha_0&\prel\cr
\lambda_i:{\phantom{=}} &(\lambda_i|\alpha_j)=\tilde d_i\delta_{ij}\ (i,j\in I_0)&\prle
\end{align*}
{\bf\underbar{Other notations}}:
\begin{align*}
\omega   ={\phantom{:}} & {\roman{primitive\ }}k^{{\roman{th}}}\ {\roman{root\ of}}\ 1 &\drrl.\didi\cr
q_i   ={\phantom{:}}  &   q^{d_i}&\drjq.\sbal\cr
I_{\Z}  ={\phantom{:}}  &   \{(i,r)\in I_0\times\Z|\tilde d_i|r\}&\drrl.\iz\cr
\varepsilon ={\phantom{:}}  &  \pm1 &\drrl.\drcpjdef,\drrl.\drdef\cr
b_{ijr}  ={\phantom{:}}  &\begin{cases}  0&{\roman{if}}\ \tilde d_{i,j}\not|r\cr
{[2r]_q(q^{2r}+(-1)^{r-1}+q^{-2r})\over r}&{\roman{if}}\ (X_{\tilde n}^{(k)},i,j)=(A_{2n}^{(2)},1,1)\cr
{[\tilde r a_{ij}]_{q_i}\over\tilde r}&{\roman{otherwise,\ with\ }}
\tilde r={r\over \tilde d_{i,j}}.\end{cases} & \drrl.\bri  \cr
[a,b]_u ={\phantom{:}}  &   ab-uba &\drrl.\vcm\cr
\1=\1_l  ={\phantom{:}}  &  (1,...,1)\in\Z^l &\dautemb.\zl\cr
\{e_1,...,e_l\}  ={\phantom{:}}  &  {\roman{canonical\ basis\ of\ }}\Z^l&\dautemb.\zl \cr
\overline{(r_1,...,r_l)}   ={\phantom{:}}&  (r_l,...,r_1) &\dautemb.\zl\cr
o  :{\phantom{=}}  &  I_0\to\{\pm1\},\ a_{ij}\neq 0\Rightarrow o(i)o(j)=-1,&\hmfr.\defo\cr
&{\roman {if}}\ k\neq 1\ {\roman{and}}\ X_{\tilde n}^{(k)}\neq A_{2n}^{(2)}\ a_{ij}=-2\Rightarrow o(i)=1
\end{align*}
{\bf\underbar{Generators of the $\C(q)$-algebras}}: 
\begin{align*}
{\Cal G}= &\{{\Cal {C}}^{\pm 1},{\Cal {K}}_{i^{\prime}}^{\pm 1},{\Cal {X}}_{i^{\prime},r}^{\pm},{\Cal {H}}_{i^{\prime},s}|i^{\prime}\in\tilde I,r,s\in\Z,s\neq 0\}&\drrl.\drcpjdef\cr
G= &\{C^{\pm 1},K_{i}^{\pm 1},X_{i,r}^{\pm},H_{i,s}|i\in I_0,r,s\in\Z,s\neq 0\} &\drrl.\zig,\drrl.\drdef\cr
G^{\prime}=&\{C^{\pm 1},K_{i}^{\pm 1},X_{i,r}^{\pm},H_{i,s}|i\in I_0,r,s\in\tilde d_i\Z,s\neq 0\} &\drrl.\zig,\drrl.\greq\cr
\bar G= &\{C^{\pm 1},K_{i}^{\pm 1},X_{i,r}^{\pm}|i\in I_0,r\in\tilde d_i\Z\} &\drcmp.\tilu\cr
G^+= & \{X_{i,r}^+|(i,r)\in I_0\times \Z\}&\drcmp.\tilpm\cr
G^-= & \{X_{i,r}^-|(i,r)\in I_0\times \Z\}&\drcmp.\tilpm\cr
G^{\prime,+}=& \{X_{i,r}^+|(i,r)\in I_{\Z}\}&\drcmp.\rtilpm\cr
G^{\prime,-}=& \{X_{i,r}^-|(i,r)\in I_{\Z}\}&\drcmp.\rtilpm
\end{align*}
{\bf\underbar{Relations in the $\C(q)$-algebras}}:
\begin{align*}
({\Cal {DR}})={\phantom{:}}&({\Cal{Z}},{\Cal{ZX}},{\Cal{C}},{\Cal{KK}},{\Cal{KX}},{\Cal{KH}},{\Cal{XX}},{\Cal{HX}},{\Cal{HH}},&\drrl.\drcpjdef\cr
&{\Cal{XFG}},{\Cal{X}}3^{\varepsilon},{\Cal{S}},{\Cal{XP}})\cr
(DR)={\phantom{:}}&(ZX,ZH,CUK,CK,KX,KH,XX,HX,HH,&\drrl.\drdef\cr
&XD,X1,X2,X3^{\varepsilon},SUL,S2,S3)\cr
&(ZX^{\pm}),(ZH)&\drrl.\zig\cr
(HX^{\pm}):{\phantom{=}}&[H_{i,r},X_{j,s}^{\pm}]=\begin{cases}  0&{{\roman{if}}\ }\tilde d_j\not|r\cr
\pm b_{ijr}C^{{r\mp|r|\over 2}}X_{j,r+s}^{\pm}&{{\roman{if}}\ }\tilde d_j|r\end{cases} &\drrl.\greq\cr
(XD^{\pm}):{\phantom{=}}&M_{(2)}^{\pm}((i,\tilde d_ir),(j,\tilde d_js))=0 &\drmr.\tcalb\cr
(X1^{\pm}):{\phantom{=}}&\sum_{\sigma\in\sy_2}\sigma.M_i^{\pm}(\tilde d_ir)=0 &\drmr.\tcalb\cr
(X2^{\pm}):{\phantom{=}}& \sum_{\sigma\in\sy_2}\sigma.M_{(2,2)}^{\pm}(r)=0 &\drmr.\tcalb\cr
(X3^{\varepsilon,\pm}):{\phantom{=}}& \sum_{\sigma\in\sy_3}\sigma.M_{(3)}^{\varepsilon,\pm}(r)=0 &\drmr.\tcalb\cr
(S(UL)^{\pm}):{\phantom{=}}&\sum_{\sigma\in\sy_{1-a_{ij}}}\sigma.X_{i,j;1-a_{ij};1}^{\pm}(r;s)=0 &\drmr.\scpm,\drmr.\relser\cr
&\sum_{\sigma\in\sy_{1-a_{ij}}}\sigma.M_{i,j;1-a_{ij};1}^{\pm}(r;s)=0 &\drmr.\sercol,\drmr.\sercom\cr
&\sum_{\sigma\in\sy_{1-a_{ij}}}\sigma.M_{i,j;1-a_{ij};1}^{\pm}(\tilde d_ir;\tilde d_js)=0 &\drmr.\tcalb
\end{align*}
\begin{align*}
&\sum_{\sigma\in\sy_{1-a_{ij}}}\sigma.X_{i,j;1-a_{ij};1}^{\pm}(\tilde d_ir;\tilde d_js)=0 &\drmr.\tcalb\cr
(S2^+):{\phantom{=}}&\sum_{\sigma\in\sy_2}\sigma.\big((q^2+q^{-2})[[X_{j,s}^{+},X_{i,r_1+1}^{+}]_{q^2},X_{i,r_2}^{+}]+&\drmr.\sercod\cr
&+q^2[[X_{i,r_1+ 1}^{+},X_{i,r_2}^{+}]_{q^{ 2}},X_{j,s}^{+}]_{q^{- 4}}\big)=0 \cr
&\sum_{\sigma\in\sy_2}\sigma.\big((q^2+q^{-2})[[X_{j,s}^{+},X_{i,r_1}^{+}]_{q^{- 2}},X_{i,r_2+1}^{+}]+&\drmr.\sercod\cr
&+[X_{j,s}^{+},[X_{i,r_2+1}^{+},X_{i,r_1}^{+}]_{q^{ 2}}]_{q^{-4}}\big)=0\cr
&\sum_{\sigma\in\sy_2}\sigma.\big([[X_{j,s}^{+},X_{i,r_1+1}^{+}]_{q^{- 2}},X_{i,r_2}^{+}]+&\drmr.\sercod\cr
&-q^{2}[X_{i,r_1+1}^{+},[X_{j,s}^{+},X_{i,r_2}^{+}]_{q^{-2}}]_{q^{-4}}\big)=0\cr
(S3^{\pm}):{\phantom{=}}& \sum_{\sigma\in\sy_2}\!\!\sigma.\big((q^2+q^{-4})[[X_{j,s}^{\pm},\!X_{i,r_1\pm 2}^{\pm}]_{q^3},\!X_{i,r_2}^{\pm}]_{q^{- 1}}+&\drmr.\sercod\cr
&+
(1-q^{-2}+q^{-4})[[X_{j,s}^{\pm},\!X_{i,r_1\pm 1}^{\pm}]_{q^3},\!X_{i,r_2\pm 1}^{\pm}]_q+\cr
&+q^2[[X_{i,r_1\pm 2}^{\pm},X_{i,r_2}^{\pm}]_{q^{ 2}}+[X_{i,r_2\pm 1}^{\pm},X_{i,r_1\pm 1}^{\pm}]_{q^{ 2}},X_{j,s}^{\pm}]_{q^{- 6}}\big)=0\cr
(Sk^{\pm}):{\phantom{=}}&\sum_{\sigma\in\sy_2}\sigma.\!\!\!\!\sum_{{u,v\geq 0\atop u+v=-1-a_{ij}}}\!\!\!\!\!\!\!\!q^{v-u}X_{i,j;2;-a_{ij}}^{\pm}(r_1\pm v,r_2\pm u;s)=0 &\drmr.\scpm\cr
&\sum_{\sigma\in\sy_2}\sigma.X_{[k]}^{\pm}(r;\tilde d s)=0&\drmr.\tcalb\cr
(T2^{\pm}):{\phantom{=}}&\sum_{\sigma\in\sy_2}\sigma.[[X_{j,s}^{\pm},X_{i,r_1\pm 1}^{\pm}]_{q^2},X_{i,r_2}^{\pm}]=0&\drmr.\skq\cr
&\sum_{\sigma\in\sy_2}\sigma.M_{[2]}^{\pm}(r;2 s)=0&\drmr.\tcalb\cr
(T3^{\pm}):{\phantom{=}}&\sum_{\sigma\in\sy_2}\sigma.((q^2+1)[[X_{j,s}^{\pm},X_{i,r_1\pm 2}^{\pm}]_{q^3},X_{i,r_2}^{\pm}]_{q^{-1}}+&\drmr.\skq\cr
&+[[X_{j,s}^{\pm},X_{i,r_1\pm 1}^{\pm}]_{q^3},X_{i,r_2\pm 1}^{\pm}]_q)=0\cr
&\sum_{\sigma\in\sy_2}\sigma.M_{[3]}^{\pm}(r;3 s)=0&\drmr.\tcalb\cr
(HXL^{\pm}):{\phantom{=}} &[H_{i,r},X_{j,s}^{\pm}]=\pm b_{ijr}C^{{r\mp|r|\over 2}}X_{j,r+s}^{\pm}\ \ \ (\tilde d_i\leq|r|\leq\tilde d_{ij}) &\drcmp.\tilu\cr
(XXD):{\phantom{=}} &[X_{i,r}^+,X_{j,s}^-]=0 \ \ (i\neq j)&\darl.\xxxtr\cr
(XXE):{\phantom{=}} &[X_{i,r}^+,X_{i,-r}^-]=
{C^{r}k_i-C^{-r}k_i^{-1}\over q_i-q_i^{-1}}&\darl.\xxxtr\cr
(XXH^+):{\phantom{=}} &[X_{i,r}^+,X_{i,s}^-]=
{C^{-s}k_i\tilde H_{i,r+s}^+ \over q_i-q_i^{-1}} \ \ (r+s>0)&\darl.\xxxtr\cr
(XXH^-):{\phantom{=}} & [X_{i,r}^+,X_{i,s}^-]=-
{C^{-r}\tilde H_{i,r+s}^-k_i^{-1}\over q_i-q_i^{-1}} \ \ (r+s<0)&\darl.\xxxtr
\end{align*}
{\bf\underbar{$\C(q)$-algebras}}: 
\begin{align*}
{\Cal U}_q
=& {\roman{Drinfeld\ and\ Jimbo\ quantum\ algebra}}&\drjq.\QNTALG \cr
{\Cal U}_q^{DJ}= &  {\Cal U}_q(\Gamma)& \drjq.\QNTALG \cr
{\Cal U}_q^{fin} = &  {\Cal U}_q(\Gamma_0)&  \drjq.\QNTALG  \cr
{\Cal U}_q^{Dr}   = &  ({\Cal G}|{\Cal {DR}})&\drrl.\drcpjdef \cr
{\Cal U}_q^{Dr}   =&(G|DR)&\drrl.\drdef ,\drrl.\greq  \cr
\tilde{\Cal U}_q^{Dr}  =&(G|ZX^{\pm},CUK,CK,KX^{\pm},XX,HXL^{\pm}) &  \drcmp.\tilu \cr
\bar{\Cal U}_q^{Dr} =&(\bar G|ZX^{\pm},CUK,CK) & \drcmp.\tilu  \cr
{\Cal F}_q^{\pm} =&(G^{\pm}|ZX^{\pm})=(G^{\prime,\pm}) &  \drcmp.\tilpm,\drcmp.\rtilpm
\end{align*}
{\bf\underbar{Elements in the $\C(q)$-algebras}}:
\begin{align*}
\tilde H_{i,\pm r}^{\pm}\ (H_{i,r}):{\phantom{:}}&\sum_{r\in\Z}\tilde H_{i,\pm r}^{\pm}u^r=exp\Big(\pm(q_i-q_i^{-1})\sum_{r>0}H_{i,\pm r}u^r\Big)&\drrl.\drdef\ (\drcmp.\bthi)
\cr
\tilde H_{i,\pm r}^{\pm}=&\begin{cases} 
(q_i-q_i^{-1})k_i^{- 1}[X_{i, r}^{+},X_{i,0}^{-}]&{\roman{ if}}\ r,\pm r>0\cr
(q_i-q_i^{-1})[X_{i,-r}^{-},X_{i,0}^{+}]k_i&{\roman{ if}}\ r>0,\pm r<0\cr
1&{\roman{ if}}\ r=0\cr 0&{\roman{ if}}\ r<0,\end{cases} &\drcmp.\bthi\cr
X_{i,j;l;a} (r;s)= &  \sum_{u=0}^l(-1)^u{l\brack u}_{q_i^a}X_{i,r_1}^{\pm}\cdot...\cdot X_{i,r_u}^{\pm}X_{j,s}^{\pm}X_{i,r_{u+1}}^{\pm}\cdot...\cdot X_{i,r_l}^{\pm}&\drmr.\scmp\cr
M_{i,j;l;a} (r;s)  = &\begin{cases}  \!X_{j,s}^{\pm}&{\roman{if}}\ 
l=0\cr
\![M_{i,j;l-1;a}^{\pm}(r_1,...,r_{l-1};s),\!X_{i,r_l}^{\pm}]_{q_i^{-a_{ij}-2a(l-1)}}&{\roman{if}}\ 
l>0\end{cases}  &\drmr.\nnmm  \cr
M_2^{\pm}((i,r),(j,s))   =& [X_{i,r\pm\tilde d_{ij}}^{\pm},X_{j,s}^{\pm}]_{q_i^{a_{ij}}}+[X_{j,s\pm\tilde d_{ij}}^{\pm},X_{i,r}^{\pm}]_{q_j^{a_{ji}}}  &\drmr.\mmnn\cr
M_i^{\pm}(r) =&[X_{i,r_1\pm\tilde d_{i}}^{\pm},X_{i,r_2}^{\pm}]_{q_i^{2}}& \drmr.\mmnn \cr
M_{(2,2)}^{\pm}(r) = &[X_{1,r_1\pm 2}^{\pm},X_{1,r_2}^{\pm}]_{q^{ 2}}-q^{ 4}[X_{1,r_1\pm 1}^{\pm},X_{1,r_2\pm 1}^{\pm}]_{q^{-6}}  & \drmr.\mmnn\cr
M_{(3)}^{\varepsilon,\pm}(r) =& [[X_{1,r_1\pm\varepsilon}^{\pm},X_{1,r_2}^{\pm}]_{q^{2\varepsilon}},X_{1,r_3}^{\pm}]_{q^{4\varepsilon}}& \drmr.\mmnn \cr
X_{[k]}^{\pm}(r;s) = & \sum_{u,v\geq 0\atop u+v=k-1}q^{v-u}
X_{i,j;2;k}^{\pm}(r_1\pm v,r_2\pm u;s) & \drmr.\mmnn\cr
M_{[2]}^{\pm}(r;s) = &M_{i,j;2;1}^{\pm}(r_1\pm 1,r_2;s)   &\drmr.\mmnn\cr
M_{[3]}^{\pm}(r;s)  = &(q^2\!+\!1)M_{i,j;2;2}^{\pm}(r_1\!\pm\! 2,r_2;s)\!+\!M_{i,j;2;1}^{\pm}(r_1\!\pm\! 1,r_2\!\pm\! 1;s)  &\drmr.\mmnn \cr
k_{\alpha} :{\phantom{:}} &k_{m\delta+\sum_{i\in I_0}m_i\alpha_i}=C^m\prod_{i\in I_0}k_i^{m_i}  &\drcmp.\qgr
\end{align*}

{\bf\underbar{Relations and ideals}}:

a) given the relations 
$$S_{\zeta}(r,s)=0\ (\zeta\in{\Cal{Z}},\ r\in\Z^l,\ s\in \Z^{\tilde l})\leqno(R)$$
denote by $(R_{\zeta})$ ($\zeta\in{\Cal{Z}}$) the relations
$$(R_{\zeta})\hskip 3truecm S_{\zeta}(r,s)=0\ (r\in\Z^l,\ s\in \Z^{\tilde l});\hskip 3truecm\eqno\drmr.\aaa$$

b) given the relations 
$$S_{\zeta}^{\pm}(r,s)=0\ (\zeta\in{\Cal{Z}},\ r\in\Z^l,\ s\in \Z^{\tilde l})\leqno(R^{\pm})$$
denote by $(R)$  the relations
$$(R)\hskip 2truecm S_{\zeta^{\prime}}(r,s)=0\ (\zeta^{\prime}\in{\Cal{Z}\times\{\pm\}},\ r\in\Z^l,\ s\in \Z^{\tilde l})\hskip 1 truecm\eqno\drmr.\aaa$$
where $S_{(\zeta,\pm)}=S_{\zeta}^{\pm}$;

c) given relations $(R)$ as in a), denote by ${\Cal{I}}(R)$ the ideal 
$${\Cal{I}}(R)=(S_{\zeta}(r,s)|\zeta\in{\Cal{Z}},\ r\in\Z^l,\ s\in \Z^{\tilde l}),\eqno\drmr.\aaa$$
by ${\Cal{I}}_{cte}(R)$ the ideal 
$${\Cal{I}}_{cte}(R)=(S_{\zeta}(r\1_l,s)|\zeta\in{\Cal{Z}},\ r\in\Z,\ s\in \Z^{\tilde l}),\eqno\dlem.\relid$$
by ${\Cal{I}}_{0}(R)$ the ideal 
$${\Cal{I}}_{0}(R)=(S_{\zeta}(\underline 0)|\zeta\in{\Cal{Z}})\eqno\dlem.\relid$$
and, if $\zeta\in{\Cal {Z}}$, $r\in\Z^l$, $s\in\Z^{\tilde l}$, denote by ${\Cal{I}}_{(r,s)}(R_{\zeta})$ the ideal 
$${\Cal{I}}_{(r,s)}(R_{\zeta})=(S_{\zeta}(r\1_l,s));\eqno\dlem.\relid$$

d) given relations $(R^{\pm})$ as in a), denote by ${\Cal{I}}_*^{\pm}(R)$ the ideals 
$${\Cal{I}}_*^{+}(R)={\Cal{I}}_*(R^{+}),\ \ \ {\Cal{I}}_*^{-}(R)={\Cal{I}}_*(R^{-})\ \ \ (*\in\{\emptyset,cte,0\});\eqno\drmr.\aaa,\dlem.\relid$$

e) given a family of relations $({}^{(h)}\!R)$ as in a) denote by ${\Cal{I}}_*({}^{(1)}\!R,...,{}^{(m)}\!R)$ the ideals 
$${\Cal{I}}_*({}^{(1)}\!R,...,{}^{(m)}\!R)=({\Cal{I}}_*({}^{(1)}\!R),...,{\Cal{I}}_*({}^{(m)}\!R))\ \ \ (*\in\{\emptyset,cte,0\}).\eqno\drmr.\aaa,\dlem.\relid$$
{\bf\underbar{(Anti)homomorphisms}}:
\begin{align*}
\bar\Omega,\tilde\Omega,\Omega:{\phantom{:}} &q\mapsto q^{-1},\ C^{\pm 1}\mapsto C^{\mp 1},\ k_i^{\pm 1}\mapsto k_i^{\mp 1},\ X_{i,r}^{\pm}\mapsto X_{i,-r}^{\mp}&\drcmp.\gendef,\darl.\boto,\dautemb.\omdef\cr
\Theta_{{\Cal{F}}}^+:{\phantom{:}} &q\mapsto q^{-1},\ X_{i,r}^+\mapsto X_{i,-r}^+&\drcmp.\gendef\cr
\Theta_{{\Cal{F}}}^-:{\phantom{:}} &q\mapsto q^{-1},\ X_{i,r}^-\mapsto X_{i,-r}^-&\drcmp.\gendef\cr
\bar\Theta,\tilde\Theta,\Theta:{\phantom{:}} &q\mapsto q^{-1},\ C^{\pm 1}\mapsto C^{\pm 1},\ k_i^{\pm 1}\mapsto k_i^{\mp 1},&\drcmp.\gendef,\darl.\thibd,\dautemb.\xidef\cr
& X_{i,r}^+\mapsto -X_{i,-r}^+k_iC^{-r},\ X_{i,r}^-\mapsto -k_i^{-1}C^{-r}X_{i,-r}^-\cr
\bar t_i,\tilde t_i,t_i:{\phantom{:}} &C^{\pm 1}\mapsto C^{\pm 1},\ k_j^{\pm 1}\mapsto (k_jC^{-\delta_{ij}\tilde d_i})^{\pm 1},\ X_{j,r}^{\pm}\mapsto X_{j,r\mp\delta_{ij}\tilde d_i}^{\pm}&\drcmp.\gendef,\darl.\ttibd,\dautemb.\tidef\cr
\bar\phi_i,\tilde\phi_i,\phi_i:{\phantom{:}} &q\mapsto q_i,\ C^{\pm 1}\mapsto C^{\pm\tilde d_i},\ k^{\pm 1}\mapsto k_i^{\pm 1},\ X_r^{\pm}\mapsto X_{i,\tilde d_i r}^{\pm}&\drcmp.\gendef,\darl.\fibdef\cr
\phi:{\phantom{:}} &K_i^{\pm 1}\mapsto k_i^{\pm 1},\ E_i\mapsto X_{i,0}^+,\ F_i\mapsto X_{i,0}^-\ (i\in I_0)&\serrel.\fimbd\cr
\tilde \psi,\psi:{\phantom{:}} &C^{\pm 1}\mapsto K_{\delta}^{\pm 1},\ k_i^{\pm 1}\mapsto K_i^{\pm 1},&\hmfr.\defpsi,\hmfr.\bdf,\hmfr.\bendef\cr
&X_{\tilde d_i r}^+\mapsto o(i)^rT_{\lambda_i}^{-r}(E_i),\ X_{\tilde d_i r}^-\mapsto o(i)^rT_{\lambda_i}^{r}(F_i)
\end{align*}

\end{document}